\documentclass[12pt, oneside, reqno, a4paper]{amsart}

\usepackage[utf8]{inputenc}
\usepackage[margin=1in]{geometry} 
\usepackage{enumerate} 
\usepackage{fancyhdr} 
\usepackage{multicol} 
\usepackage{url}
\usepackage[numbers]{natbib}
\usepackage{soul}
\usepackage[svgnames]{xcolor}

\usepackage{graphicx} 
\graphicspath{{figures/}} 
\usepackage{booktabs} 
\usepackage[font=footnotesize,labelfont=bf]{caption} 
\usepackage{amsfonts, amssymb}
\usepackage{amsmath, amsthm} 
\usepackage{wrapfig} 
\usepackage{float}
\usepackage{hyperref}

\usepackage[utf8]{inputenc}
\usepackage[english]{babel}
\usepackage{latexsym}
\usepackage{slantsc}
\usepackage{calligra}
\usepackage{yfonts}
\usepackage{times}
\usepackage{amsbsy,version}
\usepackage{bm,bbm}
\usepackage{mathrsfs}
\usepackage{color}
\usepackage{graphicx}
\usepackage{lineno}
\usepackage{psfrag}
\usepackage[subnum]{cases}
\usepackage{varioref}
\usepackage{wasysym,textcomp}
\usepackage{accents}
\usepackage{relsize}
\usepackage{xfrac}

\usepackage[T3,T1]{fontenc}
    \DeclareSymbolFont{tipa}{T3}{cmr}{m}{n}
    \DeclareMathAccent{\invbreve}{\mathalpha}{tipa}{16}

    \makeatletter
    \newcommand{\snovek}[1]{%
    \scriptsize
    \raisebox{0pt}{\mbox{$
      {\mathop{\scriptstyle #1}\limits^{\vbox to -1.5\ex@{\kern-\tw@\ex@
       \hbox{\tiny\mbox{$\nrightarrow$}}\vss}}}
       $}}
       }
    \makeatother
    
    \makeatletter
    \newcommand{\novek}[1]{%
    \raisebox{0pt}{\mbox{$
      {\mathop{#1}\limits^{\vbox to -1.5\ex@{\kern-\tw@\ex@
       \hbox{\tiny $\nrightarrow$}\vss}}}
       $}}
       }
    \makeatother
    
        \makeatletter
    \newcommand{\bnovek}[1]{%
    \raisebox{0pt}{\mbox{$
      {\mathop{\bar{#1}}\limits^{\vbox to -1.3\ex@{\kern-\tw@\ex@
       \hbox{\tiny $\nrightarrow$}\vss}}}
       $}}
       }
    \makeatother

    \makeatletter
    \newcommand{\bsim}[1]{%
    \raisebox{0pt}{\mbox{$
      {\mathop{\bar{#1}}\limits^{\vbox to -2.0\ex@{\kern-\tw@\ex@
       \hbox{\normalsize $\tilde{}$}\vss}}}
       $}}
       }
    \makeatother
    
        \makeatletter
    \newcommand{\bsimeq}[1]{%
    \raisebox{0pt}{\mbox{$
      {\mathop{#1}\limits^{\vbox to -1.0\ex@{\kern-\tw@\ex@
       \hbox{\tiny $\simeq$}\vss}}}
       $}}
       }
    \makeatother
    
    \makeatletter

    \newcommand{\sbsim}[1]{%
    \scriptsize
    \raisebox{0pt}{\mbox{$
      {\mathop{\bar{#1}}\limits^{\vbox to -1.5\ex@{\kern-\tw@\ex@
       \hbox{\small $\tilde{}$}\vss}}}
       $}}
       }
    \makeatother

\newcommand{\inner}[2]{\left( {#1} {\,,\,} {#2}\right)}

\newcommand{\vek}[1]{\mbox{$\bm{#1}$}}
\newcommand{\ten}[1]{\mbox{$\mathbf{#1}$}}

\newcommand{\T}{\mathrm{T}}

\newcommand{\bxi}{\boldsymbol{\xi}}
\newcommand{\bpsi}{\boldsymbol{\psi}}

\newcommand{\bsigma}{\boldsymbol{\sigma}}

\newcommand{\grad}{{\bf grad\,}}

\newcommand{\bepsilon}{\boldsymbol{\epsilon}}

\newcommand{\bdiv}{{\bf div\,}}

\newcommand{\curl}{{\bf curl\,}}

\renewcommand{\div}[1]{\mathrm{div} {#1} \mathrm{ }}

\def\eolqed{\hspace{\stretch1}\ensuremath\square}
\def\sym#1{\mathrm{sym\{}{#1}\mathrm{\}}}

\def\at0lim{\mbox{$\scriptstyle\partial_{\text{\tiny N}}\mathscr{ B}$} }
\def\ax0lim{\mbox{$\scriptstyle\partial_{\text{\tiny D}}\mathscr{ B}$} }
\def\a0lim{\mbox{$\scriptstyle\partial\mathscr{ B}$} }

\def\3half{{\scriptstyle\frac{3}{2}}}
\def\2third{{\textstyle\dfrac{2}{3}}}
\def\2thrd{{\scriptstyle\frac{2}{3}}}
\def\4thrd{{\scriptstyle\dfrac{4}{3}}}
\def\4third{{\textstyle\dfrac{4}{3}}}

\def\ax0lim{\mbox{$\scriptstyle\partial_{\text{\tiny D}}\mathscr{ B}$} }
\numberwithin{equation}{section}

\numberwithin{assumption}{section}

\numberwithin{notation}{section}

\newtheorem{problem}{Problem}
\numberwithin{problem}{section}

\numberwithin{proposition}{section}

\newtheorem{lemma}{Lemma}
\numberwithin{lemma}{section}
\newtheorem{theorem}{Theorem}

\newtheorem{definition}{Definition}
\numberwithin{definition}{section}

\newtheorem{remark}{Remark}
\numberwithin{remark}{section}

\usecounter{equation}
\definecolor{lightgrey}{rgb}{0.75,0.75,0.75}

\makeatletter
\newcommand{\specialnumber}[1]{%
  \def\tagform@##1{\maketag@@@{(\ignorespaces##1\unskip\@@italiccorr#1)}}%
}
\newcommand{\specialeqref}[2]{\begingroup
  \def\tagform@##1{\maketag@@@{(\ignorespaces##1\unskip\@@italiccorr#2)}}%
  \eqref{#1}\endgroup}
\makeatother

\graphicspath{ {figures} }
\DeclareGraphicsExtensions{.pdf,.png,.jpg}

\title[On pressure robustness \& independent determination in incompressible elasticity]{On pressure robustness and independent determination\\ of displacement and pressure in incompressible linear elasticity}

\author[A.~Zdunek]{Adam Zdunek}
\address[Adam Zdunek]{HB BerRit, Solhemsbackarna 73, SE-163 56 Spånga, Sweden}
\email{adam.zdunek@berrit.se}
\author[M.~Neunteufel]{Michael Neunteufel}
\address[Michael Neunteufel]{Institute for Analysis and Scientific Computing, TU Wien,\\
	Wiedner Hauptstr. 8-10, A-1040 Wien, Austria}
\email{michael.neunteufel@tuwien.ac.at}
\author[W.~Rachowicz]{Waldemar Rachowicz}
\address[Waldemar Rachowicz]{Institute of Computer Science,
	Cracow University of Technology, Pl 31-155 Cracow, Poland}
\email{wrachowicz@pk.edu.pl}

\begin{document}
\maketitle
\begin{abstract}
\normalsize
	We investigate the possibility to determine the divergence-free displacement $\vek{u}$ \emph{independently} from the pressure reaction $p$ for a class of boundary value problems in incompressible linear elasticity. If not possible, we investigate if it is possible to determine it \emph{pressure robustly}, i.e. pollution free from the pressure reaction.
 
		For convex domains there is but one variational boundary value problem among the investigated that allows the independent determination. It is the one with essential no-penetration conditions combined with homogeneous tangential traction conditions.
				
		Further, in most but not all investigated cases, the weakly divergence-free displacement can be computed pressure robustly provided the total body force is decomposed into its direct sum of divergence- and rotation-free components using a Helmholtz decomposition. The elasticity problem is solved using these components as separate right-hand sides. The total solution is obtained using the superposition principle. 
		
		We employ a $(\vek{u},p)$ higher-order finite element formulation with discontinuous pressure elements. It is \emph{inf-sup} stable for polynomial degree $p\ge 2$ but not pressure robust by itself. 
		We propose a three step procedure to solve the elasticity problem preceded by the Helmholtz decomposition of the total body force. The extra cost for the three-step procedure is essentially the cost for the Helmholtz decomposition of the assembled total body force, and the small cost of solving the elasticity problem with one extra right-hand side. The results are corroborated by theoretical derivations as well as numerical results.\\
\vspace*{0.25cm}

\end{abstract}

\section{Introduction}
	\label{sec:Introduction}
	The pressure reaction $p$ and the divergence-free displacement $\vek{u}$ are the two basic independent variables in incompressible linear elasticity. If we exchange the displacement with the velocity we may consider the flow of incompressible Newtonian fluid. Many results in incompressible linear elasticity are simply adopted from the theory for Stokes equations \cite{ChorinMarsden1992}.
	A prominent example is the so called Ladyzhenskaya–Babu\v{s}ka–Brezzi (LBB) \cite{Lad69,Bab73,Bre74} theory which establishes the conditions for the well posedness for the weak formulation of a Dirichlet boundary value problem in terms of the basic variables $\vek{u}$ and $p$ in incompressible linear elasticity, in our case. The \emph{inf-sup} condition of the LBB theory tells if an assumed mixed approximation for the pair $(\vek{u},p)$ results in a stable formulation or not. 
	For our current purposes, it is of importance that the LBB theory predicts that $\vek{u}$ and $p$ can be computed \emph{independently} in two successive steps, provided that the so-called weakly divergence-free test and trial functions for the displacement are available.
	The two equations of the abstract mixed formulation, in the two unknowns $(\vek{u},p)$, decouples then. 
	The LBB theory states conditions for each of the decoupled equations to be solved uniquely for $\vek{u}$ and $p$,
	respectively. The unique solution of $p$ requires the fulfillment of the \emph{inf-sup} condition.
	
	This elegant approach becomes formal in absence of weakly divergence-free test and trial functions, yet we register the predicted decoupling of $\vek{u}$ and $p$. In practice we are relegated to search for the weakly divergence-free displacement $\vek{u}$ in the whole displacement space using the Lagrange multiplier method. The LBB theory describes in detail how $\vek{u}$ and the pressure reaction $p$ enforcing the weak form of the divergence constraint $\div{\,\vek{u}}=0$ are determined.
	Nevertheless, in the realisation of the Lagrange multiplier method the members of the pair $(\vek{u},p)$ are determined together. 
	
	Recently \cite{Linke2014,Linke2016,JohnSIAM2017} it was found that many \emph{inf-sup}  stable mixed  $(\vek{u},p)$, Lagrange multiplier based formulations for Stokes equations in use, fail to determine the member $\vek{u}$ \emph{pressure robustly}, i.e., without pollution from the member $p$, in the sense of \cite{LedererSIAM2017}, clarified by the estimate \eqref{eq:estimate_displacement_error} herein\footnote{Pressure robustness has a physical relevance in fluid flow when Reynolds number $\mathrm{Re}\rightarrow\infty$, or equivalently the viscosity $\mu\rightarrow 0$. However, neither the Stokes equations nor the elasticity equations include the limit $\mu=0$, $\mu$ being the shear modulus in elasticity. Further, as opposed to fluid mechanics there are no phenomena or other equations governing the situation $\mu=0$ in solid mechanics. The loss of pressure robustness in incompressible elasticity simply follows by analogy with the Stokes equations. Its actual significance depends on the appearance of large contrasts between the divergence- and rotation-free body forces in solid mechanics. Answering this question is outside the scope of this investigation.},
	\begin{align}
	\|\vek{u}-\vek{u}_h\|_{\bm{H}^1(\Omega)}\leq C\Big( \inf\limits_{\bm{v}_h\in V_h}\|\vek{u}-\vek{v}_h\|_{\bm{H}^{1}(\Omega)} + \mu^{-1}\inf\limits_{q_h\in Q_h}\|p-q_h\|_{L^2(\Omega)} \Big).
	\tag{\ref{eq:estimate_displacement_error}}
	\end{align}
	Evidently, this estimate tells that the error in the displacement is polluted by the error in the pressure, scaled by the reciprocal of the shear modulus.
	
	In \cite{LedererSIAM2017} pressure robustness is restored by modifying the right-hand side, i.e., the external work, using a reconstruction operator, element by element. We use a related but different approach relying on 
	a Helmholtz decomposition of the involved global body force.
	Further, we mention the class of exactly divergence-free methods like H(div)-conforming hybrid discontinuous Galerkin (HDG) \cite{LS15} or mass conserving mixed stress (MCS) \cite{GLS19} used mainly for Stokes and Navier-Stokes problems, which in addition are pressure robust. Therein the velocity field $\vek{u}$ is discretized by H(div)-conforming Raviart-Thomas or Brezzi-Douglas-Marini elements \cite{RT77,BDM85}. An application of the H(div)-conforming HDG method in linear elasticity has been presented in \cite{FLLS21}.
	
	We will theoretically and by numerical experiments investigate if our \emph{inf-sup} stable higher-order 
	but other wise standard native mixed $(\vek{u},p)$ finite element formulation in incompressible linear elasticity is pressure robust in the sense of \cite{LedererSIAM2017} or not. We use hexahedral discontinuous pressure elements $Q_{k}/P^{disc}_{k-1}$, \mbox{$k\ge 2$} of the mapped type \cite[Section 3.6.4]{John2016}.
	
	The boundary value setting, and the background is presented in Section~\ref{sec:Goal and problem class}. The corresponding variational formulations used are introduced and a deep connection to a curl-curl formulation is presented arising e.g. for Maxwell's equation in Section~\ref{sec:The corresponding class of variational boundary value problems}.
	We also resolve the pollution problem applying Helmholtz decomposition explicitly, i.e. in a different fashion than in \cite{LedererSIAM2017}. 
	
	To that end, a first indication how to resolve the problem was provided by Leonhard Euler \cite[(1757)]{CFT60} who observed that the pressure reaction in an incompressible material is equivalent to a rotation-free\footnote{Throughout the article, the synonyms conservative and solenoidal may be used instead of rotation-free and divergence-free, respectively.} body force acting in the domain $\Omega$, for homogeneous essential conditions, $\vek{u}=\vek{0}$ on the boundary $\partial\Omega$.
	In essence his observation says that the solution to a pure Dirichlet boundary value problem with a rotation-free body force is in the form $(\vek{u},p)=(\vek{0},\phi)$ where $\phi$ is the scalar potential defining the rotation-free body force. That is,  for a rotation-free body force, the pressure is determined independently of the displacement. In particular, one of the side goals of this investigation is to see if Euler's observation 
	can be extended to other boundary conditions than clamped.

	Since we are dealing with linear problems \emph{the superposition principle} applies. To that end, \emph{the Helmholtz decomposition} \cite{Stokes1849,Helmholtz1858,GiraultRaviart86} of the body force into the direct sum of its divergence- and 
	rotation-free parts comes handy into play. It is introduced in Section~\ref{sec:Helmholtz decomposition}. The lucky strike would of course be that a divergence-free  body force
	yields the complementary solution $(\vek{u},p)=(\vek{u},0)$ for $\vek{u}=\vek{0}$ on $\partial\Omega$, which would 
	realise pressure robustness by construction.
	The Helmholtz decomposition together with the superposition principle gives the by us called three step procedure,
	introduced in Section~\ref{sec:A three step procedure relying on a Helmholtz decomposition of the loading}.
	
	We will theoretically and numerically show that the divergence-free body force in general yields a non-trivial solution $(\vek{u},p)$ where the pressure component $p\neq 0$ for $\vek{u}=\vek{0}$ on $\partial\Omega$. Further, 
	we will theoretically and numerically show that the computation with separated divergence- and rotation-free 
	right-hand sides allows a pressure robust determination of the total
        solution $(\vek{u},p)$, in the sense of \cite{LedererSIAM2017}.
        This is possible if so-called mixed-mixed no-penetration and
        homogeneous tangential traction boundary conditions
        or pure essential boundary
        conditions are prescribed on the whole boundary.
        Also the combination of no-slip and homogeneous normal traction
        boundary conditions or homogeneous pure natural boundary conditions is valid.
	
	We investigate a possible non-trivial pressure reaction  by boundary conditions for a divergence-free body
	force in Section~\ref{sec:Possible pressure reaction by boundary conditions for divergence-free body force}.
	Our investigation clarifies that the type of the boundary conditions and also the regularity of the geometry of the domain decides whether the pressure reaction $p$ is trivial or not.
	We show  that essential no-slip and natural homogeneous normal traction boundary conditions, as well as, essential no-penetration and natural homogeneous tangential traction boundary conditions yield the solution in the form
	$(\vek{u},p)=(\vek{u},0)$ for a divergence-free forcing, provided the domain is convex, respectively. This will be achieved by exploiting that for regular enough domains the elasticity and curl-curl problem yield equivalent solutions.
	On the other hand we demonstrate that these boundary conditions yield in general
	a non-trivial singular pressure solution for a divergence-free forcing in a domain with re-entrant corners.
	Non-homogeneous as well as homogeneous essential boundary conditions and homogeneous natural boundary conditions will in general induce 	a non-trivial pressure reaction, $p\neq 0$.
	
	We further investigate a possible non-trivial displacement by boundary conditions for a rotation-free body
	force, Section~\ref{sec:Possible displacement by boundary conditions for a rotation-free body force}.
	We show that Euler's observation $(\vek{u},p)=(\vek{0},\phi)$ directly can be extended to the case of no-penetration and zero tangential traction boundary conditions. For non-homogeneous essential, no-slip and zero normal traction boundary conditions, as well as for homogeneous natural boundary conditions we obtain solutions with a non-trivial displacement $\vek{u}\neq \vek{0}$.
	Finally, we show that the special case when the scalar potential $\phi$ vanishes on the boundary extends 
	Euler's observation to the no-slip and homogeneous natural boundary conditions as well.
	
	Numerical experiments are performed in Section~\ref{sec:Numerical experiments}. Subsection~\ref{sec:The finite element constructs} presents the finite element constructs used.
	Sections~\ref{sec:the cube problems} to \ref{sec:the cube with no-penetration and zero tangential traction boundary conditions} describe performed experiments on a cube domain. Subsection \ref{sec:Helmholtz decomposition of the total body force} presents the test using the Helmholtz decomposition to split a total body force into its divergence and rotation free parts respectively. Section \ref{sec:the clamped cube}
	presents the linear elastic solutions for the clamped cube. Section~\ref{sec:the cube with no-penetration and zero tangential traction boundary conditions} presents the corresponding results for the cube with no-penetration and zero tangential traction boundary conditions. Finally, in Section~\ref{sec:The three-dimensional L-shaped domain} we show that the pressure is non-trivial overall, but concentrated to the re-entrant corners,  on an three-dimensional L-shaped domain subject to no-penetration and zero tangential traction boundary conditions and a divergence-free body force.
	
	The investigation is concluded by a summary and conclusions in Section~\ref{sec:Summary and conclusions}. Supporting material is placed in \ref{sec:proof elasticity curlcurl connection}, \ref{sec:stability curlcurl problem}, and \ref{sec:The assumed body force on the cube}.
\section{The goal of the investigation and the considered problem class}
	\label{sec:Goal and problem class}
	In the following we consider a simply connected bounded Lipschitz domain $\Omega\subset\mathbbm{R}^{3}$
    with a boundary $\Gamma=\partial\Omega$ which consists of planar walls and each with an oriented unit normal $\vek{n}(\vek{x})$. The domain $\Omega$ is assumed to be filled by an homogeneous incompressible isotropic linear elastic material with constant shear modulus $\mu > 0$.
	
	The ultimate goal of this work would be to determine the divergence-free displacement vector field 
	$\vek{u}=\vek{u}(\vek{x})$ \emph{independently} of the scalar pressure reaction field  $p=p(\vek{x})$.
	Another similar goal is to determine 
	$\vek{u}$ \emph{pressure robustly},
	i.e. free from pollution of the scalar pressure reaction field  $p$, in the sense of \cite{LedererSIAM2017}, clarified by the estimate \eqref{eq:estimate_displacement_error} herein.
		Therefore, we will consider pure Dirichlet, no-penetration with zero tangential traction , no-slip with zero normal traction, or natural Neumann boundary conditions, respectively.
	Leonhard Euler observed (1757)\cite{CFT60} that the pressure reaction in an incompressible material is equivalent to a conservative body force $\vek{f}(\vek{x})=\mathrm{grad\,}\phi(\vek{x})$, in the pure Dirichlet case.
	A tempting hypothesis is, on the other hand, that a divergence-free body force $\vek{f}(\vek{x})=\curl\vek{A}(\vek{x})$, such that  $\vek{n}\cdot\vek{f}=0$ on $\Gamma$,
	does not cause any pressure reaction in the pure Dirichlet problem. It turns out to be wrong. 
	We further found that, it is not possible to compute the divergence-free displacement field $\vek{u}$ independently from $p$, i.e. in a decoupled fashion, in general. We investigate if it can be computed pressure robustly \cite{LedererSIAM2017}.

	The basic class of boundary value problems considered in isotropic incompressible linear elasticity is;
	\begin{problem}[Isotropic incompressible linear elasticity problem]
		\label{prob:Pure Dirichlet isotropic incompressible linear elasticity problem}
		\begin{equation}
		\left.
		\begin{array}{lclcl}
		\textit{Find} \quad\mbox{$(\vek{u},p)$}: &&&& \\[1ex]
		\bsigma_{\textrm{a}} &=& 2\mu\,\ten{e}(\vek{u})  &\qquad\text{in} &\Omega, \\[1ex]
		\ten{e}(\vek{u}) &:=&\sym{\grad\vek{u}}-\sfrac{1}{3}(\div{\vek{u}})\ten{1}&\qquad\text{in} &\Omega, \\[1ex] 
		-\bdiv{\bsigma_{\textrm{a}}} + \textrm{grad}\,p &=& \vek{f} &\qquad\text{in} &\Omega, \\[1ex]
		\div{\,\vek{u}}&=&0&\qquad\text{in} &\Omega, \\[1ex]
		+\mathrm{b.c.} &&  &\qquad\text{on} &\Gamma,
		\end{array}\;\right\}
		\specialnumber{a-e}
		\label{eq:problem class}
		\end{equation}
	\end{problem}
	where the boundary conditions will be specified in Section~\ref{sec:Boundary conditions investigated}.
	
	Here the classic divergence-free/traceless  strain $\ten{e}(\vek{u})\in\sym{\mathbbm{R}^{3}\otimes\mathbbm{R}^{3}}$, such that $\ten{e}(\vek{u}):\ten{1}=0$, is introduced,
	wherein $\ten{1}$ is the identity mapping in $\mathbbm{R}^{3}\otimes\mathbbm{R}^{3}$.
	The traceless stresses \specialeqref{eq:problem class}{a} are therefore symmetric \mbox{$\bsigma_{\textrm{a}}=\bsigma_{\textrm{a}}^{\T}$}.
	In passing we note that $\sym{\grad\vek{u}}=\bepsilon(\vek{u})\in\sym{\mathbbm{R}^{3}\otimes\mathbbm{R}^{3}}$ is the classic small strain tensor. 
	
	Here it is understood that the traceless stress $\bsigma_{\textrm{a}}$ in  \eqref{eq:problem class}
	is to be eliminated between \specialeqref{eq:problem class}{a} and \specialeqref{eq:problem class}{c}.
	The resulting formulation requires \mbox{$\vek{u}\in \vek{C}^{2}(\Omega)$} and 
	\mbox{$p\in C^{1}(\Omega)$} in the strong sense. 
	
	\begin{remark}[Euler's observation (1757).]
		\label{rem: Eulers observation}
		Referring to Euler's observation \cite{CFT60}, we note that a gradient of a scalar field may be added to the left and right-hand sides of the equilibrium equation \specialeqref{eq:problem class}{c} in Problem~\ref{prob:Pure Dirichlet isotropic incompressible linear elasticity problem} in combination with clamped boundary conditions without changing the divergence-free displacement $\vek{u}$ in $\Omega$.
		\eolqed  
	\end{remark}
	
	\noindent In the following we use the standard functional spaces;
	\begin{equation}
	\label{eq:basic spaces}
	\begin{array}{lcl}
	{L}^{2}(\Omega) &:=& \left\{q:\Omega\to \mathbb{R}\,\vert\; \displaystyle{\int_{\Omega} \vert q(\vek{x})\vert^2\,d\Omega} < \infty  \right\}, \\[2ex]
	\vek{H}^{1}(\Omega) &:=& \left\{\vek{u}\in\vek{L}^{2}(\Omega)\,\vert\;
	\grad\vek{u}\in\ten{L}^{2}(\Omega)\right\}, \\[2ex]
	\vek{H}(\textrm{curl},\Omega) &:=& \left\{\vek{u}\in\vek{L}^{2}(\Omega)\,\vert\;
	\curl\vek{u}\in\vek{L}^{2}(\Omega)\right\}, \\[2ex]
	\vek{H}(\textrm{div},\Omega) &:=& \left\{\vek{u}\in\vek{L}^{2}(\Omega)\,\vert\;
	\div{\,\vek{u}}\in{L}^{2}(\Omega)\right\}, 
	\end{array}
	\end{equation}

	where $\ten{L}^{2}(\Omega)=[L^{2}(\Omega)]^{3\times 3}$ and $\vek{L}^{2}(\Omega)=[L^{2}(\Omega)]^{3}$.
	Further, we denote 
	scalar fields $q\in{L}_{0}^{2}(\Omega)$ such that $q\in{L}^{2}(\Omega)$ with zero mean value,
	\begin{equation}
	\displaystyle{\int_{\Omega}q(\vek{x})\,d\Omega}=0.
	\label{eq:mean value}
	\end{equation}
	Likewise, $\vek{H}^{1}(\Omega)=[H^{1}(\Omega)]^{3}$, and with $\|\cdot\|_{L^2(\Omega)}$, $\|\cdot\|_{H^1(\Omega)}$, and $\|\cdot\|_{\bm{H}(\mathrm{curl},\Omega)}$ we denote the usual Sobolev norms. We henceforth denote vector fields  $\vek{u}\in\vek{H}^{1}_{0}(\Omega)$ such that \mbox{$\vek{u}\in\vek{H}^{1}(\Omega)$} with $\vek{u}=\vek{0}$ for $\vek{x}\in\Gamma$. 
	Similarly, vector fields $\vek{u}\in\vek{H}_{0}(\textrm{curl},\Omega)$ satisfy the tangential essential boundary condition,
	\mbox{$\vek{n}(\vek{x})\times\vek{u}(\vek{x})=\vek{0}$} for $\vek{x}\in\Gamma$, and vector fields $\vek{u}(\vek{x})\in\vek{H}_{0}(\textrm{div},\Omega)$ satisfy the normal essential boundary condition,
	\mbox{$\vek{n}(\vek{x})\cdot\vek{u}(\vek{x})={0}$} for $\vek{x}\in\Gamma$. 
	\begin{remark}
		Strictly speaking, the trace operators for $\vek{H}^1$, $\vek{H}({\rm curl})$ and $\vek{H}({\rm div})$ should be denoted by $\mathrm{tr\,}u:\vek{H}^1(\Omega)\to \vek{H}^{1/2}(\Gamma)$, $\mathrm{tr}_{\bm{t}}\vek{u}:\vek{H}({\rm curl},\Omega)\to \vek{H}^{-1/2}(\Gamma)$, and $\mathrm{tr}_{\bm{n}}\vek{u}:\vek{H}({\rm div},\Omega)\to H^{-1/2}(\Gamma)$, respectively. For continuous functions the trace operators translate to point, normal and tangential evaluation, respectively. Nevertheless, for simpler notation, we will not use the more abstract trace operators indicating Dirichlet data.\eolqed
	\end{remark}

	Depending on properties of the domain $\Omega$ the solution of the Poisson problem may be more regular than $H^1(\Omega)$. 
	\begin{definition}[s-regularity of Poisson's problem]
		\label{def:Poissons problem}
		The scalar Poisson problem reading:
		\begin{equation}
		\left.
		\begin{array}{lclcl}
		\textit{For given } f\in L^2(\Omega) &\textit{ find }& \quad&\mbox{${u\in H^1(\Omega)}$}: & \\[1ex]
		-\div{\,\mathrm{grad\,}} u &=& {f} &\qquad\text{in} &\Omega, \\[1ex]
		{u} &=&  {0} &\qquad\text{on} &\Gamma,
		\end{array}\;\right\}
		\specialnumber{a,b}
		\label{eq:Poissons problem}
		\end{equation}
		is called \emph{s-regular} with $s\in(0,1]$ if, 
		\[
		\|u\|_{H^{1+s}(\Omega)} \le {c}\,\|f\|_{L^{2}(\Omega)},\qquad \text{ where } c=c(\Omega,s).
		\]
		
		For domains $\Omega$ without re-entrant corners, i.e., $\Omega$ is convex,
		the Poisson-problem is $s$-regular with $s=1$. If not otherwise stated it will henceforth be assumed that $\Omega$ does not have re-entrant corners.
		
		In general, if \mbox{$\vek{u}\in\vek{H}(\textrm{curl},\Omega)\cap\vek{H}(\textrm{div},\Omega)$} it does not follow that $\vek{u}$ is also in $\vek{H}^1(\Omega)$, see e.g. \cite[Lemma 3.56]{Mon03}.
		 If the Poisson problem is $s$-regular, however, the following holds as proved e.g. in \cite[Theorem 3.50]{Mon03}.
	\end{definition}

	\begin{lemma}
		\label{lem:regularity hcurl hdiv h1}
		Let $\vek{u}\in \vek{H}_{0}({\rm curl},\Omega)\cap\vek{H}({\rm div},\Omega)$ or 
		$\vek{H}({\rm curl},\Omega)\cap\vek{H}_{0}({\rm div},\Omega)$. If the Poisson problem is $s$-regular there holds with $\delta:=\max\{1/2,s\}$ and a constant $c>0$ independent of $\vek{u}$
		\begin{align}
		\|\vek{u}\|_{{\bm H}^{\delta}(\Omega)}\le c\,\left(\|\curl{\vek{u}}\|_{{\bm L}^2(\Omega)}+\|\div{\,\vek{u}}\|_{{L}^2(\Omega)}\right).\label{eq:s_reg_hcurl_hdiv}
		\end{align}
		\eolqed
	\end{lemma}
	This result will be used in Sections~\ref{sec:Possible pressure reaction by boundary conditions for divergence-free body force} and \ref{sec:Possible displacement by boundary conditions for a rotation-free body force} with $s=1$. Our focus now turns to an equivalent weak formulation of \eqref{eq:problem class}.
	
	\section{The corresponding class of variational boundary value problems}
	\label{sec:The corresponding class of variational boundary value problems}
	The corresponding class of variational boundary value problems equivalent to \eqref{eq:problem class} is obtained eliminating the traceless stress $\bsigma_{\rm a}=\bsigma_{a}(\vek{u})$ between \specialeqref{eq:problem class}{a} and \specialeqref{eq:problem class}{c} and testing the resulting equation
	with the vector field \mbox{$\vek{v}\in \vek{V}(\Omega)$}. Further, equation
	\specialeqref{eq:problem class}{d} is tested with the scalar field \mbox{$q\in Q(\Omega)$}. Here 
	$Q(\Omega)\subset L^2(\Omega)$ is given by \eqref{eq:basic spaces}, and \mbox{$\vek{V}(\Omega)\subset\vek{H}^1(\Omega)$} which may include boundary conditions, will be specified later.
	Using Gauss theorem and integrating by parts yields,
	\begin{equation}
	\displaystyle{\int_{\Gamma}\vek{v}\cdot\vek{s}\,d\Gamma}=\displaystyle{\int_{\Omega}\div{(\bsigma^{\T}\vek{v})}\,d\Omega} = \displaystyle{\int_{\Omega}\grad\vek{v}:    \bsigma\,d\Omega} + 
	\displaystyle{\int_{\Omega}\vek{v}\cdot\bdiv{\bsigma}\,d\Omega}. 
	\label{eq:Gauss theorem + integrating by parts}
	\specialnumber{a}
	\end{equation}
	Here, and in the sequel we denote the traction vector given by the Cauchy stress principle,
	\begin{equation}
	\vek{s}:=\bsigma\vek{n}=\bsigma_{\rm a}\vek{n}-p\vek{n} = \vek{s}_{\rm a}-p\vek{n} \quad\text{where}\quad \vek{s}_{\rm a}:=2\mu\,\ten{e}(\vek{u})\vek{n}.
	\tag{\ref{eq:Gauss theorem + integrating by parts}b}
	\end{equation}
	Substituting $\bsigma=-p\ten{1}$ in \specialeqref{eq:Gauss theorem + integrating by parts}{a} and simplifying yields the following useful identity,
	\begin{equation}
	\displaystyle{\int_{\Gamma}\vek{v}\cdot p\vek{n}\,d\Gamma} = \displaystyle{\int_{\Omega}\div{(p\,\vek{v})}\,d\Omega} = \displaystyle{\int_{\Omega}p\,\div{\vek{v}}\,d\Omega} + 
	\displaystyle{\int_{\Omega}\vek{v}\cdot\mathrm{grad\,}p\,d\Omega}.
	\tag{\ref{eq:Gauss theorem + integrating by parts}c}
	\end{equation}
	Using \specialeqref{eq:Gauss theorem + integrating by parts}{a}, \specialeqref{eq:Gauss theorem + integrating by parts}{c} the stress symmetry $\bsigma^{\T}=\bsigma$ and $\ten{1}:\ten{e}=0$ our basic mixed $(\vek{u},p)$ variational boundary value formulation of isotropic incompressible linear elasticity becomes: 
	\begin{problem}[Isotropic incompressible linear elasticity]
		\label{prob:LE}
		\begin{equation}
		\left.
		\begin{array}{lcl}
		\textit{For given }\mu > 0, \vek{f}\in \vek{L}^2(\Omega) \textit{ and } \vek{s}\in \vek{L}^2(\Gamma)\textit{ find }&&\mbox{$\vek{u}\in \vek{V}(\Omega)$ and $p\in Q(\Omega)$}: \\[2ex]
		\phantom{-} \displaystyle{\int_{\Omega}\ten{e}(\vek{v}): [2\mu\,\ten{e}(\vek{u})]\,d\Omega }-\displaystyle{\int_{\Omega}p\,\div{\vek{v}}\,d\Omega} &=& \displaystyle{\int_{\Omega}\vek{v}\cdot\vek{f} \,d\Omega}
		+\int_{\Gamma}\vek{v}\cdot\vek{s}\,d\Gamma,
		 \\[2ex]
		-\displaystyle{\int_{\Omega}q\,\div{\vek{u}}\,d\Omega}&=&0,
		\end{array}\;\right\}
		\label{eq:weak-formulation}
		\tag{LE}
		\end{equation}
		for any $\vek{v}\in \vek{V}(\Omega)$ and $q\in Q(\Omega)$.

	The behavior with respect to decoupling of displacement and pressure as well as pressure robustness for different boundary conditions will be investigated in Sections~\ref{sec:Possible pressure reaction by boundary conditions for divergence-free body force} and \ref{sec:Possible displacement by boundary conditions for a rotation-free body force}. 
	The test and trial set $\vek{V}(\Omega)$  is a subset of $\vek{H}^1(\Omega)$ and $Q(\Omega)$ is a subset of $L^{2}(\Omega)$ respectively. They will be specified in Section~\ref{sec:Boundary conditions investigated}. Among the investigated cases we will have those
	with a non-trivial traction vector $\vek{s}\neq\vek{0}$ on $\Gamma$. We therefore keep the boundary term in the right-hand side of \eqref{eq:weak-formulation}\textsubscript{1}. \eolqed
	\end{problem}
	\begin{remark}[Realisation of Problem~\ref{prob:LE}]
	\label{rem:Perturbed Lagrangian}
	Equations \eqref{eq:weak-formulation} correspond to a pure saddle point problem. In the numerical experiments we use a so-called Perturbed Lagrangian or Projection Penalty approach \cite{Bercovier78}. It implies that
	that the term $\int_\Omega{q}{\varepsilon\,p}\,d\Omega$ is added in the left-hand side of equation \eqref{eq:weak-formulation}\textsubscript{2}, where volumetric compliance $\varepsilon = 1/\kappa \ge 0$ is the reciprocal of the bulk modulus $0 < \kappa < \infty$. Near incompressibility is usually described by a ratio $\kappa/\mu \ge 10^{3}$. 
	In finite precision arithmetic the Projection Penalty approach allows about half the actual precision to be used to suppress volumetric changes.
	 \eolqed
	\end{remark}

	\begin{remark}[Stokes equations]
		Equations \eqref{eq:weak-formulation} correspond to the Stokes equations in fluid mechanics
		replacing $\ten{e}(\vek{u})$ by the traceless part of the strain rate tensor, i.e., $\dot{\ten{e}}(\vek{u})$,  cf. \cite[Eq.(1.3.4)]{ChorinMarsden1992}.\eolqed
	\end{remark}
	
	If  $\curl\vek{u}\times \vek{n}=\vek{0}$ or the normal component, tangential components, or both for $\vek{u}$ and $\vek{v}$ are zero, the left-hand side of \eqref{eq:weak-formulation}\textsubscript{1} involving only the displacement fields can be reformulated using the curl operator.
	\begin{lemma}[The elasticity - curl-curl connection]
		\label{th:The elasticity - curl-curl connection}
		Let $\vek{u},\vek{v}\in\vek{H}^1(\Omega)$ and $\vek{u}$ be divergence-free, $\div{\vek{u}}\equiv0$. 
		\begin{enumerate}
			\item If $\vek{n}\times\curl\vek{u}=\vek{0}$ on $\Gamma$ there holds
		\begin{align}
		\int_{\Omega} \mu\,\curl \vek{u}\cdot \curl\vek{v}\,d\Omega = \int_{\Omega}2\mu\,\ten{e}(\vek{u}):\ten{e}(\vek{v})\,d\Omega+\int_{\Gamma}2{\mu}\,\ten{e}(\vek{u})\vek{n}\cdot \vek{v}\,d\Gamma,
		\label{eq:curlcurl_epseps_bnd}
		\end{align}
		where $\mu\in\mathbbm{R}^{+}$.
		\item If $\vek{u},\vek{v}\in\vek{H}_{0}({\rm div},\Omega)$ or $\vek{u},\vek{v}\in\vek{H}_{0}({\rm curl},\Omega)$, i.e., $\vek{n}\cdot\vek{u}=\vek{n}\cdot\vek{v}=0$ or $\vek{n}\times \vek{u}=\vek{n}\times \vek{v}=\vek{0}$ on $\Gamma$, respectively, there holds,
		\begin{align}
		\int_{\Omega} \mu\,\curl \vek{u}\cdot \curl\vek{v}\,d\Omega = \int_{\Omega}2\mu\,\ten{e}(\vek{u}):\ten{e}(\vek{v})\,d\Omega
		\quad\text{where}\quad \mu \in \mathbbm{R}^{+}.
		\label{eq:curlcurl_epseps}
		\end{align}
		\end{enumerate} 
		The proof is moved to \ref{sec:proof elasticity curlcurl connection} for better readability.
	\end{lemma}
	
	\begin{remark}[Shifting regularity between displacement and pressure field]
		\label{rem:regularity of p and q}
		We will see that a well posed elasticity - curl-curl formulation can be achieved increasing
		the regularity requirement on the trial and test functions $p$ and $q$, from $L^{2}(\Omega)$ to $H^{1}(\Omega)$ while  lowering the regularity of $\vek{u}$ and $\vek{v}$ from $\vek{H}^1(\Omega)$ to $\vek{H}(\mathrm{curl},\Omega)$.
		\eolqed
	\end{remark}
	
	Motivated by Lemma~\ref{th:The elasticity - curl-curl connection} we consider the following curl-curl problem arising e.g. for Maxwell's equations. We note, that due to the lowered regularity of the displacement field only its tangential or normal component can be controlled, not both simultaneously. 

	\begin{problem}[Curl-curl elasticity problem]
		\label{prob:AUX-LE}
		Using the elasticity -- curl-curl connection, see the identity \eqref{eq:curlcurl_epseps} in the Lemma~\ref{th:The elasticity - curl-curl connection} and the integration by parts identity \specialeqref{eq:Gauss theorem + integrating by parts}{c}
		we obtain the following auxiliary elasticity formulation:
		\begin{equation}
		\left.
		\begin{array}{lcl}
		\textit{For given } \vek{f}\in\vek{L}^2(\Omega) \textit{ and } \vek{s}_{\rm a}\in \vek{L}^2(\Gamma) \textit{ find}&& \mbox{$\vek{u}\in \vek{W}(\Omega)$ and $\pi\in P(\Omega)$}:  \\[2ex]
		\displaystyle{\int_{\Omega}\!\mu\,\curl\vek{u}\cdot\curl\vek{v}\,d\Omega}  + \displaystyle{\int_{\Omega}\vek{v}\cdot\mathrm{grad\,}\pi\,d\Omega} &=& \displaystyle{\int_{\Omega}\vek{v}\cdot\vek{f} \,d\Omega} \\[2.5ex]
		&+& \displaystyle{\int_{\Gamma} \vek{n}\times\vek{v}\cdot\vek{s}_{\rm a}\,d\Gamma}
		, \\[3ex]
		\displaystyle{\int_{\Omega}\!\vek{u}\cdot\mathrm{grad\,}\varpi\,d\Omega}&=&0,
		\end{array}\right\}
		\label{eq:auxiliary curl-curl elasticity formulation} 
		\tag{AUX-LE}
		\end{equation}
		for any $\vek{v}\in\vek{W}(\Omega)$ and $\varpi\in P(\Omega)$,
		where the partial traction vector $\vek{s}_{\rm a}$ is defined in \specialeqref{eq:Gauss theorem + integrating by parts}{b}. Here, $\vek{W}(\Omega) \subset \vek{H}({\rm curl},\Omega)$ while
		$P(\Omega) \subset H^{1}(\Omega)$. They will be specified for each given problem. For completeness a proof that the problem is well-posed will be given in \ref{sec:stability curlcurl problem}.
		\eolqed
	\end{problem}
	
	\begin{remark}[Realisation of Problem~\ref{prob:AUX-LE})]
		\label{rem:Realisation of Problem AUX-LE}
		The key to the construction of the mixed method \eqref{eq:auxiliary curl-curl elasticity formulation} is the following compatibility condition between the spaces $\vek{W}(\Omega)$ and $P(\Omega)$ inherited from electromagnetics \cite{WRLDEM00} used in the proof for well-posedness of Problem~\ref{prob:AUX-LE},
		\[
		\curl\vek{u}=\vek{0} \quad\Leftrightarrow\quad\exists\,\varpi\in P(\Omega):\; \mathrm{grad\,}\varpi =\vek{u}\in\vek{W}(\Omega) .
		\]
		It is henceforth understood that the $(\vek{u},\pi)\in\vek{W}(\Omega)\times P(\Omega)$  mixed method \eqref{eq:auxiliary curl-curl elasticity formulation} is realised using $\vek{H}({\rm curl},\Omega)$ conforming vector test and trial functions for the displacement and scalar $H^{1}(\Omega)$ conforming test and trial functions for the pressure reaction including no constant functions except zero, respectively.
		In this realisation, no-penetration boundary conditions $\vek{n}\cdot\vek{u}=0$ on $\Gamma$ will be imposed 
		in the sense of a natural condition by equation \eqref{eq:auxiliary curl-curl elasticity formulation}\textsubscript{2}. The exact definition of test and trial spaces $\vek{W}(\Omega)$ and $P(\Omega)$ are given in Section~\ref{sec:Boundary conditions investigated} where boundary conditions are discussed.
		\eolqed
	\end{remark} 
	
	We now show the major impact of the curl-curl elasticity Problem~\ref{prob:AUX-LE} namely that the pressure reaction $\pi$ is trivial for a divergence-free body force and homogeneous Neumann data.
	\begin{theorem}
		\label{th:0}
		The Lagrange multiplier $\pi \in P(\Omega)$ in problem \eqref{eq:auxiliary curl-curl elasticity formulation} is trivial for a divergence-free body force, \mbox{$\vek{f}=\curl\vek{A}$} that is parallel to the boundary $\Gamma$, $\vek{n}\cdot\curl\vek{A}=0$ on $\Gamma$ and if the boundary term on the right-hand side of 
		\eqref{eq:auxiliary curl-curl elasticity formulation}\textsubscript{1} vanishes.
		\begin{proof}\!\!:
			Use the test function $\vek{v}=\mathrm{grad\,}\pi\in \vek{W}(\Omega)$ with $\pi\in P(\Omega)$ the solution of \eqref{eq:auxiliary curl-curl elasticity formulation}. The first integral on the left-hand side of 
			\eqref{eq:auxiliary curl-curl elasticity formulation}\textsubscript{1} vanishes since $\curl\mathrm{grad\,}\pi\equiv\vek{0}$,  while the volume integral on the right-hand side of 
			\eqref{eq:auxiliary curl-curl elasticity formulation}\textsubscript{1} vanishes by $L^{2}(\Omega)$-orthogonality, see \eqref{eq:L2-orthogonality version 1}. Equation \eqref{eq:auxiliary curl-curl elasticity formulation}\textsubscript{1} simplifies to,
			\begin{equation*}
			\int_{\Omega} \mathrm{grad\,}\pi\cdot \mathrm{grad\,}\pi\,d\Omega = \|\mathrm{grad\,}\,\pi\|^{2}_{L^{2}}= 0,
			\end{equation*}
			and the claim follows.
		\end{proof}
	\end{theorem}
	
	The opposite, that a rotation-free 
	body force combined with homogeneous Neumann data will only introduce pressure, is also partly true as stated in the next theorem.
	\begin{theorem}[Trivial displacement solutions]
		A rotation-free body force, $\vek{f}=\mathrm{grad\,}\,\phi\in L^{2}(\Omega)$, $\phi\in H^1(\Omega)\backslash\mathbb{R}$, together with $\vek{s}_{\rm a}=\vek{0}$ is equivalent to the gradient of the Lagrange multiplier,  $\;\mathrm{grad\,}\,\pi$ in problem \eqref{eq:auxiliary curl-curl elasticity formulation}  for pure no-penetration with zero tangential traction boundary conditions for all $\phi\in H^1(\Omega)\backslash\mathbb{R}$, or no-slip with zero normal-traction if $\phi\in H^1_0(\Omega)$. Moreover, the displacement field $\vek{u}\in \vek{W}(\Omega)$ is trivial (see Section~\ref{sec:Boundary conditions investigated} for the boundary conditions).
		\begin{proof}\!\!:
			Using  the  rotation-free test function $\vek{v}=\mathrm{grad\,}(\pi-\phi)\in\vek{W}(\Omega)$ 
			in \eqref{eq:auxiliary curl-curl elasticity formulation}\textsubscript{1} we obtain,
			\[
			\int_{\Omega}\vert\mathrm{grad\,}\,\pi-\mathrm{grad\,}\,\phi\vert^2\, d\Omega = 0.
			\]
			Thus, the claim follows, $\mathrm{grad\,}\,\pi = \mathrm{grad\,}\,\phi$,  and $(\vek{u},\pi)=(\vek{0},\phi)$ solves Problem~\ref{prob:AUX-LE}.
		\end{proof}
	\end{theorem}
	
	\subsection{Discussion of boundary conditions}
	\label{sec:Boundary conditions investigated}
	
	Apart from the prescribed displacement boundary conditions we consider three other types. The complete list reads,
	\begin{enumerate}
		\item No-penetration with zero tangential traction boundary conditions on $\Gamma$: 
		\begin{equation}
		u_n:=\vek{n}\cdot\vek{u}=0 \quad\text{and}\quad  \vek{n}\times\vek{s}=\vek{n}\times\vek{s}_{\rm a}=\vek{0}. 
		\label{eq:boundary conditions}
		\specialnumber{a}
		\end{equation}
		For Problem~\ref{prob:LE} the sets of test and trial functions are: $\vek{U}(\Omega)=\vek{V}(\Omega)=$ \\$\left\{\vek{w}\in \vek{H}^{1}(\Omega)\,\vert\; w_n=0\text{ on }\Gamma\right\}$ and $Q(\Omega)=L^{2}_{0}(\Omega)$. For Problem~\ref{prob:AUX-LE} they formally are: $\vek{W}(\Omega)=\vek{H}(\textrm{curl},\Omega)$ and $P(\Omega)=H^{1}(\Omega)\backslash\mathbbm{R}$, where the zero normal component arises naturally, see Remark~\ref{rem:Realisation of Problem AUX-LE}. 
		\item No-slip and zero normal traction on $\Gamma$:  
		\begin{equation}
		\vek{n}\times\vek{u}=\vek{0} \quad \text{and} \quad \vek{n}\cdot\vek{s}=0.
		\tag{\ref{eq:boundary conditions}b}
		\end{equation}
		For Problem~\ref{prob:LE} the sets of test and trial functions are: $\vek{U}(\Omega)=\vek{V}(\Omega)=$\\ $\left\{\vek{w}\in \vek{H}^{1}(\Omega)\,\vert\; \vek{n}\times\vek{w}=\vek{0}\text{ on }\Gamma\right\}$ and $Q(\Omega)=L^{2}(\Omega)$. For Problem~\ref{prob:AUX-LE} they formally are: $\vek{W}(\Omega)=\vek{H}_{0}(\textrm{curl},\Omega)$ and $P(\Omega)=H^{1}_0(\Omega)$. 
		\item Clamped or prescribed displacement on $\Gamma$:
		\begin{equation}
		\vek{u}=\bar{\vek{u}} \quad\text{such that}\quad \int_{\Gamma} \bar{\vek{u}}\cdot\vek{n}\,d\Gamma = 0,
		\tag{\ref{eq:boundary conditions}c}
		\end{equation}
		where $\bar{\vek{u}}=\bar{\vek{u}}(\vek{x})$ the prescribed displacement field is consistent with the divergence-free constraint. 
		
		For Problem~\ref{prob:LE} test and trial spaces are: $\vek{U}(\Omega)=\vek{V}(\Omega)=\vek{H}_{0}^{1}(\Omega)$ 
		for clamped conditions, $\vek{u}=\vek{0}$, while $\vek{U}(\Omega)=\{\tilde{\vek{u}}\} + \vek{H}_{0}^{1}(\Omega)$ and $\vek{V}(\Omega)=\vek{H}_{0}^{1}(\Omega)$ for prescribed conditions with $\tilde{\vek{u}}\in \vek{H}^1(\Omega)$ any extension of $\bar{\vek{u}}$ into the domain $\Omega$, $\tilde{\vek{u}}=\bar{\vek{u}}\neq\vek{0}$ on $\Gamma$, respectively.\\
		For both these cases, $Q(\Omega)=L^{2}_{0}(\Omega)$.\\
		We note, that the consistency condition on $\bar{\vek{u}}$ is needed to enable a point-wise divergence-free solution $\vek{u}$. If violated, the problem is still solvable, however, the solution only fulfills the property $\int_{\Omega}\div{\vek{u}}\,q\,d\Omega=0$ for all $q\in L^2_0(\Omega)$, but not point-wise, $\div{\vek{u}}\not\equiv 0$.
		\item Homogeneous Neumann or zero traction on $\Gamma$: 
		\begin{equation}
		\vek{s}=\vek{0}.
				\tag{\ref{eq:boundary conditions}d}
		\end{equation}
		For Problem~\ref{prob:LE} test and trial spaces are: $\vek{U}(\Omega)=\vek{V}(\Omega)=\vek{H}^1(\Omega)$, $Q(\Omega)=L^2(\Omega)$. In order to remove any rigid displacement we add mean value side conditions on the displacement trial- and test functions, 
		$\int_{\Omega} \vek{w}\, d\Omega =  \vek{0}$ and  $\int_{\Omega} \vek{r}\times\vek{w}\, d\Omega =  \vek{0}$  for $\vek{w}\in\{\vek{u},\vek{v}\}$, respectively. Here, $\vek{r}$ is a fixed location vector. 
	\end{enumerate}
	
	\begin{remark}
		Note, that both, no-penetration with zero tangential traction and no-slip with zero normal traction boundary condition, are of the so-called \emph{mixed-mixed} type. They consist of a homogeneous essential part for the normal (tangential) and a homogeneous natural tangential (normal) part, respectively. 
		\eolqed
	\end{remark}
	
		As already mentioned above, the no-penetration boundary condition $\vek{n}\cdot\vek{u}=0$ is enforced  in a weak sense in\eqref{eq:auxiliary curl-curl elasticity formulation}\textsubscript{2}. Furthermore it turns out that $\vek{u}$ is not only in $\vek{H}(\mathrm{curl},\Omega)$ but also in $\vek{H}(\mathrm{div},\Omega)$
		\begin{lemma}
			\label{lem:curlcurl solution hdiv}
			The solution $\vek{u}\in\vek{W}(\Omega)$ of \eqref{eq:auxiliary curl-curl elasticity formulation} is in $\vek{H}(\mathrm{div},\Omega)$ and additionally in $\vek{H}_0(\mathrm{div},\Omega)$ if $P(\Omega)=H^1(\Omega)\backslash\mathbb{R}$.
			\begin{proof}\!\!:
				The solution $\vek{u}$ has zero weak divergence as stated by \eqref{eq:auxiliary curl-curl elasticity formulation}\textsubscript{2} and thus $\vek{u}\in\vek{H}(\mathrm{div},\Omega)$. As already noted in Remark~\ref{rem:Realisation of Problem AUX-LE} for $P(\Omega)=H^1(\Omega)\backslash\mathbb{R}$ there directly follows that $\vek{u}\in\vek{H}_0(\mathrm{div},\Omega)$.\eolqed
			\end{proof}
	\end{lemma}
	
	Using the auxiliary elasticity formulation in Problem~\ref{prob:AUX-LE} we will show that the no-penetration with zero tangential traction boundary conditions, for \emph{a convex domain} admit a decoupled determination of the divergence-free displacement $\vek{u}$ and the pressure reaction $p$ of the elasticity problem, while other boundary conditions prevent such a decoupling.
 	We will also use the three-dimensional L-shaped domain to numerically illustrate that singular pressure fields are induced also in this case.

	The interest is now turned to  an $L^{2}(\Omega)$-orthogonal decomposition of the body force
	\mbox{$\vek{f}\in\vek{L}^2(\Omega)$} into divergence-free and rotation-free components. The ideal situation would be if this orthogonal decomposition 
	implied that the displacement $\vek{u}$ and pressure $p$ could be computed independently. That is,
	if one of the components generated pressure $p$ only and the other was responsible for the the divergence-free displacement $\vek{u}$. We will show that the boundary conditions will not allow such an 
	independent determination of $\vek{u}$ and $p$, in general.

	\section{The Helmholtz decomposition}
	\label{sec:Helmholtz decomposition}
	In the following we assume that $\vek{f}\in \vek{L}^2(\Omega)$ is a given field of body forces per unit volume in $\Omega$.
	Following Helmholtz \cite{Stokes1849,Helmholtz1858} we assume the following decomposition of the body force into its rotation-free and divergence-free parts as,
	\begin{equation}
	\vek{f} = \mathrm{grad\,}\phi + \curl\vek{A}, \quad\text{with}\quad 
	\div{\vek{A}} = 0.
	\label{eq:Helmholtz-decomposition}
	\end{equation}
	Here, $\phi\in H^1(\Omega)\backslash\mathbb{R}$ and $\vek{A}\in \vek{H}_0(\mathrm{curl}, \Omega)$. 
From the latter we obtain $\vek{n}\cdot\curl\vek{A}=0$ on $\Gamma$, i.e. the divergence-free part of the body force is parallel to the boundary \eqref{eq:Helmholtz-decomposition}\textsubscript{3}. Finally, \eqref{eq:Helmholtz-decomposition}\textsubscript{2} is recognised as the Coulomb gauge\footnote{The left-hand side of \eqref{eq:Helmholtz-decomposition}\textsubscript{1} corresponds to three independent components while the right-hand side introduces four independent quantities. The Coulomb gauge \eqref{eq:Helmholtz-decomposition}\textsubscript{2} balances the system.}.
	
	\begin{theorem}
		\label{th:1}
		Divergence free and rotation-free body forces are $L^{2}(\Omega)$-orthogonal provided the divergence-free body force is parallel to the boundary.
		\begin{proof}\!:\;
			Let any divergence-free body force be given as the curl of a vector potential 
			$\vek{A}\in \vek{H}_0(\mathrm{curl}, \Omega)$, i.e. as $\vek{f}_{1}=\curl\vek{A}$ such that $\vek{n}\cdot\curl\vek{A}=0$ on $\Gamma$. Further, let any rotation-free body force be given as the gradient of a scalar potential 
			$\phi\in H^1(\Omega)\backslash\mathbb{R}$, i.e. as $\vek{f}_{2}=\mathrm{grad\,}\phi$. 
			The claim follows by partial integration of the following $L^{2}(\Omega)$-inner product,
			\begin{equation}
			\int_{\Omega} \curl\vek{A} \cdot \mathrm{grad\,}\phi \,d\Omega =
			-\int_{\Omega} \phi\, \underbrace{\div{\,\curl\vek{A}}}_{\equiv 0}  \,d\Omega + \int_{\Gamma}\phi\,
			\underbrace{\vek{n}\cdot\curl\vek{A}}_{=0}\,d\Gamma = 0.
			\label{eq:L2-orthogonality version 1}
			\end{equation}
		\end{proof}
	\end{theorem}
	Alternatively, the rotation-free part of the body force can be considered perpendicular to the boundary, i.e.\linebreak[4]  
	\mbox{$\vek{n}\times\mathrm{grad\,}\phi=\vek{0}$} on $\Gamma$. There are two classic and frequently used methods for the Helmholtz decomposition \eqref{eq:Helmholtz-decomposition}, one for the vector potential $\vek{A}$
	and the other for the scalar potential $\phi$, see for example  \cite{Bathiaetal2013}.
	Girault and Raviart \cite{GiraultRaviart86} propose a vector Laplacian formulation for the determination of the vector potential.

	Our preferred Helmholtz decomposition relies on the setting \eqref{eq:Helmholtz-decomposition}
	and the $L^{2}(\Omega)$-orthogonality provided by Theorem~\ref{th:1}.
	It corresponds to the following strong formulation for the determination of the vector potential $\vek{A}$; 
	\begin{equation}
	\left.
	\begin{array}{lclcl}
	\textit{For given } \vek{f}\in\vek{L}^2(\Omega) &\textit{ find }& \quad\mbox{$\vek{A}\in \vek{H}(\mathrm{curl},\Omega)$}: && \\[1ex]
	\curl\curl\vek{A} &=& \curl\vek{f} &\qquad\text{in} &\Omega, \\[1ex]
	\div{\,\vek{A}}&=&0&\qquad\text{in} &\Omega, \\[1ex]
	\vek{n}\times\vek{A} &=& \vek{0} &\qquad\text{on} &\Gamma. 
	\end{array}\;\right\}
	\specialnumber{a,b,c}
	\label{eq:strong Helmholtz formulation A}
	\end{equation}
	Note, that for $\vek{f}\in \vek{L}^2(\Omega)$ its curl lies in the topological dual space of $\vek{H}_0(\mathrm{curl},\Omega)$, $\curl\vek{f}\in \vek{H}_0(\mathrm{curl},\Omega)^*=\vek{H}^{-1}(\mathrm{div},\Omega)=\left\{\vek{w}\in \vek{H}^{-1}(\Omega)\,\vert\; \div{\,\vek{w}}\in H^{-1}(\Omega)\right\}$, $H^{-1}(\Omega)=H^1_0(\Omega)^*$ \cite[Lemma 2.1]{Pechstein11}, such that \specialeqref{eq:strong Helmholtz formulation A}{a} has to be understood in the sense of functionals (or as operators).
	Using a Lagrange multiplier $\pi\in H_{0}^{1}(\Omega)$ to enforce the Coulomb gauge 
	\specialeqref{eq:strong Helmholtz formulation A}{b}
	the corresponding $(\vek{A},\pi)$-mixed weak formulation becomes:
	\begin{problem}[Helmholtz decomposition]
		\label{prob:HD}
		\begin{equation}
		\left.
		\begin{array}{lcl}
		\textit{Find} \quad\mbox{$\vek{A}\in \vek{H}_{0}(\textrm{curl},\Omega)$ and $\pi\in H_{0}^{1}(\Omega)$}: && \\[2ex]
		\displaystyle{\int_{\Omega}\curl\vek{A}\cdot\curl\vek{v}\,d\Omega}  + \displaystyle{\int_{\Omega}\vek{v}\cdot\mathrm{grad\,}\pi\,d\Omega} &=& \displaystyle{\int_{\Omega}\curl\vek{v}\cdot\vek{f} \,d\Omega},\\[3ex]
		\displaystyle{\int_{\Omega}\vek{A}\cdot\mathrm{grad\,}\varpi\,d\Omega}&=&0,
		\end{array}\;\right\}
		\label{eq:weak Helmholtz formulation A}
		\tag{HD}
		\end{equation}
		for any $\vek{v}\in \vek{H}_{0}(\textrm{curl},\Omega)$ and $\varpi\in H_{0}^{1}(\Omega)$.
		Here we have used the definition of the duality pairing between the functional $\curl \vek{f}\in \vek{H}_0(\mathrm{curl},\Omega)^*$ and $\vek{v}\in\vek{H}_0(\mathrm{curl},\Omega)$ 
		\begin{equation}
		\langle\curl\vek{f}, \vek{v}\rangle_{\bm{H}_0(\mathrm{curl},\Omega)^*} =
			\int_{\Omega} \curl\vek{v}\cdot\vek{f} \,d\Omega, 
		\end{equation}
		where no boundary term appears as $\vek{v}\in\vek{H}_{0}(\textrm{curl},\Omega)$.
		It is noted that the formulation relies on test and trial functions $\varpi,\pi\in H_{0}^{1}(\Omega)$. We further note that the choice
		$\pi\in H_{0}^{1}(\Omega)$ implies a continuous Lagrange multiplier $\pi$ when discretizing with Lagrangian and with \emph{hp}-adaptive finite elements \cite{WRLDEM00}.\eolqed
	\end{problem}
	\begin{theorem}
		\label{th:2}
		The Lagrange multiplier $\pi$ in problem \eqref{eq:weak Helmholtz formulation A} is trivial.
		\begin{proof}\!\!:
			See the proof of Theorem~\ref{th:0}.
		\end{proof}
	\end{theorem}
	Having determined the vector potential $\vek{A}$ we compute its curl and compute the 
	rotation-free part of the body force as, $\textrm{grad\,}\phi=\vek{f}-\curl\vek{A}$. The  
	determination  of the scalar potential $\phi$ is cheaper and does not involve the solution of a saddle point problem.
	The approaches are equal on the continuous level. On the discrete level, however, they differ significantly as $\vek{f}-\textrm{grad\,} \phi_h$ is in general not a divergence-free function. This might lead to problems for not pressure-robust methods, see Remark~\ref{rem:neglection_of_u2}. $\vek{f}-\curl\vek{A}_h$ is also in general not a gradient field. Nevertheless, the approximation error is not critical in the corresponding step of the three step procedure described in the next section. Therefore, we emphasize to determine the vector potential $\vek{A}$ instead of the scalar potential $\phi$.
	
	Having decomposed the body force according to \eqref{eq:Helmholtz-decomposition} the linear form in the elasticity problem~\eqref{eq:weak-formulation} is divided into two $L^{2}(\Omega)$-orthogonal linear forms. These will in general induce two mixed solutions $(\vek{u}_{1},p_{1})$ and $(\vek{u}_{2},p_{2})$, corresponding to the divergence-free
	and rotation-free part of the body force, respectively. These can be computed separately on account of the linear superposition principle. We will show in Sections~\ref{sec:Possible pressure reaction by boundary conditions for divergence-free body force} and \ref{sec:Possible displacement by boundary conditions for a rotation-free body force} that $p_{1}\neq 0$ as well as $\vek{u}_2\neq 0$, in general, clearly depending on the boundary conditions.
	
	\begin{remark}[Alternative Helmholtz decomposition.]
	Instead of solving \eqref{eq:weak Helmholtz formulation A} for the vector potential we could have solved a simpler scalar Poisson problem for the 
	potential $\phi$. The approaches are seemingly equivalent on the continuous level. But computing the divergence free 
	part of the body force as $\curl\vek{A}=\vek{f}-\textrm{grad\,}\phi$ the gauge $\div{\,\vek{A}}=0$ of the vector potential is lost. The gauge condition is needed for the left-hand and right-hand side in $\vek{f}=\curl\vek{A}+\textrm{grad\,}\phi$ to be balanced.
	 \eolqed
	\end{remark}

	\section{A three step procedure relying on a Helmholtz decomposition of the loading}
	\label{sec:A three step procedure relying on a Helmholtz decomposition of the loading}
	
	The initial step (Step 0) in the three step procedure corresponds to decomposing the body force into its rotation-free and divergence-free parts according to \eqref{eq:Helmholtz-decomposition}. It is described in Section~\ref{sec:Helmholtz decomposition}. Assuming that the decomposition \eqref{eq:Helmholtz-decomposition} is available the finalising two steps are implied by the superposition principle valid for linear elasticity.
	
	Step 1 and Step 2 will therefore correspond to solving the pure Dirichlet incompressible linear elasticity variational boundary value problem \eqref{eq:weak-formulation} with two different right-hand sides given for Step 1 by:
	$\vek{f}_{1}=\curl\vek{A}$ and for Step 2 by:
	$\vek{f}_{2}=\vek{f}-\curl\vek{A}$, respectively. The solutions are denoted  
	$(\vek{u}_{i},p_{i})$, $i=1,2$, respectively. 
	The three step procedure can be summarised as follows;
	\begin{enumerate}
		\item[]Step 0: given body force $\vek{f}$ solve \eqref{eq:weak Helmholtz formulation A} for vector potential $\vek{A}$. 
		Compute $\curl\vek{A}$ and $\mathrm{grad\,}\phi=\vek{f}-\curl\vek{A}$.
		\item[]Step 1: given the solenoidal part of the body force $\vek{f}_{1}=\curl\vek{A}$ solve \eqref{eq:weak-formulation} for $(\vek{u}_{1},p_{1})$. 
		\item[]Step 2: given the rotation-free part of the body force $\vek{f}_{2}=\vek{f}-\curl\vek{A}$ solve \eqref{eq:weak-formulation} for $(\vek{u}_{2},p_{2})$.
	\end{enumerate}
	On account of linearity, the total solution is simply the sum 
	$(\vek{u},p)=(\vek{u}_{1}+\vek{u}_{2},p_{1}+p_{2})$.
	
	\begin{remark}
		In the case of a prescribed non-zero displacement field $\bar{\vek{u}}$ or non-homogeneous Neumann boundary data $\vek{s}$ they have to be integrated in Step 1 of the three step procedure. As we will see in Sections~\ref{sec:Clamped boundary conditions} and \ref{sec:Homogeneous Neumann boundary conditions} the possibly induced pressure will not be harmful with respect to pressure robustness. Further, more critical, when handled in Step 2 of the procedure the resulting displacement field $\vek{u}_2$ might be overlapped by ``pollution'' induced by the not pressure-robustness treated gradient fields.
		\eolqed
	\end{remark}
	\begin{remark}
		\label{rem:neglection_of_u2}
		If only no-penetration with zero tangential traction or Dirichlet boundary conditions are given, the displacement field $\vek{u}_2$ might contain only non-zero ``noise'' due to lack of pressure robustness, where $\vek{u}_2=\vek{0}$ would be the exact solution. Therefore, the unusable displacement field will be neglected and the final solution takes the form $(\vek{u},p)=(\vek{u}_{1},p_{1}+p_{2})$.
		\eolqed
	\end{remark}
		\section{Possible pressure reaction by boundary conditions for a divergence-free body force}
	\label{sec:Possible pressure reaction by boundary conditions for divergence-free body force}
	In this section we will investigate if a given divergence-free body force $\vek{f}=\curl\vek{A}$, parallel to the boundary $\Gamma$, 
	$\vek{n}\cdot\curl\vek{A}=0$ on a domain $\Omega$ with piece wise flat boundary, may induce a non-trivial pressure reaction, depending on  the boundary condition of the type defined in Section \ref{sec:Boundary conditions investigated}.

	We will see that, for a \emph{convex domain}, for no-penetration with zero tangential traction and no-slip with zero normal-traction boundary conditions no pressure will be induced, whereas Dirichlet and homogeneous Neumann boundary conditions will give a non-trivial pressure.

	In the remaining parts of this section we will investigate each boundary condition separately. The overall idea, however, is the same: First, we will solve the auxiliary elasticity problem \eqref{eq:auxiliary curl-curl elasticity formulation}  using the curl of the vector  potential $\vek{f}=\curl\vek{A}$,  determined by solving problem \eqref{eq:weak Helmholtz formulation A} or given analytically, to construct a divergence-free displacement candidate $\vek{u}$. Then, depending on the specific chosen boundary condition, a second auxiliary problem might be needed to be solved to correct the boundary data of $\vek{u}$. This second step will also answer the question if a non-trivial pressure reaction $p$ is induced.
	
	\subsection{No-penetration with zero tangential traction boundary conditions}
	\label{sec:No-penetration and zero tangential traction boundary conditions}
	We solve the auxiliary problem \eqref{eq:auxiliary curl-curl elasticity formulation} to obtain a candidate 
	$(\vek{u},0)\in \vek{H}({\rm curl},\Omega)\times H^{1}(\Omega)$ fulfilling the boundary conditions \specialeqref{eq:boundary conditions}{a}. It remains to show that $(\vek{u},0)$ also solves elasticity problem \eqref{eq:weak-formulation} with boundary conditions \specialeqref{eq:boundary conditions}{a}. Lemma~\ref{lem:curlcurl solution hdiv} states that the trial function 
	$\vek{u}\in\vek{H}({\rm curl},\Omega)\cap\vek{H}_{0}({\rm div},\Omega)$ and it follows from \eqref{eq:s_reg_hcurl_hdiv} that $\vek{u}\in\vek{H}^1(\Omega)$ if the requirements of Lemma~\ref{lem:regularity hcurl hdiv h1} are fulfilled with $s=1$. By choosing $\vek{v}\in \vek{H}^1(\Omega)$ such that $v_n=0$ on $\Gamma$, test function $\vek{v}$ in \specialeqref{eq:auxiliary curl-curl elasticity formulation}{a} is also evidently in 
	$\vek{H}({\rm curl},\Omega)\cap\vek{H}_{0}({\rm div},\Omega)$, and further using  by Theorem~\ref{th:0}  that the Lagrange multiplier is zero, we obtain,
	\begin{align}
	\label{eq:symmetry curlcurl to elasticity}
	\int_{\Omega}\vek{f}\cdot\vek{v}\,d\Omega&=\int_{\Omega}\mu\,\curl\vek{v}\cdot\curl\vek{u}\,d\Omega\overset{\eqref{eq:curlcurl_epseps}}{=}\int_{\Omega}\ten{e}(\vek{v}):[2\mu\,\ten{e}(\vek{u})]\,d\Omega,
	\end{align}
	concluding the proof, for $\vek{f}$ being divergence-free and parallel to the boundary.
	
	\begin{remark}
		The three-dimensional L-shaped domain problem will show numerically in Section~\ref{sec:The three-dimensional L-shaped domain} that the regularity requirement of Lemma~\ref{lem:regularity hcurl hdiv h1} is a necessary condition. Otherwise a non-zero pressure might be induced by the singularity arising from the geometry. It follows from the idea of taking  $\vek{w}\in \vek{H}_0(\mathrm{curl},\Omega)$ with $\div{\vek{w}}=0$. Thus, there exists a potential $\vek{z}\in \vek{H}^1(\Omega)$ such that $\curl\vek{z}=\vek{w}$ and $\vek{f}=\curl \vek{w}$ is a valid divergence-free body force. The solution of the curl-curl problem is given by $\vek{u}=\nabla\varphi +\vek{z}$, where $\varphi\in H^1(\Omega)\backslash\mathbb{R}$ solves the Poisson problem $-\Delta \varphi=\div{\vek{z}}$ and $\frac{\partial \varphi}{\partial n} = -z_n$, such that $\vek{u}\in \vek{H}_0(\mathrm{div},\Omega)\cap \vek{H}(\mathrm{curl},\Omega)$ with $\div{\vek{u}}\equiv0$, but not in $\vek{H}^1(\Omega)$, as in general, $\varphi$ is not in $H^2(\Omega)$, cf. the s-regularity of the Poisson problem in Definition~\ref{def:Poissons problem}.
		\eolqed
	\end{remark}
	
	\subsection{No-slip and zero normal traction boundary conditions}
	\label{sec:No-slip and zero normal traction boundary conditions}
	The case of no-slip and zero normal traction boundary conditions follows the same idea as in the previous subsection. This time, however, we solve the auxiliary problem \eqref{eq:auxiliary curl-curl elasticity formulation}  to obtain the candidate $(\vek{u},0)\in \vek{H}_{0}({\rm curl},\Omega)\times H^{1}_0(\Omega)$ fulfilling the boundary conditions \specialeqref{eq:boundary conditions}{b}. Again, by Lemma~\ref{lem:curlcurl solution hdiv} the trial function $\vek{u}\in\vek{H}_{0}({\rm curl},\Omega)\cap\vek{H}({\rm div},\Omega)$ and it follows from Lemma~\ref{lem:regularity hcurl hdiv h1}, with $s=1$, that $\vek{u}\in\vek{H}^1(\Omega)$. By choosing $\vek{v}\in \vek{H}^1(\Omega)$ such that $\vek{n}\times\vek{v}=\vek{0}$ on $\Gamma$, 
	test function $\vek{v}$ in \eqref{eq:auxiliary curl-curl elasticity formulation}\textsubscript{1}
	is evidently also in $\vek{H}_{0}({\rm curl},\Omega)\cap\vek{H}({\rm div},\Omega)$,
	and further by Theorem~\ref{th:0} using that the Lagrange multiplier is zero we obtain
	\begin{align*}
	\int_{\Omega}\vek{f}\cdot\vek{v}\,d\Omega&=\int_{\Omega}\mu\,\curl\vek{v}\cdot\curl\vek{u}\,d\Omega\overset{\eqref{eq:curlcurl_epseps}}{=}\int_{\Omega}\ten{e}(\vek{v}):[2\,\mu\ten{e}(\vek{u})]\,d\Omega,
	\end{align*}
	concluding that $(\vek{u},0)$ solves elasticity problem \eqref{eq:weak-formulation} with boundary conditions  \specialeqref{eq:boundary conditions}{b}.

	\begin{remark}
		Also for no-slip with zero traction boundary conditions, the regularity assumption is necessary. For a non-convex domain with non-smooth boundary a non-zero pressure might be induced for the elasticity problem following the idea: Take a potential $\vek{z}\in \vek{H}^2_0(\Omega)$, {where additionally to the trace also all first derivatives vanish at the boundary}, and
		that has a non-vanishing divergence, $\div{\vek{z}}\neq 0$. The force $\vek{f}=\curl\curl \vek{z}$ fulfills the divergence-free assumption. The solution of the curl-curl problem is given by \mbox{$\vek{u}= \nabla \varphi + \vek{z}\in \vek{H}_0(\mathrm{curl},\Omega)$}, where $\varphi\in H^1_0(\Omega)$ solves the Poison problem $-\Delta \varphi = -\div{\vek{z}}$, such that $\div{\vek{u}}\equiv 0$ fulfills the divergence-free condition. Due to the non-convexity of the domain, $\varphi$, is in general, not in $H^2(\Omega)$. Thus $\vek{u}$ is not in $\vek{H}^1(\Omega)$.\eolqed
	\end{remark}

	\subsection{Clamped boundary conditions}
	\label{sec:Clamped boundary conditions}
	In the two previous subsections we observed that homogeneous essential normal and tangential boundary conditions lead to a trivial pressure. One could think that fixing all components of the displacement at the boundary, also results in zero pressure. This, however, is not true in general as we will show now.
	
	Again, we start by solving the  auxiliary problem \eqref{eq:auxiliary curl-curl elasticity formulation} with no-penetration with zero tangential traction boundary conditions, \specialeqref{eq:boundary conditions}{a} obtaining the trial solution $(\vek{u}_{0},0)$. That is, we control the normal component $\vek{u}_{0}\cdot\vek{n}=0$ on $\Gamma$
	while the tangential component $\vek{n}\times\vek{u}_{0}$ is free and in general not zero.
	Therefore, we need to solve the elasticity problem \eqref{eq:weak-formulation} with zero right-hand side $\vek{f}=\vek{0}$ 
	prescribing the computed displacement $\bar{\vek{u}}=\vek{u}_{0}$ on $\Gamma$  to correct the boundary conditions.
	We designate the obtained correction $(\vek{u},p) \rightarrow (\vek{w},r)$.
	The pressure reaction $r$ will be non-zero if the divergence constraint needs to be enforced.
	The difference $\vek{u}_{0}-\vek{w}$ has both zero normal and zero tangential components on $\Gamma$, and is (weakly) divergence-free in $\Omega$. The correction problem \eqref{eq:weak-formulation} can be interpreted as the following extension problem:\\[1ex] Find $\vek{w}\in\vek{H}^1(\Omega)$ with $\vek{w}=\bar{\vek{u}}$ on $\Gamma$ minimizing the constrained energy
	\begin{align}
	\vek{w}=\arg\min\limits_{\substack{\bm{v}\in \bm{H}^1(\Omega)\\\bm{v}=\bar{\bm{u}}\text{ on }\Gamma\\\div{\bm{v}}=0}}\left\{\int_{\Omega}\mu\|\ten{e}(\vek{v})\|^2\,d\Omega\right\}.
	\end{align}

	We claim, that $(\vek{u},p):=(\vek{u}_{0}-\vek{w}, -r)$ solves the original elasticity problem \eqref{eq:weak-formulation} for clamped boundary conditions \specialeqref{eq:boundary conditions}{c}, $\vek{u}=\vek{0}$ on $\Gamma$, i.e., $p\neq 0$, in general.
	
	For $\vek{v}\in \vek{H}^{1}_{0}(\Omega)$ the right-hand side of \eqref{eq:auxiliary curl-curl elasticity formulation},  
	\begin{align*}
	\int_{\Omega} \vek{v}\cdot\vek{f} \,d\Omega &\overset{\eqref{eq:symmetry curlcurl to elasticity}}{=}\int_{\Omega}\left\{ \ten{e}(\vek{v}):\left[2\mu\,\ten{e}(\vek{u}_{0}) - 2\mu\,\ten{e}( \vek{w})\right]+2\mu\,\ten{e}( \vek{v}):\ten{e}(\vek{w})\right\} \,d\Omega\\
	&= \int_{\Omega}\left\{ \ten{e}(\vek{v}):[2\mu\,\ten{e}( \vek{u})] +\ten{e}( \vek{v}): [2\mu\,\ten{e}(\vek{w})]\right\}\,d\Omega\\
	&= \int_{\Omega}\left\{ \ten{e}( \vek{v}):[2\mu\,\ten{e}(\vek{u})] +  r\,\div{\vek{v}}\right\} \,d\Omega,
	\end{align*}
	coincides with the left-hand side of \eqref{eq:weak-formulation}, which proves the conjecture. For the last equality we used that $\int_{\Omega}\ten{e}( \vek{v}):[2\mu\,\ten{e}(\vek{w})]-r\,\div{\vek{v}}\,d\Omega=0$ for all $\vek{v}\in \vek{H}_{0}(\Omega)$ from the correction problem.\\

	One important question with respect to the correction problem \eqref{eq:weak-formulation}, with zero body force $\vek{f}=\vek{0}$ and non-homogeneous Dirichlet condition $\vek{w}=\bar{\vek{u}}$ on $\Gamma$, is how the induced pressure behaves and if it might lead to problems for $\mu\to 0$. The following lemma shows, that the pressure is bounded by $\mu$ times the norm of the extension.\\
	
	\begin{lemma}
		\label{th:5.1}
		Let $(\vek{w},r)$ be the solution of \eqref{eq:weak-formulation}, with zero body force $\vek{f}=\vek{0}$ and non-homogeneous Dirichlet condition $\vek{w}=\bar{\vek{u}}$ on $\Gamma$. Then there holds the estimate with a constant $c>0$ independent of $\mu$, $r$, and $\vek{w}$
		\begin{align}
		\|r\|_{L^2(\Omega)}\leq c \mu \|\vek{w}\|_{\bm{H}^1(\Omega)}.\label{eq:est_pressure_ext_prob}
		\end{align}
		\begin{proof}\!\!:
			Following the proofs of the Stokes LBB condition, see e.g. \cite[Theorem 6.3]{Braess07} there exists a $\vek{v}\in \vek{H}^{1}_{0}(\Omega)$ such that $\div{\vek{v}}=-r$ and $\|\vek{v}\|_{\bm{H}^1}\leq c \|r\|_{L^2}$. Using this $\vek{v}$ as test function of the weak formulation of \eqref{eq:weak-formulation}, with zero body force $\vek{f}=\vek{0}$ and non-homogeneous Dirichlet condition $\vek{w}=\bar{\vek{u}}$ on $\Gamma$,
			\begin{align*}
			\int_{\Omega} \left\{ 2\mu\, \ten{e}(\vek{w}):\ten{e}(\vek{v}) - r\,\div{\vek{v}} \right\}\,d\Omega=0,
			\end{align*} 
			results in,
			\begin{align*}
			\|r\|^2_{L^2(\Omega)} + 2\mu \int_{\Omega}\ten{e}(\vek{w}):\ten{e}(\vek{v})\,d\Omega = 0, \qquad\mu={\rm const.},
			\end{align*}
			and further with the Cauchy-Schwarz inequality,
			\begin{align*}
			\|r\|^2_{L^2(\Omega)} \leq 2\mu \|\vek{w}\|_{\bm{H}^1(\Omega)}\|\vek{v}\|_{\bm{H}^1(\Omega)} \leq 2 c \mu \|\vek{w}\|_{\bm{H}^1(\Omega)}\|r\|_{L^2(\Omega)}.
			\end{align*}
			Dividing through $\|r\|_{L^2(\Omega)}$ yields the desired result.\eolqed
		\end{proof}
	\end{lemma}

	Lemma~\ref{th:5.1} states that the pressure is proportional to $\mu$. Therefore, the elasticity extension problem of boundary data without sources does not induce a problem using a discretization method which is not pressure robust, i.e. in which the displacement error depends on the pressure error,
	\begin{align}
	\|\vek{u}-\vek{u}_h\|_{\bm{H}^1(\Omega)}\leq C\Big( \inf\limits_{\bm{v}_h\in V_h}\|\vek{u}-\vek{v}_h\|_{\bm{H}^{1}(\Omega)} + \mu^{-1}\inf\limits_{q_h\in Q_h}\|p-q_h\|_{L^2(\Omega)} \Big).
	\label{eq:estimate_displacement_error}
	\end{align}
	The reciprocal $\mu^{-1}$ in front of the pressure error gets compensated by the $\mu$ arising on the right-hand side of estimate \eqref{eq:est_pressure_ext_prob}. Numerical examples confirm this observation.\\
	
	As a result we proved that clamped boundary conditions may introduce non-zero pressure, however, without causing problems with respect to pressure robustness. 
	\begin{remark}[Non-homogeneous prescribed displacements]
		\label{rem:Non-homogeneous prescribed displacements}
		A further simple observation is that this procedure can be directly extended to non-homogeneous prescribed displacements  
		$\bar{\vek{u}}$ on $\Gamma$, condition \specialeqref{eq:boundary conditions}{c}, by solving the elasticity problem \eqref{eq:weak-formulation} for $(\vek{u},p)\rightarrow(\vek{w},r)$, with zero body force $\vek{f}=\vek{0}$ and the non-homogeneous Dirichlet condition, 
		$\vek{w}=\vek{u}_{0}-\bar{\vek{u}}$, where $\vek{u}_{0}$ is the displacement
		solution of the auxiliary elasticity problem  \eqref{eq:auxiliary curl-curl elasticity formulation} with no-penetration with zero tangential traction boundary conditions \specialeqref{eq:boundary conditions}{a}. 
		\eolqed
	\end{remark}
	
	\subsection{Homogeneous Neumann boundary conditions}
	\label{sec:Homogeneous Neumann boundary conditions}
	In the case of clamped boundary conditions the solution of the first auxiliary elasticity problem 
	only fulfilled no-penetration with zero tangential traction boundary conditions.
	For homogeneous Neumann boundary conditions \specialeqref{eq:boundary conditions}{d} the situation is exactly the reverse. The solution to the auxiliary elasticity problem with no-penetration with zero tangential traction boundary conditions only  fulfills the tangential zero traction condition. 
	There will be a non-trivial normal traction reaction corresponding to the 
	prohibited normal displacement on $\Gamma$. 
	
	Again, we start by solving the auxiliary elasticity problem  \eqref{eq:auxiliary curl-curl elasticity formulation} with no-penetration and zero tangential traction boundary conditions, \specialeqref{eq:boundary conditions}{a}. Its solution is denoted $(\vek{u}_{0},0)$ as before.
	We need to correct this solution for the traction reaction, $\bsigma\vek{n}\cdot\vek{n}=(2\mu\,\ten{e}(\vek{u}_{0}))\vek{n}\cdot\vek{n}\neq 0$ on
	$\Gamma$. The correction problem will therefore be a non-homogeneous elasticity Neumann problem \eqref{eq:weak-formulation} 
	without body force $\vek{f}=\vek{0}$ in $\Omega$, imposing the traction $\vek{s}=[2\mu\,\ten{e}(\vek{u}_{0})]\vek{n}$ on $\Gamma$. The trial solution is denoted $(\vek{w},r)$. 
	The corrected solution $\vek{u}_{0}-\vek{w}$ will be determined up to a rigid displacement. The non-uniqueness 
	of the homogeneous Neumann problem is removed imposing the following zero mean value constraints for the correction problem,
	\[
	\int_{\Omega}\vek{u}_{0}-\vek{w}\,d\Omega = \vek{0} \quad\text{and}\quad \int_{\Omega}\vek{r}\times(\vek{u}_{0}-\vek{w})\,d\Omega = \vek{0},
	\]
	where, $\vek{r}$ is a location vector.
	We claim that $(\vek{u},p)= (\vek{u}_{0}-\vek{w}, -r)$ solves the elasticity problem \eqref{eq:weak-formulation} with free instead of clamped boundary conditions.
	For $\vek{v}\in \vek{H}^{1}(\Omega)$ we use \eqref{eq:curlcurl_epseps_bnd} of Lemma~\ref{th:The elasticity - curl-curl connection}
	\begin{align*}
	\int_{\Omega}\vek{v}\cdot\vek{f} \,d\Omega &=
	\int_{\Omega} (\ten{e}(\vek{v}):[2\mu\,\ten{e}(\vek{u}_{0}) - 2\mu\,\ten{e}(\vek{w})]+ \ten{e}(\vek{v}):2\mu\,\ten{e}(\vek{w})) \,d\Omega - \int_{\Gamma} (\vek{v}\cdot[2\mu\,\ten{e}(\vek{u}_{0})-0\ten{1}]\vek{n})\,d\Gamma\\
	&= \int_{\Omega} (\ten{e}(\vek{v}):[2\mu\,\ten{e}(\vek{u})])\,d\Omega + \int_{\Omega}( \ten{e}(\vek{v}):[2\mu\,\ten{e}(\vek{w})] )\,d\Omega 
	- \int_{\Gamma} (\vek{v}\cdot[2\mu\,\ten{e}(\vek{u}_{0})-0\ten{1}]\vek{n}) \,d\Gamma\\
	&= \int_{\Omega} (\ten{e}(\vek{v}): [2\mu\,\ten{e}(\vek{u})]+r\,\div{\vek{v}} )\,d\Omega,
	\end{align*}
	which again is the left-hand side of \eqref{eq:weak-formulation}.\\
	
	With the same technique as for clamped boundaries we can prove that the pressure norm is bounded by $\mu$ times the norm of the displacement solution. One has only to replace $\vek{v}\in \vek{H}^{1}_{0}(\Omega)$ by $\vek{v}\in \vek{H}^{1}(\Omega)$ with zero mean value for the test function. Thus, the in general non-trivial pressure will not give problems with respect to pressure robustness.
	
	\section{Possible displacement caused by boundary conditions for a rotation-free body force}
	\label{sec:Possible displacement by boundary conditions for a rotation-free body force}
	In accordance with the previous section we will now investigate in which situation as far as the boundary conditions are concerned a rotation-free body force, on the other hand, induces a non-zero \emph{weakly divergence-free displacement field}.
	
	Throughout this section we assume that the rotation-free body force is given by the gradient of a scalar potential field,
	\begin{align}
	\vek{f} = \mathrm{grad\,}\phi,\qquad \phi\in H^1(\Omega)\backslash\mathbb{R}.
	\label{eq:rotation-free-body-force}
	\end{align}
	
	For no-penetration with zero tangential traction and for clamped boundary conditions, no displacement will be induced. In the case of no-slip with zero normal traction or homogeneous Neumann boundary conditions, however, a nonzero displacement field is expected, except when the potential is prescribed zero at the boundary, i.e., for $\phi\in H^{1}_{0}(\Omega)$. Then, no displacement will appear independently of the prescribed boundary conditions, as long as homogeneous data is given. In the following we again investigate each of the four boundary conditions, introduced in Section~\ref{sec:Boundary conditions investigated} in more detail.

	\subsection{No-penetration with zero tangential traction boundary conditions and clamped boundary conditions}
	\label{sec:No-penetration and clamped boundary conditions}
	Assuming no-penetration with zero tangential traction boundary conditions, \specialeqref{eq:boundary conditions}{a} and the given rotation-free body force \eqref{eq:rotation-free-body-force} we directly claim that the solution of the elasticity problem \eqref{eq:weak-formulation}\textsubscript{1} is given by $(\vek{u},p)=(\vek{0},\phi)$.
	Noting that for the test functions $\vek{v}$ such that  $\vek{n}\cdot\vek{v}=0$ elaborating on the right-hand side 
	of \eqref{eq:weak-formulation}, by partial integration, there holds,
	\begin{align*}
	\int_{\Omega}\vek{f}\cdot\vek{v}\,d\Omega = \int_{\Omega}\mathrm{grad\,}\phi\cdot\vek{v}\,d\Omega = -\int_{\Omega}\phi\,\div{\vek{v}}\,d\Omega + \int_{\Gamma}\phi\,\vek{n}\cdot\vek{v}\,d\Gamma = -\int_{\Omega}\phi\,\div{\vek{v}}\,d\Omega,
	\end{align*}
	which coincides with the left-hand side of \eqref{eq:weak-formulation}\textsubscript{1} for $(\vek{u},p)=(\vek{0},\phi)$.\\

	Also for clamped boundary conditions $\vek{u}=\vek{0}$ we have with the same arguments as before that $(\vek{u},p)=(\vek{0},\phi)$ is the solution of \eqref{eq:weak-formulation}.\\
	
	\begin{remark}
		We emphasize that Euler's observation, Remark~\ref{rem: Eulers observation}, can thus directly be extended to the case of no-penetration with zero tangential traction boundary conditions.
		\eolqed
	\end{remark}
	
	In the case of non-homogeneous prescribed boundary conditions \specialeqref{eq:boundary conditions}{c}, the situation of course changes. Since the trial solution $(\vek{u},p)=(\vek{0},\phi)$ for 
	clamped conditions obviously  does not fulfill the non-homogeneous boundary condition 
	$\vek{u}=\bar{\vek{u}}\neq \vek{0}$ on $\Gamma$, the correction problem  
	described in Remark~\ref{rem:Non-homogeneous prescribed displacements} has to be considered. 
	We claim that $(\vek{u},p)=(\vek{w},\phi+r)$ is the solution of the governing elasticity problem \eqref{eq:weak-formulation}, which is confirmed by the straight forward computation, for $\vek{v}\in\vek{H}^1_0(\Omega)$,
	\begin{equation*}
	\int_{\Omega}\vek{f}\cdot\vek{v}\,d\Omega = -\int_{\Omega} (\phi+r-r)\div{\vek{v}} \,d\Omega = \int_{\Omega}(r\,\div{\vek{v}} -p\,\div{\vek{v}} )\,d\Omega = \int_{\Omega}(\ten{e}(\vek{v}):[2\mu\,\ten{e}(\vek{w})] - p\,\div{\vek{v}})\,d\Omega,
	\end{equation*}
	where we have used the correction problem in the last equality. Note, that the mean value constraint,
	\[
	\int_{\Omega}(\phi+r)\,d\Omega = 0,
	\]
	holds.

	\subsection{No-slip and zero normal traction boundary conditions}
	\label{sec:No-slip boundary conditions rotation-free body force}
	In the case of divergence-free body forces and no-slip boundary conditions, \specialeqref{eq:boundary conditions}{b},  we obtained a trivial pressure, see Theorem~\ref{th:0}. Given the rotation-free body force \eqref{eq:rotation-free-body-force}, we cannot expect $\vek{s}=\vek{0}$ on the boundary, unless the potential explicitly vanishes there, i.e., $\phi\in H^1_0(\Omega)$. Thus, an elasticity corrective problem \eqref{eq:weak-formulation}\footnote{Problem \ref{prob:AUX-LE} cannot be used since we need to impose a non-homogeneous normal surface traction.} has to be solved for $(\vek{w},r)$. In this problem, on $\Gamma$ we impose the mixed mixed boundary conditions $\vek{n}\times\vek{w}=\vek{0}$ with non-homogeneous Neumann data $\vek{s}=\phi\vek{n}$, wherein $\phi$ is the scalar potential defining the rotation-free body force.
	As a result the displacement $\vek{w}$ will be non-zero in general.
	We claim that the solution of \eqref{eq:weak-formulation} with boundary conditions \specialeqref{eq:boundary conditions}{b} has the form $(\vek{u},p)=(\vek{w},\phi+r)$, which is confirmed by the following straight forward computation,
	\begin{equation*}
	\int_{\Omega}\vek{f}\cdot\vek{v}\,d\Omega = -\int_{\Omega}\phi\,\div{\vek{v}}\,d\Omega + 
	\int_{\Gamma}\phi\,\vek{n}\cdot\vek{v}\,d\Gamma \overset{\eqref{eq:weak-formulation}}{=} \int_{\Omega}( \ten{e}(\vek{v}):[2\mu\,\ten{e}(\vek{w})] - (\phi+r)\div{\vek{v}})\,d\Omega
	\end{equation*}
	for all $\vek{v}\in\vek{H}^{1}(\Omega)$ such that $\vek{n}\times\vek{v}=\vek{0}$ on $\Gamma$.
	
	\subsection{Homogeneous Neumann boundary conditions}
	\label{sec:Homogeneous Neumann boundary conditions rotation-free body force}
	For the homogeneous Neumann boundary conditions \specialeqref{eq:boundary conditions}{d} and the given rotation-free body force \eqref{eq:rotation-free-body-force} we solve a basic elasticity problem \eqref{eq:weak-formulation} with no-penetration with zero tangential traction boundary conditions, \specialeqref{eq:boundary conditions}{a} see Section~\ref{sec:No-penetration and zero tangential traction boundary conditions}, and obtain the solution $(\vek{u},p)=(\vek{0},\phi)$. The no-penetration part induces a normal traction react		
	Therefore,  we need to solve a corrective elasticity problem \eqref{eq:weak-formulation} with non-homogeneous Neumann data $\vek{s}=\phi\vek{n}$ on $\Gamma$. In order to remove the rigid translation and rotation we 
	add the following side conditions,
	\[
	\int_{\Omega} \vek{w} \, d\Omega = \vek{0} \quad \text{and}\quad \int_{\Omega} \vek{r}\times\vek{w} \, d\Omega = \vek{0},
	\]
	where $\vek{r}$ is a fixed radius vector.
	
	Then, analogously to the no-slip with zero normal traction case we obtain the total solution $(\vek{u},p)= (\vek{w},\phi+r)$, which is confirmed by the following straight forward computation, with $\vek{v}\in\vek{H}^{1}(\Omega)$,
	\begin{equation*}
	\int_{\Omega}\vek{f}\cdot\vek{v}\,d\Omega = -\int_{\Omega}\phi\,\div{\vek{v}}\,d\Omega + \int_{\Gamma}\phi\,\vek{n}\cdot\vek{v}\,d\Gamma  \overset{\eqref{eq:weak-formulation}}{=} \int_{\Omega}(\ten{e}(\vek{v}):[2\mu\,\ten{e}(\vek{w})] - (\phi+r)\div{\vek{v}})\,d\Omega.
	\end{equation*}
	
	\subsection{Special case of potential with zero trace}
	\label{sec:special rotation-free body force}
	If the potential $\phi$ is zero on the boundary, i.e., $\phi\in H^1_0(\Omega)$, there always holds, independently of (homogeneous) boundary conditions that,
	\begin{equation*}
	\int_{\Omega}\vek{f}\cdot\vek{v}\,d\Omega = -\int_{\Omega}\phi\,\div{\vek{v}}\,d\Omega + \int_{\Gamma}\phi\,\vek{n}\cdot\vek{v}\,d\Gamma = -\int_{\Omega}\phi\,\div{\vek{v}}\,d\Omega.
	\end{equation*}
	Therefore, the solution is given by $(\vek{u},p) = (\vek{0},\phi)$ and no displacement will be induced, extending Euler's observation for this specific subspace of scalar potentials
	to homogeneous Neumann and no-slip with zero normal traction boundary conditions.

	\section{Numerical experiments}
	\label{sec:Numerical experiments}
	
	\subsection{The finite element constructs}
	\label{sec:The finite element constructs}
	We use two finite element formulations. {For simplicity we assume homogeneous Dirichlet boundary conditions in this subsection.} The main corresponds to the $(\vek{u},p)\in\vek{H}^{1}_{0}(\Omega)\times L^{2}(\Omega)$ formulation for linear 
	elasticity problem \eqref{eq:weak-formulation} and the other is the $(\vek{A},\pi)\in \vek{H}_{0}(\mathrm{curl},\Omega)\times H^1_0(\Omega)$ formulation for the Helmholtz decomposition \eqref{eq:weak Helmholtz formulation A}.
	
	In the first formulation we use \emph{inf-sup} stable, variable order, fully parametric nodal elements with discontinuous pressure of the type $Q^{\phantom{d}}_{p}/P^{disc}_{p-1}$, with polynomial degree $p\ge 2$ \cite[Section 3.2.6]{John2016}.
    For hexahedra the reference element is $\hat{\Omega}=[0,1]^{3}$. On $\hat{\Omega}$ the usual space of tensorial product of polynomials of degree $p$ is,
    \begin{equation}
     Q_{p}(\hat{\Omega})=\textrm{span}\{\xi^{i}\cdot\eta^{j}\cdot\zeta^{k}\,:\, 0 \le i,j,k \le p\}, \quad 
     \dim Q_{p}(\hat{\Omega}) = (p+1)^{3},
     \label{eq:nodal-polynomial-space}
    \end{equation}
     Further,
    \begin{equation}
     P^{disc}_{p-1}(\hat{\Omega}) = \textrm{span}\{\xi^{i}\cdot\eta^{j}\cdot\zeta^{k}\,:\,i,j,k \ge 0, \; i+j+k \le p-1\}, \quad
     \dim P^{disc}_{p-1}(\hat{\Omega}) = {1}/{6}\,p(p+1)(p+2),
     \label{eq:discontinuous-pressure-space}
    \end{equation}
	where \eqref{eq:nodal-polynomial-space}\textsubscript{2}  and \eqref{eq:discontinuous-pressure-space}\textsubscript{2} tell the maximum dimensions, for spatially uniform $p$-extension, respectively.
	The domain $\Omega$ is partitioned into $\mathrm{nel}$ elements $\Omega_{e}$, such that $\Omega=\bigcup_{e=1}^{\mathrm{nel}}\Omega_{e}$. We use the bijective mapping $\bpsi_{e}:\hat{\Omega}\mapsto\Omega_{e}$ such that functions 
	$\vek{v}_{e}\in[Q_{p}(\Omega_{e})]^{3}$ are given by $\vek{v}_{e}(\vek{x})=\hat{\vek{v}}(\bxi)\circ\bpsi_{e}^{-1}:\Omega_{e}\mapsto\mathbbm{R}^{3}$ with $\hat{\vek{v}}(\bxi)\in[Q_{p}(\hat{\Omega})]^{3}$. The discrete displacement space $\vek{U}_{h}(\Omega)\subset\vek{U}(\Omega)$ is constructed as,
	\begin{gather}
	 \vek{U}_{h}(\Omega) =\left\{\vek{u}_{h}\in[C^{0}(\Omega)]^{3}\,\vert\;\vek{u}_{h}\vert_{\Omega_{e}}=\hat{\vek{u}}\circ\bpsi^{-1}_{e}; \quad \hat{\vek{u}} \in[Q_{p}(\hat{\Omega})]^{3}; \quad \forall\Omega_{e}\in\Omega\,\right\}.
	\end{gather}
	The discrete pressure space $Q_{h}(\Omega)\subset Q(\Omega)$ is constructed similarly.
	The geometry, the displacement and the pressure  are using the same parametric mapping. This gives us an iso-parametric element.

	In the second formulation we use \emph{inf-sup} stable fully parametric edge-elements which are variable order generalisations of $\vek{H}(\curl)$-conforming N{\'e}d{\'e}lec elements \cite{Nedelec1980,Nedelec1986} described in \cite{WRLDEM00,LDEM2007,LDEM2008,FalkGattoMonk2011,Zaglmayr06}.
	The construction of the discrete $\vek{H}(\curl)$-conforming displacement vector space for the hexahedron of 
	the N{\'e}d{\'e}lec first type can be found in 
    Sections 2.1.2 and 2.1.5 of \cite{LDEM2008}. Here we denote the vector space of tensorial product of
    polynomials on the reference element $\vek{Q}_{p}(\hat{\Omega})$ \cite[Eq. (2.7)\textsubscript{2}]{LDEM2008}. The Lagrange multiplier $\pi$ is approximated 
	using $Q_{p+1}(\hat{\Omega})$ polynomials according to  \eqref{eq:nodal-polynomial-space}. The mixed finite elements for the Helmholtz problem \eqref{eq:weak Helmholtz formulation A} may therefore
	be described being of the type $\vek{Q}_{p}/Q_{p+1}$.
	
	\subsection{The cube problems}
	\label{sec:the cube problems}
	We consider two problems on the cube geometry $\Omega = [-1,+1]^{3}$. The cube is assumed to be filled with an isotropic linear elastic, incompressible material with shear modulus 
	$\mu > 0$.  The two cases are,
	\begin{enumerate}
	\item[a)] clamped boundary conditions, and 
	\item[b)] no-penetration and zero  tangential traction boundary conditions.
	\end{enumerate}
	The results for case a) are presented in Section~\ref{sec:the clamped cube} and the results for case b) are presented in Section~\ref{sec:the cube with no-penetration and zero tangential traction boundary conditions}. 
	
	Each case is considered using the three step procedure described in Section~\ref{sec:A three step procedure relying on a Helmholtz decomposition of the loading}. Step 0 is common to both cases. The following analytically given \emph{total body force} is used,
	\begin{equation}
	\vek{f}=\mu\vek{f}_{1} + \vek{f}_{2}
	\label{eq:Body force with a solenoidal and a rotation-free part}
	\end{equation}
	where $\vek{f}_{1}$  and $\vek{f}_{2}$ are given in \ref{sec:The assumed body force on the cube} by \eqref{eq:closed-form-divergence-free-body-force}
	 and by \eqref{eq:closed-form-rotation-free-body-force}, respectively.

	The two elasticity boundary value problems a) and b) are solved using the formulation in Problem~\ref{prob:LE}, first for the total body force \eqref{eq:Body force with a solenoidal and a rotation-free part} and then for the right hand sides
	specified by Step 1 and Step 2 in the three step procedure, described in Section~\ref{sec:A three step procedure relying on a Helmholtz decomposition of the loading}.
	Each case is solved for $\mu=1$ and $\mu=10^{-4}$. Using the Perturbed Lagrangian approach realising Problem~\ref{prob:LE}, see Remark~\ref{rem:Perturbed Lagrangian}, near incompressibility 
	is henceforth represented  by the constant bulk modulus to shear modulus ratio $\kappa/\mu=10^{7}$.
	
	\subsubsection{The mesh:}
	For both cases a) and b) a uniform mesh with $10\times 10\times 10$ elements with uniform polynomial extension of order four ($p=4$) with respect to the displacement is used.
	The pressure is approximated using one order lower discontinuous shape functions. The Helmholtz decomposition Problem~\ref{prob:HD} 
	has totally 270.641 $\vek{A}$ plus $\pi$ degrees of freedom.    The elasticity Problems~\ref{prob:LE} have totally 206.763 displacement degrees of freedom\footnote{Including the clamped degrees of freedom.}. 
	
	\subsection{Helmholtz decomposition of the total body force}
	\label{sec:Helmholtz decomposition of the total body force}
	This is the initial step in the three step procedure described in 
	Section~\ref{sec:A three step procedure relying on a Helmholtz decomposition of the loading}. Problem~\ref{prob:HD} is solved for the vector potential
	${\vek{A}}$ and the Lagrange multiplier $\pi$ which enforces the divergence-free constraint on $\vek{A}$. The divergence-free part ${\curl\vek{A}}$ of the total body force $\vek{f}$ is then computed, and the rotation-free part as $\vek{f}-{\curl\vek{A}}$. These are used as right-hand sides 
	in Step~1 and Step~2, respectively. The result for the clamped cube example described in Section~\ref{sec:the clamped cube} is given in Table~\ref{tab:Step0-results}
	for $\mu=1$. That is, for no contrast between the components of the body force \eqref{eq:Body force with a solenoidal and a rotation-free part}. 
	
	We summarise the observations regarding the solution $({\vek{A}},\pi)$ as,
	\begin{itemize}
		\item[$\bullet$] The numerical results for $\mu=10^{-4}$ are visually the same as those 
		for $\mu=1$, and are therefore not shown.
		\item[$\bullet$] Problem~\ref{prob:HD} is realised as a pure saddle point problem.
		\item[$\bullet$] We note that the Lagrange multiplier $\pi$ is very cleanly numerically zero.
	\end{itemize}

	\subsection{The clamped cube}
	\label{sec:the clamped cube}
	In Section~\ref{sec:Results with a total body force} we consider results with a total body force.
	In Section~\ref{sec:Divergence free body force} we discuss the results corresponding to Step 1,
	and  in Section~\ref{sec:Rotation free body force} we analyze the results corresponding to Step 2, in the three step procedure respectively.

	\subsubsection{Results with a total body force}
	\label{sec:Results with a total body force}
		
	The $(\vek{u},p)$  solutions for a clamped cube subject to a total body force $\vek{f}$ given by \eqref{eq:Body force with a solenoidal and a rotation-free part} given in closed form in \ref{sec:The assumed body force on the cube}, are shown in Table~\ref{tab:total-body-force} (of figures). 
	The left column of the table shows results for $\mu=1$ and the right column shows those for $\mu=10^{-4}$. The value of $\mu$
	also implies a contrast between the divergence-free part and the rotation-free part according to equation \eqref{eq:Body force with a solenoidal and a rotation-free part}.
	Comparing the determined displacement component $u_{x}$, shown as an iso-level map  in Figure (a) and and as a three-dimensional vector graph in Figure (b), with the corresponding results in Figures  (e) and (f) reveals the
	loss of pressure robustness changing form $\mu=1$ to $\mu=10^{-4}$. It is seen as a blurriness.
	A similar blurriness is seen in the displacement component $u_{y}$. There is a visible difference also in the pressure reaction $p$. In Figures (c) and (d) we can anticipate a pressure reaction caused by the divergence-free body force component. While in Figures (g) and (h) this contribution is completely suppressed. Only the pressure corresponding to the rotation-free body force is visible. In the following sections we will show that our three step procedure enables the retrieval of the divergence-free displacement solution without pollution. It also allows 
	one to separate the pressure reaction with respect to its sources.
	
	The main observation is,
	\begin{itemize}
	 \item[$\bullet$] our higher order mixed $(\vek{u},p) \in \vek{H}^{1}_{0}(\Omega)\times L^{2}_{0}(\Omega)$ formulation is not pressure robust in the sense \cite{LedererSIAM2017}, i.e. it looses its ability to compute the weakly divergence-free displacement solution $\vek{u}$ pollution free at high contrasts between the divergence-free part and the rotation-free part of the body force. Here the contrast is caused scaling the divergence-free component by $\mu$, according to \eqref{eq:Body force with a solenoidal and a rotation-free part}. The contrast is of four magnitudes. The observation concerns $p=4$ for the displacements and one order lower for the pressure. 
	 We note that the constant shear to bulk modulus ratio $\mu/\kappa=10^{-7}$, we use for both cases, corresponds to a constant Poisson's ratio\footnote{$\nu=(1/2-1/3\,\mu/\kappa)/(1+1/3\,\mu/\kappa)$.} of $\nu=0.49999995$. This is how we approximate a near incompressible behaviour, see also Remark~\ref{rem:Perturbed Lagrangian}.
	\end{itemize}

	\subsubsection{Step 1: Excitation by a divergence-free body force}
	\label{sec:Divergence free body force}
	We determine the solution $(\vek{u}_{1},p_{1})$ to problem \eqref{eq:weak-formulation} for the cube with clamped 
	boundary conditions subject to the divergence-free body force $\vek{f}_{1}=\curl\vek{A}_{h}$ determined from 
	the numerically computed vector potential $\vek{A}_{h}$ in Step 0 of the three step procedure.
	Iso maps and three-dimensional vector graphs of the solution components $({u}_{1\,x},p_{1})$ for $\mu=1$ and $\mu=10^{-4}$ are shown in the left and right columns of Table~\ref{tab:solenoidal-body-force} (of figures), respectively.
	
	The solution $(\vek{u}_{1},p_{1})$ can be commented in the following way,
	\begin{itemize}
		\item[$\bullet$] The displacement solution $\vek{u}_{1}$ is computed \emph{pressure robustly} in the sense of \cite{LedererSIAM2017}.  The three step procedure, described in Section~\ref{sec:A three step procedure relying on a Helmholtz decomposition of the loading} is used. The displacement solution does not become polluted as $\mu$ is changed four orders of magnitude. The shown component $u_{1\,x}$ in Figures (a) and (b) and  (e) and (f), corroborate the conclusion.
		The component $u_{1\,y}$ is also computed without pollution for $\mu=10^{-4}$.
		\item[$\bullet$] The pressure reaction $p_{1}$ shown in Figures (c) and (d) and in (g) and (h) is non-trivial.    
		 It is concluded that the divergence-free body force $\curl\vek{A}_{h}$  \eqref{eq:closed-form-divergence-free-body-force} induces a non-trivial traction reaction on the clamped boundaries. Its normal part $-p_{1}\vek{n}$ corresponds to the pressure reaction.	
	\end{itemize}

	\subsubsection{Step 2: Excitation by a rotation-free body force}
	\label{sec:Rotation free body force}
	We determine $(\vek{u}_{2},p_{2})$ of problem \eqref{eq:weak-formulation} for the excitation corresponding to 
	the numerically computed rotation-free body force $\vek{f}_{2}=\vek{f}-\curl\vek{A}_{h}$ solving problem \eqref{eq:weak Helmholtz formulation A} for the vector potential $\vek{A}_{h}$.
	The solutions $p_{2}$ for $\mu=1$ and $\mu=10^{-4}$ are shown in the left and right columns in Table~\ref{tab:rotation-free-body-force} (of figures), respectively.
	
	The solution $(\vek{u}_{2},p_{2})$ can be commented in the following way,
	\begin{itemize}
		\item[\textbullet] The displacement $\vek{u}_{2}$ is numerically trivial and therefore not shown. 
		The computed rotation-free body force $\vek{f}-\curl\vek{A}_{h}$ is very close to $\textrm{grad\,}\phi$, and therefore there is no observed pollution of the numerically computed rotation-free body force.
		\item[\textbullet] The pressure $p_{2}$ is clearly non-trivial. 
		We note that since we determine the vector potential $\vek{A}_{h}$ by the Helmholtz decomposition \eqref{eq:Helmholtz-decomposition} we cannot 
		compute the corresponding scalar potential $\phi$ by it self to confirm Euler's observation. That would require another Helmholtz decomposition solving Poisson's problem 
		for the scalar potential $\phi$.
	\end{itemize}
	
	\subsection{The cube with no-penetration and zero tangential traction boundary conditions}
	\label{sec:the cube with no-penetration and zero tangential traction boundary conditions}

	In Section~\ref{sec:Results with a total body force no-penetration} we consider results with a total body force.
	In Section~\ref{sec:Divergence free body force no-penetration} we discuss the results corresponding to Step 1,
	and  in Section~\ref{sec:Rotation free body force no-penetration}) we analyze the results corresponding to Step 2, in the three step procedure respectively.	
	
	\subsubsection{Results with a total body force}
	\label{sec:Results with a total body force no-penetration}
	The $(\vek{u},p)$  solutions for a cube with no-penetration and zero tangential traction boundary conditions subject to a total body force $\vek{f}$ given by \eqref{eq:Body force with a solenoidal and a rotation-free part}, and in 
	closed form in the \ref{sec:The assumed body force on the cube}, are shown in Table~\ref{tab:total-body-force-no-penetration} (of figures). Also here, the left column of the table shows results for $\mu=1$ and the right column shows those for $\mu=10^{-4}$. The contrast between the divergence-free part and the rotation-free part is achieved according to equation \eqref{eq:Body force with a solenoidal and a rotation-free part}.
	Also in this case, cf. Section~\ref{sec:Results with a total body force}, we loose pressure robustness changing from 
	$\mu=1$ to $\mu=10^{4}$, compare Figures (a) and (b) with (e) and (f), respectively. In this case we do not encounter any indication of a pressure reaction induced by the divergence-free part, see Figures (c) and (d), and (g) and (h).
	
	The main observation concerning loss of pressure robustness for this case is shared with that in Section~\ref{sec:Results with a total body force}.
	
	\subsubsection{Step 1: Excitation by a divergence-free body force}
	\label{sec:Divergence free body force no-penetration}
	We determine the solution $(\vek{u}_{1},p_{1})$ to problem \eqref{eq:weak-formulation} for the cube with 
	no-penetration and zero tangential traction boundary conditions subject to the divergence-free body force $\vek{f}_{1}=\curl\vek{A}_{h}$ determined from 	the numerically computed vector potential $\vek{A}_{h}$ in Step 0 of the three step procedure.
	Iso maps and three-dimensional vector graphs of the solution components $({u}_{1\,x},p_{1})$ for $\mu=1$ and $\mu=10^{-4}$ are shown in the left and right columns of Table~\ref{tab:solenoidal-body-force-no-penetration} (of figures), respectively.
	
	The results can be commented in the following way,
	\begin{itemize}
		\item[$\bullet$] The displacement solution $\vek{u}_{1}$ is computed \emph{pressure robustly} in the sense of \cite{LedererSIAM2017}.  Using the three step procedure, described in Section~\ref{sec:A three step procedure relying on a Helmholtz decomposition of the loading}. It does not become polluted as $\mu$ is changed four orders of magnitude. The shown component $u_{1\,x}$ in Figures (a) and (b) and  (e) and (f), corroborate the conclusion.
		The component $u_{1\,y}$ is also computed without pollution for $\mu=10^{-4}$.
		\item[$\bullet$] We underline that the pressure reaction $p_{1}$ shown in Figures (c) and (d) and in (g) and (h) 
		 is numerically trivial. As opposed to clamped boundary conditions, cf. Section~\ref{sec:Divergence free body force}, it is concluded that the divergence-free body force $\curl\vek{A}_{h}$  \eqref{eq:closed-form-divergence-free-body-force} which is parallel to the boundary, $\vek{n}\cdot\vek{f}_{1}=0$, does not induce any reaction in case of no-penetration and zero tangential traction boundary conditions.
	\end{itemize}

	\subsubsection{Step 2: Excitation by a rotation-free body force}
	\label{sec:Rotation free body force no-penetration}
	We determine the solution $(\vek{u}_{1},p_{1})$ to problem \eqref{eq:weak-formulation} for the cube with 
	no-penetration and zero tangential traction boundary conditions subject to the rotation-free body force $\vek{f}_{2}=\vek{f}-\curl\vek{A}_{h}$ determined using the numerically computed vector potential $\vek{A}_{h}$ in Step 0 of the three step procedure.
	The solutions $p_{2}$ for $\mu=1$ and $\mu=10^{-4}$ are shown in Table~\ref{tab:rotation-free-body-force-no-penetration} (of figures).
	
	The observations concerning the results for this case are similar to those for the clamped cube,
	see Section~\ref{sec:Rotation free body force}
	
	\subsection{The three-dimensional L-shaped domain}
	\label{sec:The three-dimensional L-shaped domain}
	Consider a material three-dimensional L-shaped domain $\Omega = [-1,+1]^{3}\setminus[0,+1]^{3}$, also known as the Fichera problem.
	The three-dimensional L-shaped domain is not convex and that it has re-entrant corners.
	We apply homogeneous no-penetration and zero tangential traction boundary conditions, \mbox{$\vek{n}\cdot\vek{u}=\vek{0}$} and \mbox{$\vek{n}\times\vek{s}=\vek{0}$} for $\vek{x}\in\Gamma$, respectively. We know from the theory and the reported numerical results that these boundary conditions do not induce a pressure reaction in combination with a divergence-free body force on a convex domain.  
	
	Here, due to the re-entrant corners, the solution $(\vek{u},p)$ is expected to be non-regular although the forcing is regular. We	are especially interested to see if we obtain a non-trivial singular pressure reaction.
	
	To this end, we assume the vector potential $\vek{A}=100\cdot[p(y)p(z),\;p(x)p(z),\;p(x)p(y)]$ where	$p(r)=r^{2}(r^{2}-1)^{2}$.	It is readily checked that $\vek{A}$ fulfills the Coulomb gauge, $\div{\,\vek{A}}=0$. Furthermore, it induces a divergence-free body force, computed as, 
	\begin{equation}
	 \vek{f}_{1}=\curl{\vek{A}} = 200\left[-p(x)\,
   q(z,y),\; p(y)\,q(z,x),\; -p(z)\,q(y,x) \right]^{\T}
   \label{eq:divergence-free body force L-shaped domain}
	\end{equation}
   where $q(s,t)=3\,s^5-4\,s^3+s-3\,t^5+4\,t^3-t$. Furthermore, it is readily checked that this body force is parallel to the all boundary walls, $\vek{n}\cdot\curl\vek{A}=0$. Finally, we let the three-dimensional L-shaped domain be filled with an isotropic linear elastic, incompressible material with shear modulus $\mu > 0$.
   
   We have no ambitions to resolve the present singularities only to investigate
   if the pressure reaction will be non-trivial. We therefore allow us to use a simple ignorant mesh without grading. It is composed of 7 cubical blocks. Each block is subdivided uniformly into 10$\times$10$\times$10 hexes. Thus the model contains 7000 hexes. We likewise use an ignorant uniform polynomial extension $p=4$. The model has totally 1.402.323 degrees of freedom.
   
   The result for the pressure reaction is shown in Table~\ref{tab:Lshaped-divergence-free-body-force-no-penetration}
   (of figures) as an iso-level map in Figure (a) and as a three-dimensional vector graph in Figure (b). 
   It is evident that the pressure reaction is non-trivial. It is singular along the re-entrant corners. Figure (b)
   shows its singular nature. The 7$\times$10$\times$10$\times$10 mesh is the finest we used. The pressure reaction 
   converged as the mesh was refined.  We conclude that not only the boundary conditions, but also the regularity of the domain decides if a divergence-free body force will induce a pressure reaction.

	\section{Summary and conclusions}
	\label{sec:Summary and conclusions}
	In isotropic incompressible linear elasticity the divergence-free displacement and the pressure reaction solution components
	$(\vek{u},p)$ cannot be determined independently in general. The majority of investigated boundary conditions do not allow such a decoupling.
	The situation further depends on the nature of the body force, i.e., if it is divergence-free or rotation-free, respectively.
	Moreover, the out-come also depends on the regularity of the solution $(\vek{u},p)$. For example, for the standard 
	three-dimensional L-shaped domain with re-entrant corners we encounter a non-trivial singular $p$ solution for the case with no-penetration boundary conditions. Discontinuities in the imposed forcing may also prevent an independent determination of $\vek{u}$ and $p$.
	
	For a convex cube domain $\Omega$ we showed that for a divergence-free body force we can determine the weakly divergence-free displacement $\vek{u}$ independent of the pressure reaction for mixed-mixed boundary conditions of no-penetration with zero tangential traction as well as for no-slip with zero normal traction. On the other hand for a rotation-free body force we can determine the pressure reaction $p$ independently of the displacement for clamped and for no-penetration with zero tangential traction boundary conditions. 
	That is, for the single case of no-penetration with zero tangential traction  boundary conditions a divergence-free body force gives the solution 
	$(\vek{u},0)$ while the rotation-free body force gives the solution $(\vek{0},p)$. We note that the no-slip with homogeneous normal traction boundary conditions and a rotation free body force in general gives a non-trivial displacement as well as a pressure reaction. Thus preventing an independent determination of $\vek{u}$ and $p$. 
	
	An $L^{2}(\Omega)$-orthogonal Helmholtz decomposition  of the given body force into its divergence- and rotation-free parts and the superposition principle provides the means for a \emph{pressure robust} determination of $\vek{u}$ for the following boundary conditions:
	\begin{enumerate}
	 \item mixed-mixed, no-penetration and homogeneous tangential traction boundary conditions, or 
	 \item mixed-mixed, no-slip and homogeneous normal traction boundary conditions, or
	 \item pure essential boundary conditions, or
	 \item pure homogeneous natural boundary conditions,
	\end{enumerate}
	are prescribed on the whole boundary.
	An \emph{independent} determination of $\vek{u}$ and $p$ in the case of no-penetration with zero tangential traction boundary conditions. Here we have chosen to determine the vector potential $\vek{A}$ imposing the Coulomb gauge 
	$\div{\,\vek{A}}=0$ weakly, using a continuous $H^{1}_{0}(\Omega)$-conforming Lagrange multiplier and $\vek{H}_{0}({\rm curl},\Omega)$-conforming formulation for $\vek{A}$, i.e., equation \eqref{eq:weak Helmholtz formulation A}. This formulation is shown to have a trivial Lagrange multiplier which was observed in electromagnetics \cite{WRLDEM00} and moreover it becomes analogous to the auxiliary elasticity formulation \eqref{eq:auxiliary curl-curl elasticity formulation} for no-slip with homogeneous traction boundary conditions and a divergence-free body force $\vek{f}=\curl\vek{A}$. It is worth noting that the elasticity inner product $\inner{\ten{e}(\vek{v})}{\ten{e}(\vek{u})}$ and $\inner{\curl\vek{v}}{\curl\vek{u}}$ analogy is restricted to convex domains. The trivial pressure in elasticity
	for these boundary conditions and loading, is valid for convex domains.
	
	All other cases of boundary conditions yield in general non-trivial displacement and pressure for a divergence- or a rotation-free body force, preventing an independent determination. 
	But importantly, we may still determine the displacement 
	\emph{pressure robustly} in the sense of \cite{LedererSIAM2017} for the stated enumerated boundary conditions, using the mentioned $L^{2}(\Omega)$-orthogonal Helmholtz decomposition  of the given body force. Without that decomposition our \emph{inf-sup} stable higher-order but otherwise \emph{standard mixed displacement-pressure method} with discontinuous pressure elements of the type $Q_{k}/P^{disc}_{k-1}$, \mbox{$k\ge 0$} with $(\vek{u},p) \in \vek{H}^{1}(\Omega)\times L^{2}(\Omega)$  will fail to determine the displacement pressure robustly at high contrasts between the divergence- and a rotation-free components of the body force.
	
	In particular, clamped boundary conditions yield a non-trivial partial solution $(\vek{u},p)$ for a divergence-free body force.  For a rotation-free body force we obtain a partial solution $(\vek{0},p)$ with clamped conditions.
	Further, for a rotation-free body force we notably obtain a non-trivial weakly divergence-free displacement for mixed-mixed no-slip and homogeneous normal traction boundary conditions. On the other hand, a divergence-free body force and clamped boundary conditions give a non-trivial pressure as already pointed-out. These results underline that Euler's observation \cite{CFT60}(1757) is restricted to clamped and no-penetration with zero tangential traction  boundary conditions and a rotation-free body force. We have already pointed out that no-slip boundary conditions and homogeneous normal traction conditions and a rotation free body force induces a non-trivial displacement and a pressure. Only in the special case of a gradient potential which vanishes on the boundary does not give rise to a displacement field independently of the prescribed (homogeneous) boundary conditions.
	
	Prescribing no-penetration with zero tangential traction and no-slip with homogeneous normal traction boundary conditions on separate parts of the boundary still yields a non-zero pressure in case of a divergence-free body force. For clamped and no-slip with homogeneous normal traction boundary conditions a rotation-free body force does not induce a non-zero displacement field. Other combinations, however, will induce a non-zero pressure and displacement field in general for both, divergence- or rotation-free body forces.
	
	\section*{Acknowledgments}
Michael Neunteufel acknowledges support by the Austrian Science Fund (FWF) project F65. The authors are also indebted to Joachim Sch\"oberl for interesting discussions.
	
	\bibliographystyle{elsarticle-num}
	\bibliography{cites_paper}

\begin{thebibliography}{10}
\expandafter\ifx\csname url\endcsname\relax
  \def\url#1{\texttt{#1}}\fi
\expandafter\ifx\csname urlprefix\endcsname\relax\def\urlprefix{URL }\fi
\expandafter\ifx\csname href\endcsname\relax
  \def\href#1#2{#2} \def\path#1{#1}\fi

\bibitem{ChorinMarsden1992}
A.~J. Chorin, J.~E. Marsden, A mathematical introduction to fluid mechanics,
  Vol.~3, Springer, 1990.

\bibitem{Lad69}
O.~A. Ladyzhenskaya, The mathematical theory of viscous incompressible flow,
  Vol.~2, Gordon and Breach New York, New York, 1969.

\bibitem{Bab73}
I.~Babu{\v{s}}ka, The finite element method with {L}agrangian multipliers,
  Numerische Mathematik 20~(3) (1973) 179--192.
\newblock \href {https://doi.org/10.1007/BF01436561}
  {\path{doi:10.1007/BF01436561}}.

\bibitem{Bre74}
F.~Brezzi, On the existence, uniqueness and approximation of saddle-point
  problems arising from {L}agrangian multipliers, R.A.I.R.O. Analyse
  Num\'erique 8~(R2) (1974) 129--151.
\newblock \href {https://doi.org/10.1051/m2an/197408R201291}
  {\path{doi:10.1051/m2an/197408R201291}}.

\bibitem{Linke2014}
A.~Linke, On the role of the {H}elmholtz decomposition in mixed methods for
  incompressible flows and a new variational crime, Computer Methods in Applied
  Mechanics and Engineering 268 (2014) 782--800.
\newblock \href {https://doi.org/https://doi.org/10.1016/j.cma.2013.10.011}
  {\path{doi:https://doi.org/10.1016/j.cma.2013.10.011}}.

\bibitem{Linke2016}
A.~Linke, C.~Merdon, Pressure-robustness and discrete {H}elmholtz projectors in
  mixed finite element methods for the incompressible {N}avier-{S}tokes
  equations, Computer Methods in Applied Mechanics and Engineering 311 (2016)
  304--326.
\newblock \href {https://doi.org/https://doi.org/10.1016/j.cma.2016.08.018}
  {\path{doi:https://doi.org/10.1016/j.cma.2016.08.018}}.

\bibitem{JohnSIAM2017}
V.~John, A.~Linke, C.~Merdon, M.~Neilan, L.~G. Rebholz, On the divergence
  constraint in mixed finite element methods for incompressible flows, SIAM
  Review 59~(3) (2017) 492--544.
\newblock \href {https://doi.org/10.1137/15M1047696}
  {\path{doi:10.1137/15M1047696}}.

\bibitem{LedererSIAM2017}
P.~L. Lederer, A.~Linke, C.~Merdon, J.~Sch\"{o}berl, Divergence-free
  reconstruction operators for pressure-robust {S}tokes discretizations with
  continuous pressure finite elements, SIAM Journal on Numerical Analysis
  55~(3) (2017) 1291--1314.
\newblock \href {https://doi.org/10.1137/16M1089964}
  {\path{doi:10.1137/16M1089964}}.

\bibitem{LS15}
C.~Lehrenfeld, J.~Sch{\"o}berl, High order exactly divergence-free hybrid
  discontinuous {G}alerkin methods for unsteady incompressible flows, Computer
  Methods in Applied Mechanics and Engineering 307 (2016) 339--361.
\newblock \href {https://doi.org/10.1016/j.cma.2016.04.025}
  {\path{doi:10.1016/j.cma.2016.04.025}}.

\bibitem{GLS19}
J.~Gopalakrishnan, P.~L. Lederer, J.~Sch\"oberl, A mass conserving mixed stress
  formulation for the {S}tokes equations, IMA Journal of Numerical Analysis
  40~(3) (2019) 1838--1874.
\newblock \href {https://doi.org/10.1093/imanum/drz022}
  {\path{doi:10.1093/imanum/drz022}}.

\bibitem{RT77}
P.-A. Raviart, J.-M. Thomas, A mixed finite element method for 2-nd order
  elliptic problems, in: Mathematical Aspects of Finite Element Methods,
  Vol.~66, Springer, 1977, pp. 292--315.
\newblock \href {https://doi.org/10.1007/BFb0064470}
  {\path{doi:10.1007/BFb0064470}}.

\bibitem{BDM85}
F.~Brezzi, J.~Douglas, L.~D. Marini, Two families of mixed finite elements for
  second order elliptic problems, Numerische Mathematik 47~(2) (1985) 217--235.
\newblock \href {https://doi.org/10.1007/BF01389710}
  {\path{doi:10.1007/BF01389710}}.

\bibitem{FLLS21}
G.~Fu, C.~Lehrenfeld, A.~Linke, T.~Streckenbach, Locking-free and
  gradient-robust $\mathbf{H}(\,\mathrm{div}\,)$-conforming {HDG} methods for
  linear elasticiy, Journal of Scientific Computing 86~(39) (2021) 1--30.
\newblock \href {https://doi.org/10.1007/s10915-020-01396-6}
  {\path{doi:10.1007/s10915-020-01396-6}}.

\bibitem{John2016}
V.~John, Finite element methods for incompressible flow problems, Vol.~51,
  Springer, 2016.

\bibitem{CFT60}
C.~Truesdell, R.~Toupin, The classical field theories, in: Principles of
  classical mechanics and field theory/Prinzipien der Klassischen Mechanik und
  Feldtheorie, Springer, 1960, pp. 226--858.

\bibitem{Stokes1849}
G.~G. Stokes, On the dynamical theory of diffraction, Transactions of the
  Cambridge Philosophical Society 9 (1849) 1--48.

\bibitem{Helmholtz1858}
H.~v. Helmholtz, {\"U}ber {I}ntegrale der hydrodynamischen {G}leichungen,
  welche den {W}irbelbewegungen entsprechen., Crelle, 55, p. 25-55 (1858).

\bibitem{GiraultRaviart86}
V.~Girault, P.-A. Raviart, Finite element approximation of the
  {N}avier-{S}tokes equations, Vol. 749, Springer Berlin, 1986.

\bibitem{Mon03}
P.~Monk, Finite element methods for {M}axwell's equations, Numerical
  Mathematics and Scientific Computation, Oxford University Press, New York,
  2003.
\newblock \href {https://doi.org/10.1093/acprof:oso/9780198508885.001.0001}
  {\path{doi:10.1093/acprof:oso/9780198508885.001.0001}}.

\bibitem{Bercovier78}
M.~Bercovier, Perturbation of mixed variational problems. application to mixed
  finite element methods, RAIRO. Analyse num{\'e}rique 12~(3) (1978) 211--236.

\bibitem{WRLDEM00}
W.~Rachowicz, L.~Demkowicz, An hp-adaptive finite element method for
  electromagnetics: Part 1: Data structure and constrained approximation,
  Computer Methods in Applied Mechanics and Engineering 187~(1) (2000)
  307--335.
\newblock \href {https://doi.org/https://doi.org/10.1016/S0045-7825(99)00137-1}
  {\path{doi:https://doi.org/10.1016/S0045-7825(99)00137-1}}.

\bibitem{Bathiaetal2013}
H.~Bhatia, G.~Norgard, V.~Pascucci, P.-T. Bremer, The {H}elmholtz-{H}odge
  decomposition--a survey, IEEE Transactions on Visualization and Computer
  Graphics 19~(8) (2013) 1386--1404.
\newblock \href {https://doi.org/10.1109/TVCG.2012.316}
  {\path{doi:10.1109/TVCG.2012.316}}.

\bibitem{Pechstein11}
A.~Pechstein, J.~Sch{\"o}berl, Tangential-displacement and normal-normal-stress
  continuous mixed finite elements for elasticity, Math. Models Methods Appl.
  Sci. 21~(8) (2011) 1761--1782.
\newblock \href {https://doi.org/10.1142/S0218202511005568}
  {\path{doi:10.1142/S0218202511005568}}.

\bibitem{Braess07}
D.~Braess, Finite elements: {T}heory, fast solvers, and applications in
  elasticity theory, 3rd Edition, Cambridge University Press, 2007.

\bibitem{Nedelec1980}
J.~C. N{\'e}d{\'e}lec, Mixed finite elements in {R}3, Numerische Mathematik
  35~(3) (1980) 315--341.
\newblock \href {https://doi.org/10.1007/BF01396415}
  {\path{doi:10.1007/BF01396415}}.

\bibitem{Nedelec1986}
J.~C. N{\'e}d{\'e}lec, A new family of mixed finite elements in {R}3,
  Numerische Mathematik 50~(1) (1986) 57--81.
\newblock \href {https://doi.org/10.1007/BF01389668}
  {\path{doi:10.1007/BF01389668}}.

\bibitem{LDEM2007}
L.~Demkowicz, Computing with hp-adaptive finite elements: volume 1 one and two
  dimensional elliptic and Maxwell problems, 1st Edition, Chapman and Hall/CRC,
  New York, 2006.
\newblock \href {https://doi.org/10.1201/9781420011685}
  {\path{doi:10.1201/9781420011685}}.

\bibitem{LDEM2008}
J.~Kurtz, D.~Pardo, M.~Paszynski, W.~Rachowicz, A.~Zdunek, Computing with
  Hp-Adaptive Finite Elements: Volume 2, Frontiers: Three Dimensional Elliptic
  and Maxwell Problems with Applications. Applied Mathematics and Nonlinear
  Science, CRC Press, 2008.

\bibitem{FalkGattoMonk2011}
R.~S. Falk, P.~Gatto, P.~Monk, Hexahedral \mbox{$\mathbf{H}(\mathrm{div})$} and
  \mbox{$\mathbf{H}(\mathrm{curl})$} finite elements, ESAIM: Mathematical
  Modelling and Numerical Analysis 45~(1) (2011) 115–143.
\newblock \href {https://doi.org/10.1051/m2an/2010034}
  {\path{doi:10.1051/m2an/2010034}}.

\bibitem{Zaglmayr06}
S.~Zaglmayr,
  \href{https://www.numerik.math.tugraz.at/~zaglmayr/pub/szthesis.pdf}{High
  order finite element methods for electromagnetic field computation}, Ph.D.
  thesis, Johannes Kepler Universit{\"a}t Linz (2006).
\newline\urlprefix\url{https://www.numerik.math.tugraz.at/~zaglmayr/pub/szthesis.pdf}

\bibitem{Schoeberl01}
J.~Sch\"{o}berl, Commuting quasi-interpolation operators for mixed finite
  elements, Preprint ISC-01-10-MATH, Texas A\&M University, College Station, TX
  (2001).

\end{thebibliography}
	
	\clearpage

		\begin{table}[tbh]
		\begin{tabular}{|c|c|}\hline
		 \includegraphics[height=7cm, angle=90]{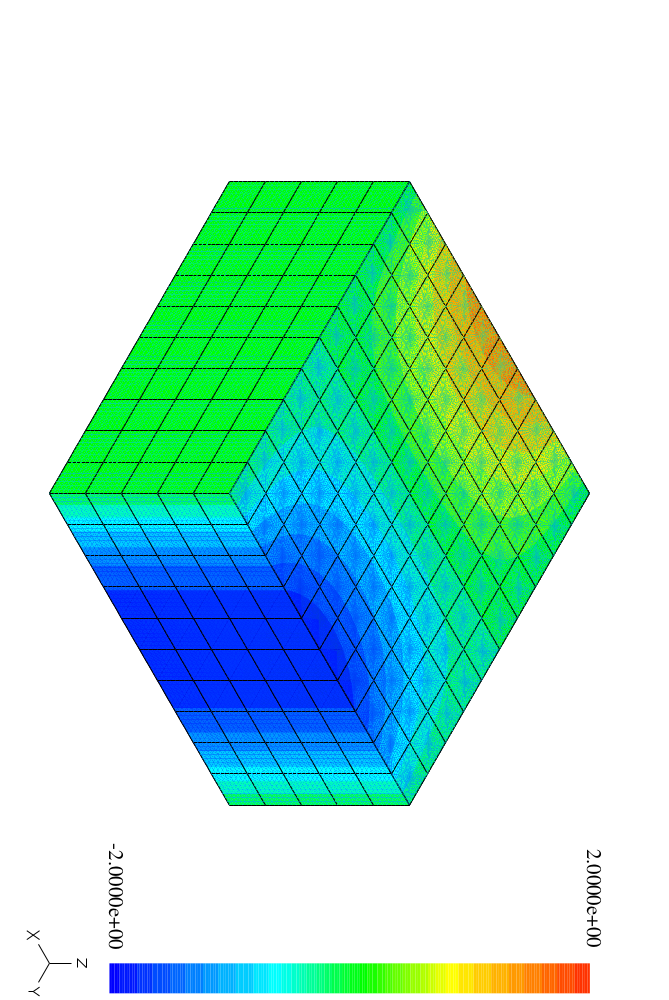}  & 
		 \includegraphics[height=7cm, angle=90]{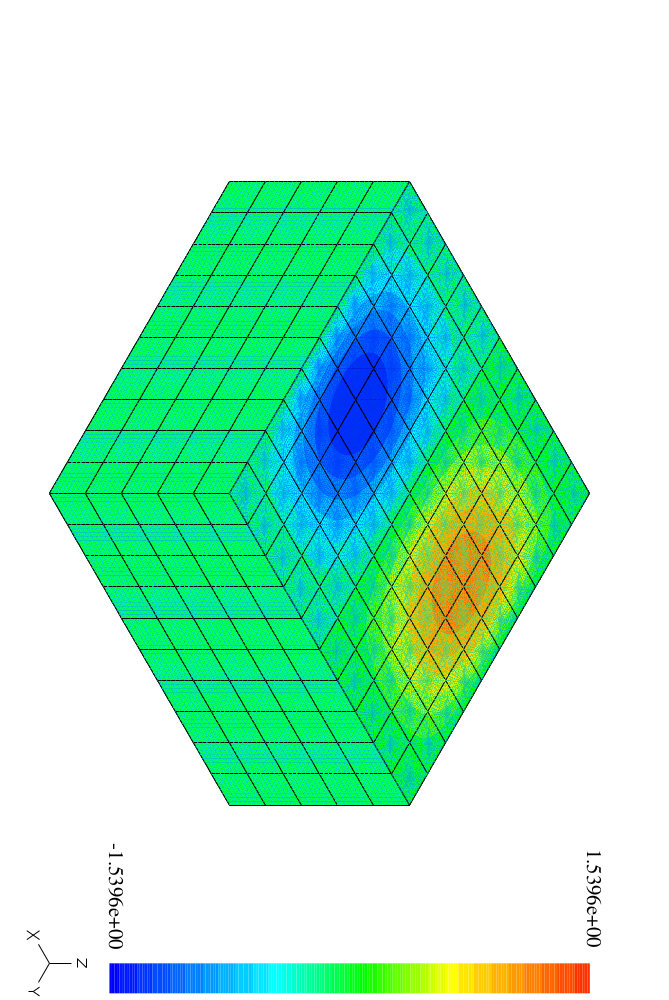}\\
		 (a) & (d) \\[-0.25ex] \hline
		 \includegraphics[height=7cm, angle=90]{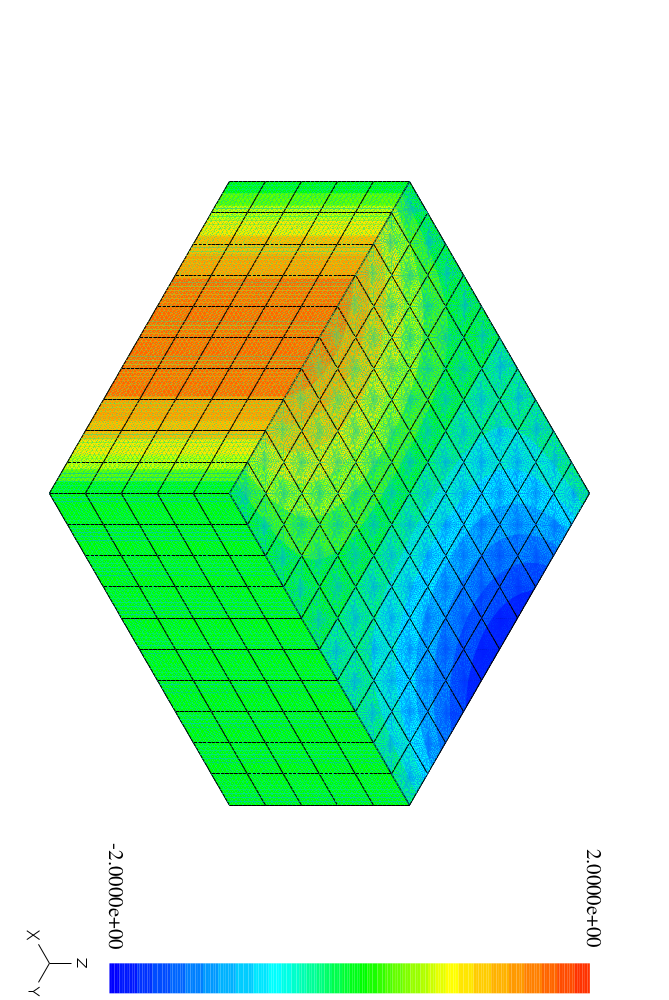}  &  
		 \includegraphics[height=7cm, angle=90]{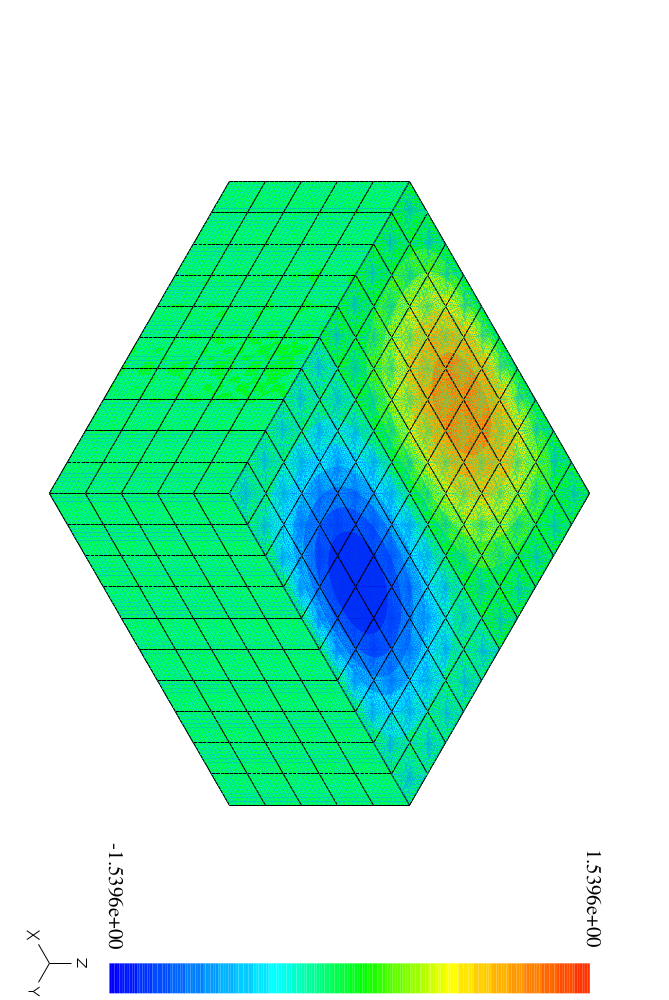}  \\
		 (b) & (e) \\[-0.25ex] \hline
		 \includegraphics[height=7cm, angle=90]{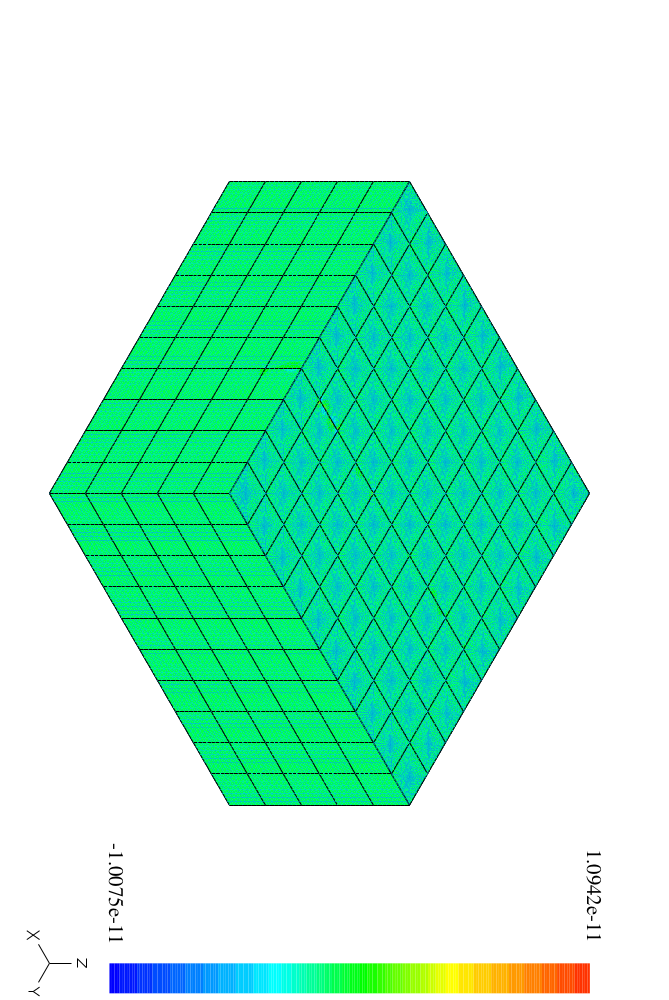}  & 
		 \includegraphics[height=7cm, angle=90]{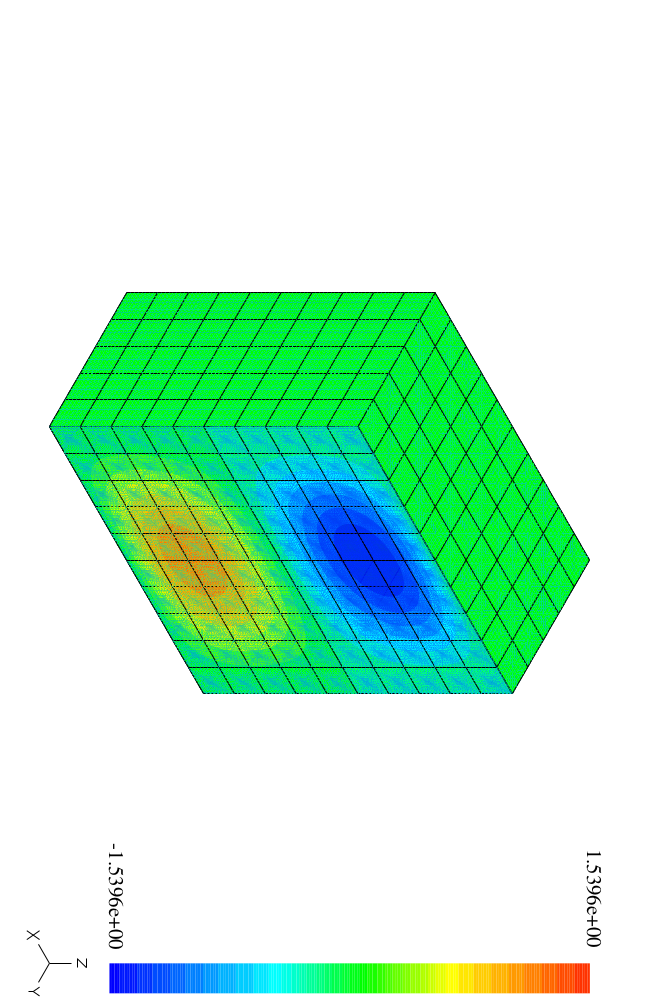}\\
		   (c) & (f) \\[-0.25ex] \hline 
		\end{tabular}
		\caption{Cube $[-1,+1]^{3}$ subject to the body force \eqref{eq:Body force with a solenoidal and a rotation-free part} including a solenoidal and a rotation-free part. Helmholtz decomposition, Step 0: computation of the vector potential $\vek{A}(\vek{x})$ and the scalar Lagrange multiplier field $\pi(\vek{x})$ solution to Problem~\ref{prob:HD}.
			Shear modulus $\mu=1$. In Figures (a), to (c) iso-maps of the components of $\curl\vek{A}$, the divergence-free part, and in (d) to (e) 
			the components of the rotation-free part $\mathrm{grad\,}\phi=\vek{f}-\curl\vek{A}$ are shown.
			Figures (a) to (e) are drawn on the plane $z=0$. Figure (f) is drawn on the plane $y=0$. All these figures are iso-level maps.\newline
			(a) component $x$ of $\curl\vek{A}$, $\min([\textrm{curl\,}A]-{x})=-2.0000$, $\max([\textrm{curl\,}A]_{x})=2.0000$.\\
			(b) component $y$ of $\curl\vek{A}$, $\min([\textrm{curl\,}A]_{y})=-2.0000$, $\max([\textrm{curl\,}A]_{y})=2.0000$.\\
			(c) component $z$ of $\curl\vek{A}$, \mbox{$\min([\textrm{curl\,}A]_{z})=-1.0075e-11$}, 
			\mbox{$\max([\textrm{curl\,}A]_{z})=1.0094e-11$}.\\
			(d) component $x$ of $\textrm{grad\,}\phi$,\;  $\min([\mathrm{grad\,}\phi]_{x})=-1.5396$, $\max([\mathrm{grad\,}\phi]_{x})=1.5396$.\\
			(e) component $y$ of $\textrm{grad\,}\phi$,\; $\min([\mathrm{grad\,}\phi]_{y})=-1.5396$, $\max([\mathrm{grad\,}\phi]_{y})=1.5396$\\  
			(f) component $z$ of $\textrm{grad\,}\phi$,\;  $\min([\mathrm{grad\,}\phi]_{z})=-1.5396$, $\max([\mathrm{grad\,}\phi]_{z})=1.5396$. \\
			Here we interpret the component $[\textrm{curl\,}A]_{z}$ as numerically zero, consistent with the analytic expression \eqref{eq:closed-form-divergence-free-body-force}\textsubscript{1}.
			\label{tab:Step0-results}}
	\end{table}

	\begin{table}[tbh]
		\begin{tabular}{|c|c|}\hline
		 \includegraphics[height=7cm, angle=90]{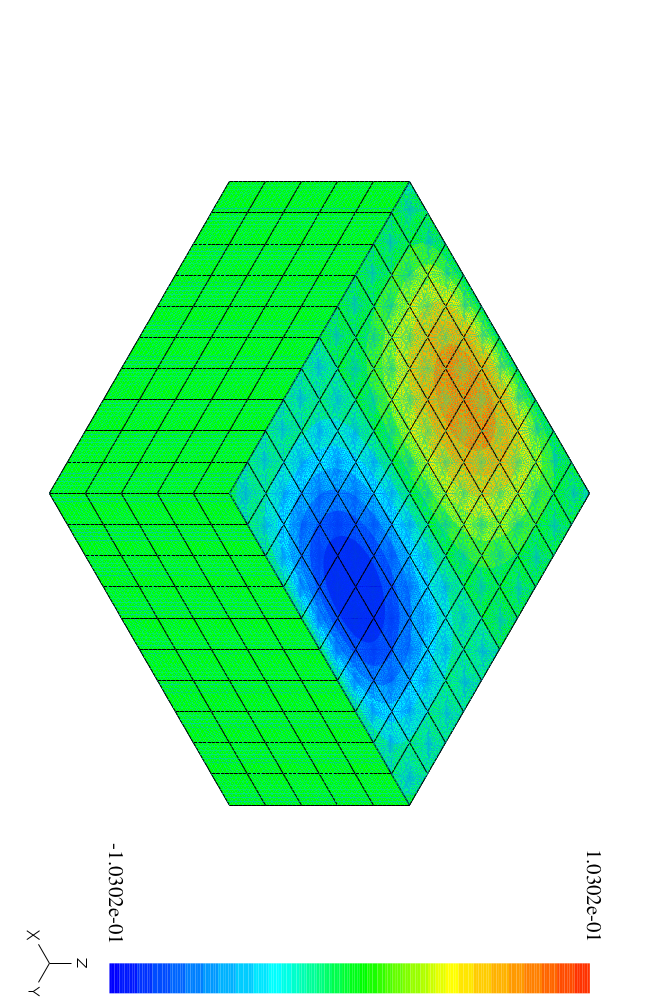}  & 
		 \includegraphics[height=7cm, angle=90]{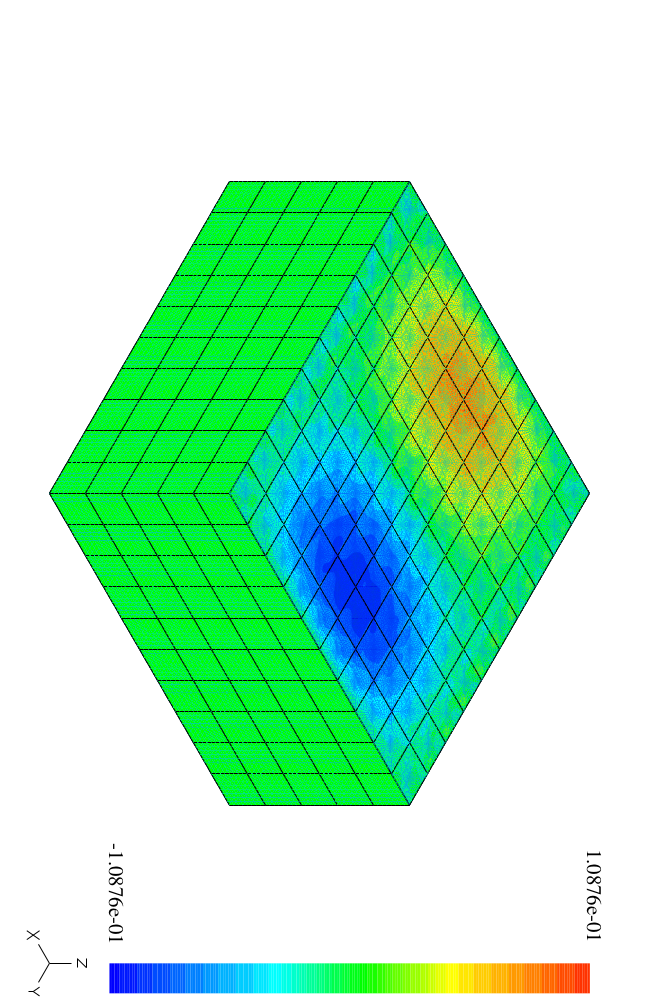}\\
		 (a) & (e) \\[-0.25ex] \hline
		 \includegraphics[height=7cm, angle=90]{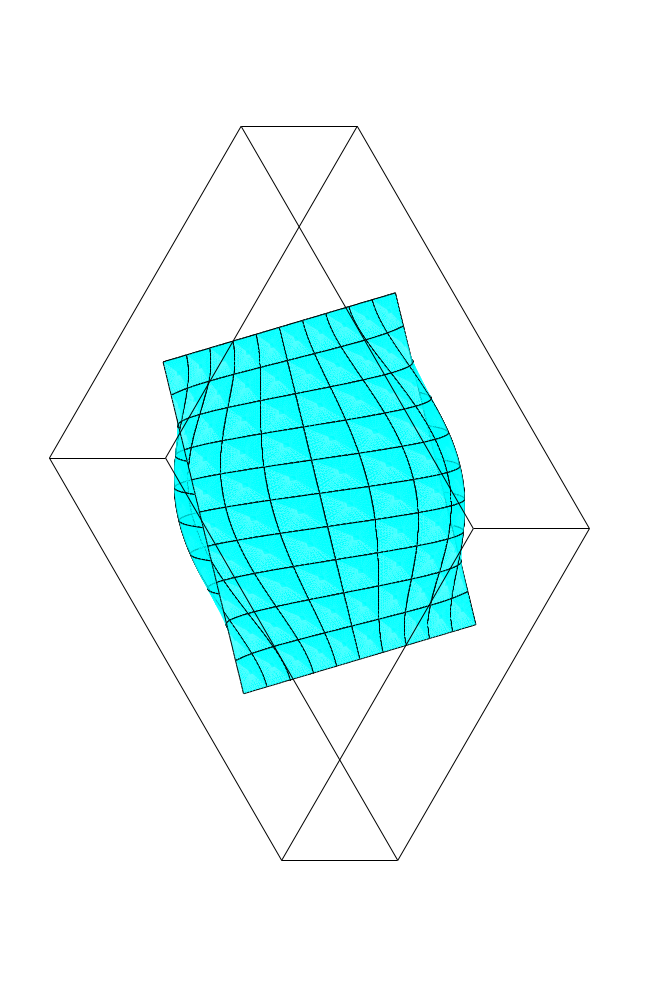}  & 
		 \includegraphics[height=7cm, angle=90]{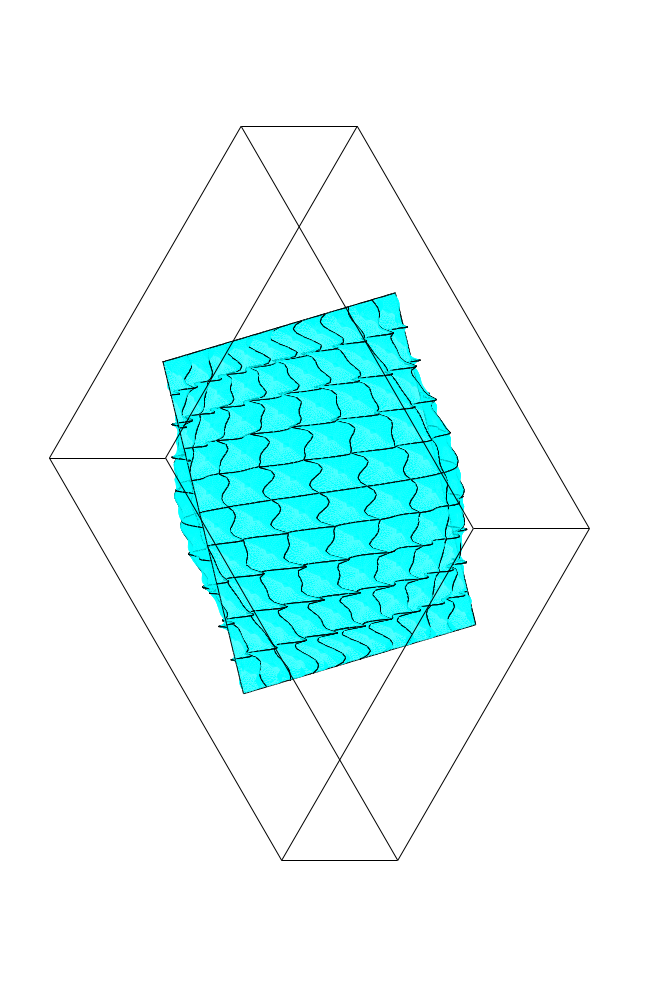}\\
		 (b) & (f) \\[-0.25ex] \hline
		 \includegraphics[height=7cm, angle=90]{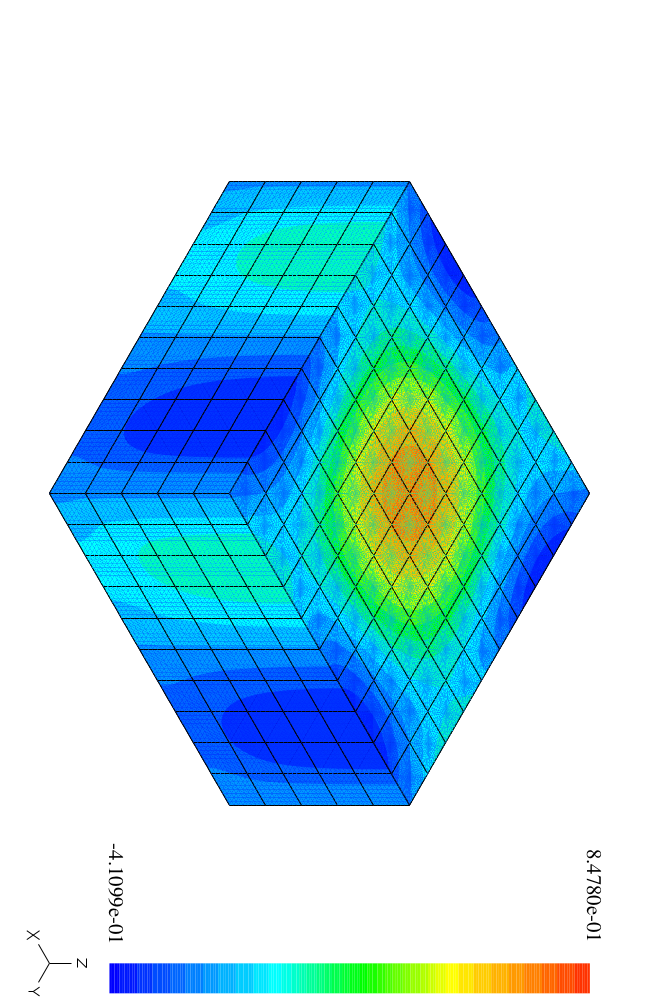}  & 
		 \includegraphics[height=7cm, angle=90]{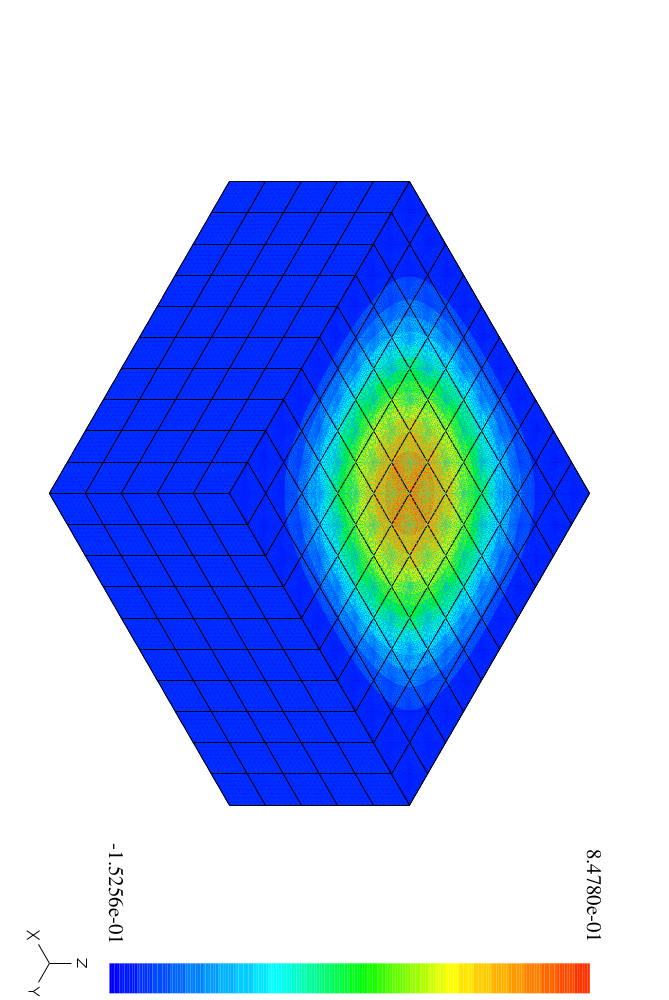} \\
		 (c) & (g) \\[-0.25ex] \hline
		 \includegraphics[height=7cm, angle=90]{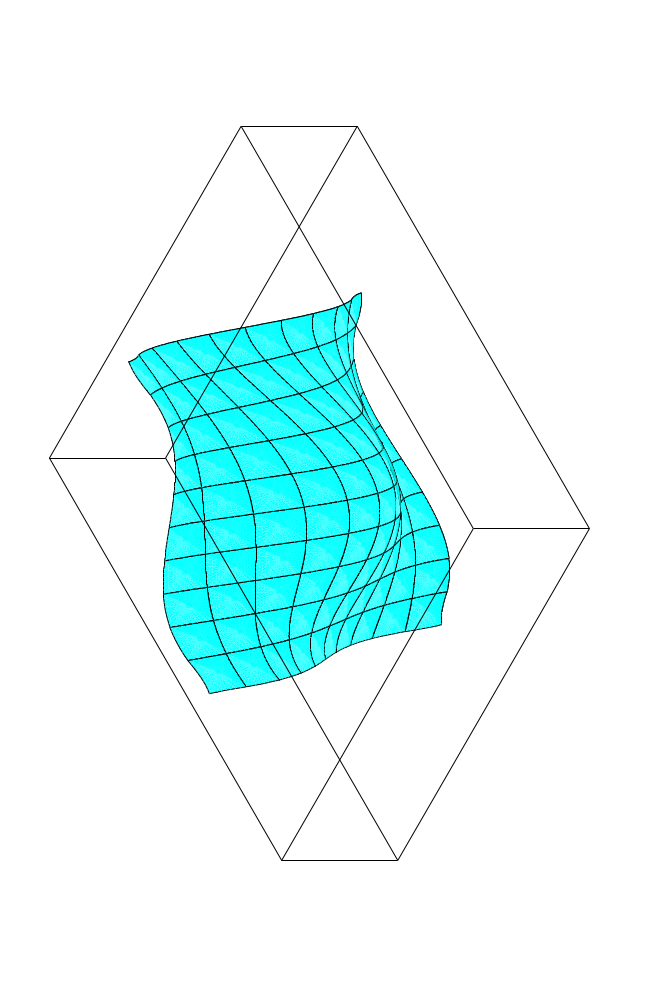}  & 
		 \includegraphics[height=7cm, angle=90]{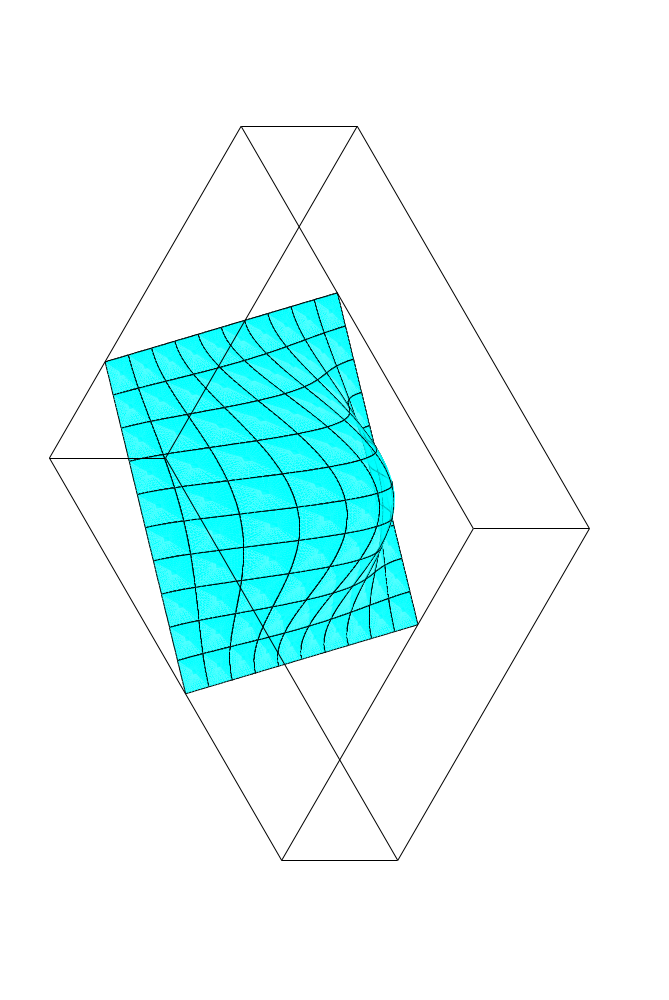}\\
		 (d) & (h) \\[-0.25ex] \hline
		\end{tabular}
		\caption{\scriptsize Cube $[-1,+1]^{3}$ with \emph{clamped boundary conditions} subject to the \emph{total body force}  \eqref{eq:Body force with a solenoidal and a rotation-free part}. The left column shows results for $\mu=1$  while results for $\mu=10^{-4}$ are shown in the right column.
		All results are drawn on the mid plane $z=0$. The solution values in the three-dimensional vector graphs (b),(d) (f) and (h) are scaled with the factor 0.25. The loss of pressure robustness computing the displacement solution is illustrated comparing Figures (a) and (b) with Figures (e) and (f), respectively.
		\newline 
		(a) iso-level map, displacement component $u_{x}$, $\min{(u_{x})}=-1.030\cdot10^{-1}$ and $\max{(u_{x})}=1.030\cdot10^{-1}$. \\
		(b) 3D three-dimensional vector graph, displacement component $u_{x}$,\\
		(c) iso-level map, pressure reaction $p$,$\min{(p)}=-4.109\cdot10^{-1}$ and $\max{(p)}=8.478\cdot10^{-1}$,\\
		(d) 3D three-dimensional vector graph, pressure reaction $p$,\\
		(e) iso-level map, displacement component $u_{x}$, $\min{(u_{x})}=-1.088\cdot10^{-1}$ and $\max{(u_{x})}=1.088\cdot10^{-1}$,\\
		(f) 3D three-dimensional vector graph, displacement component $u_{x}$,\\
		(g) iso-level map, pressure reaction $p$, $\min{(p)}=-1.526\cdot10^{-2}$ and $\max{(p)}=8.478\cdot10^{-1}$.\\
		(h) 3D three-dimensional vector graph, pressure reaction $p$.\\
					\label{tab:total-body-force}}
	\end{table}

	\begin{table}[tbh]
		\begin{tabular}{|c|c|}\hline
		 \includegraphics[height=7cm, angle=90]{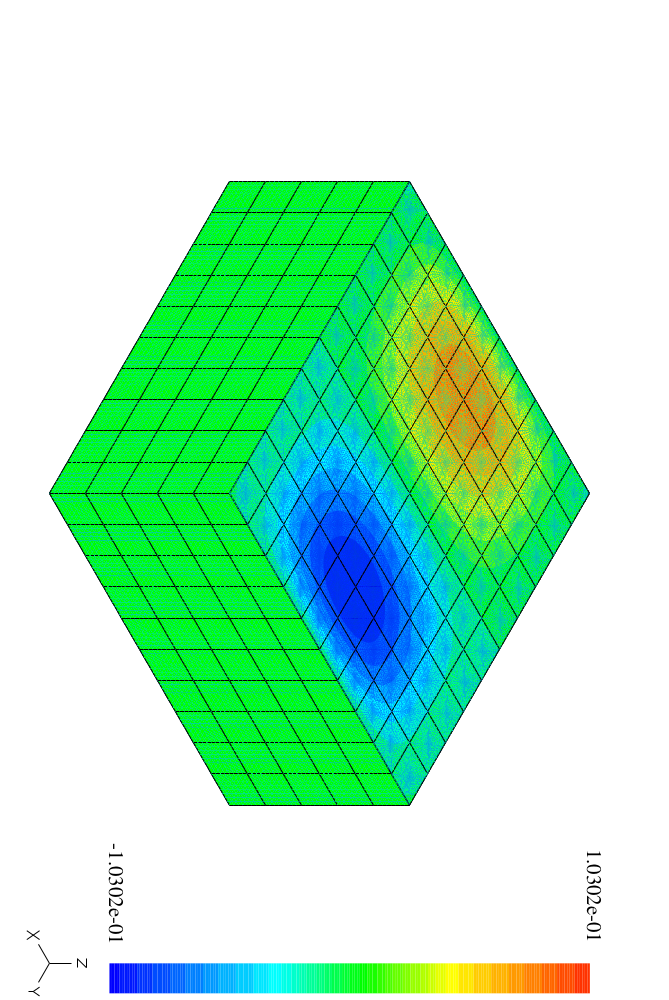}  & 
		 \includegraphics[height=7cm, angle=90]{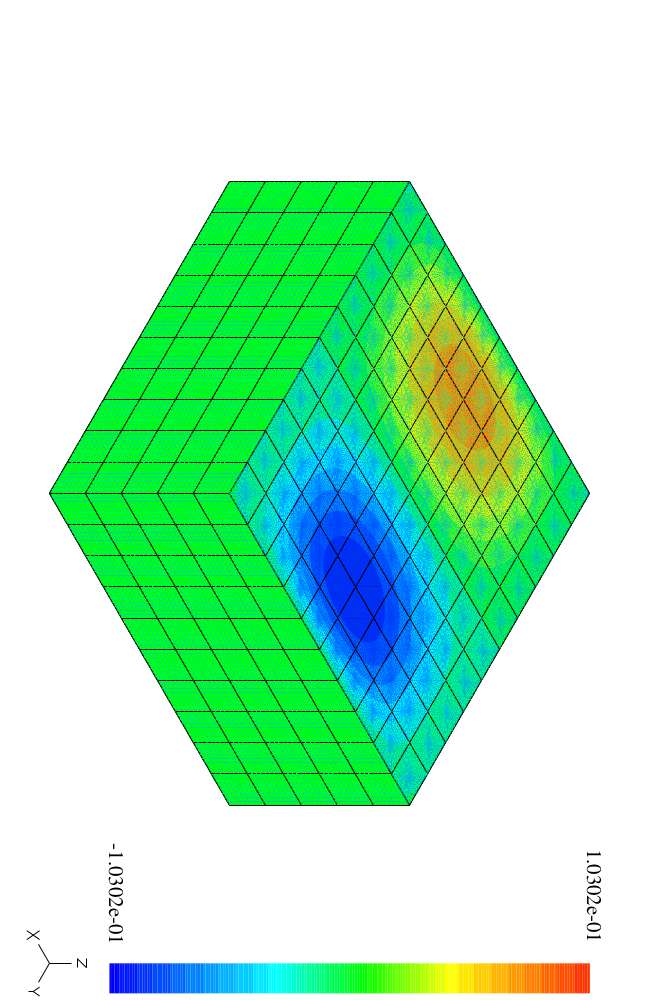}\\
		 (a) & (e) \\[-0.25ex] \hline
		 \includegraphics[height=7cm, angle=90]{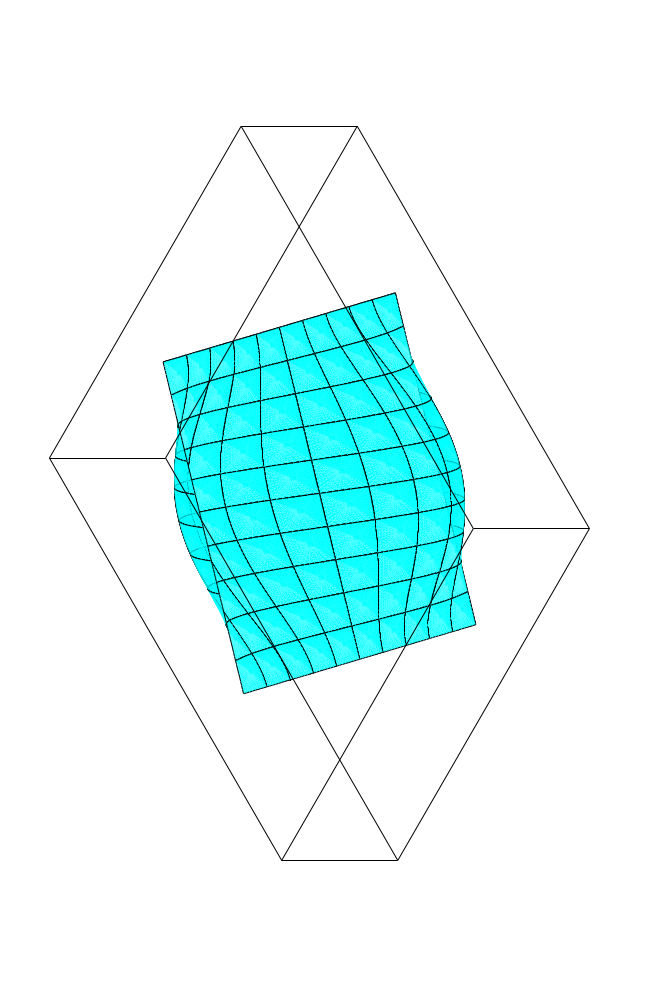}  & 
		 \includegraphics[height=7cm, angle=90]{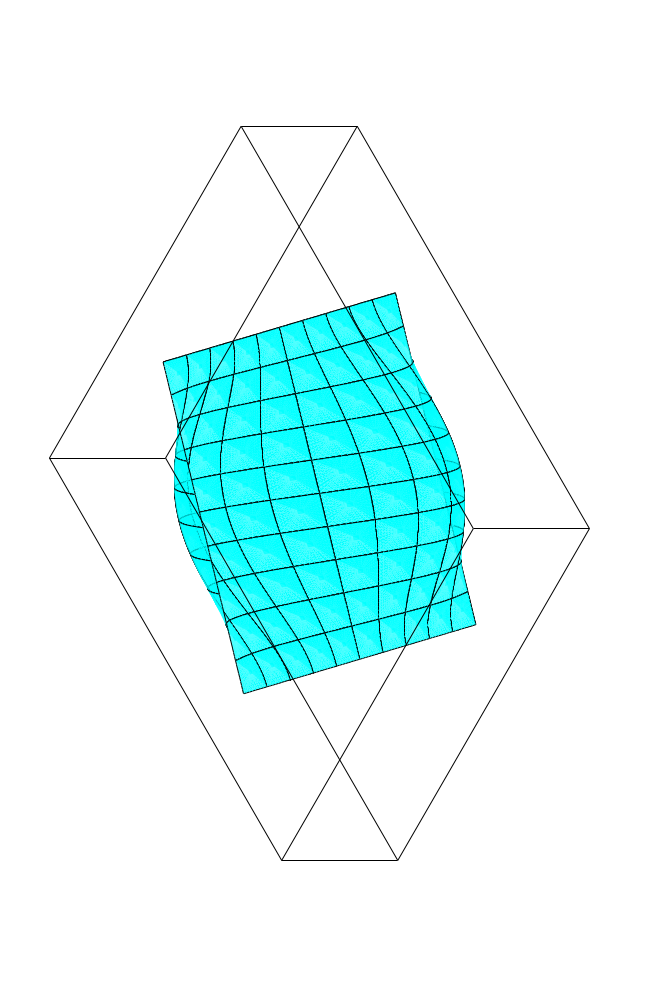}\\
		 (b) & (f) \\[-0.25ex] \hline
		 \includegraphics[height=7cm, angle=90]{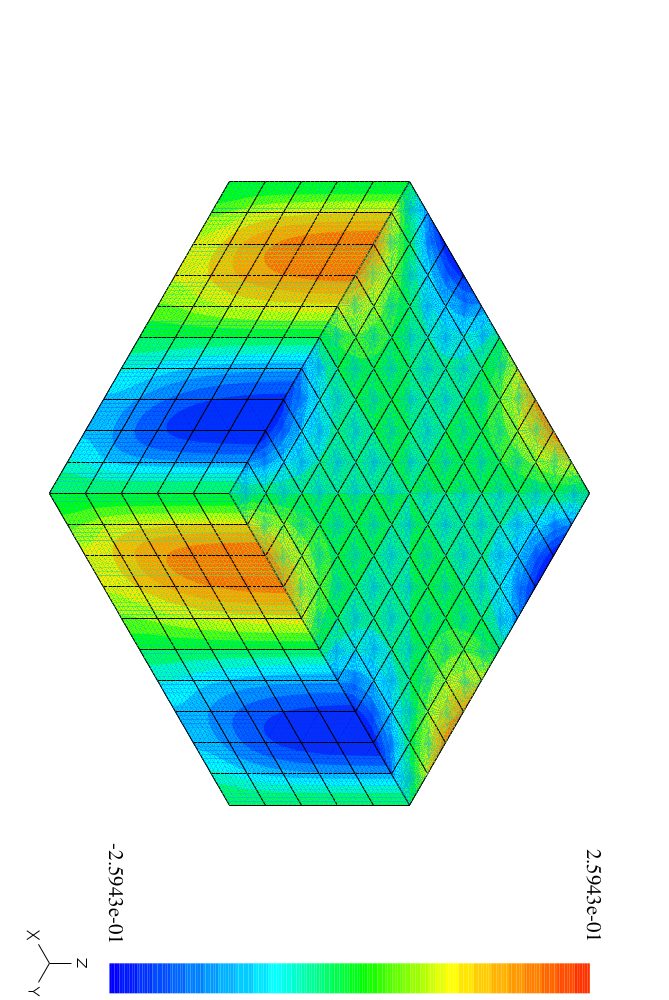}  & 
		 \includegraphics[height=7cm, angle=90]{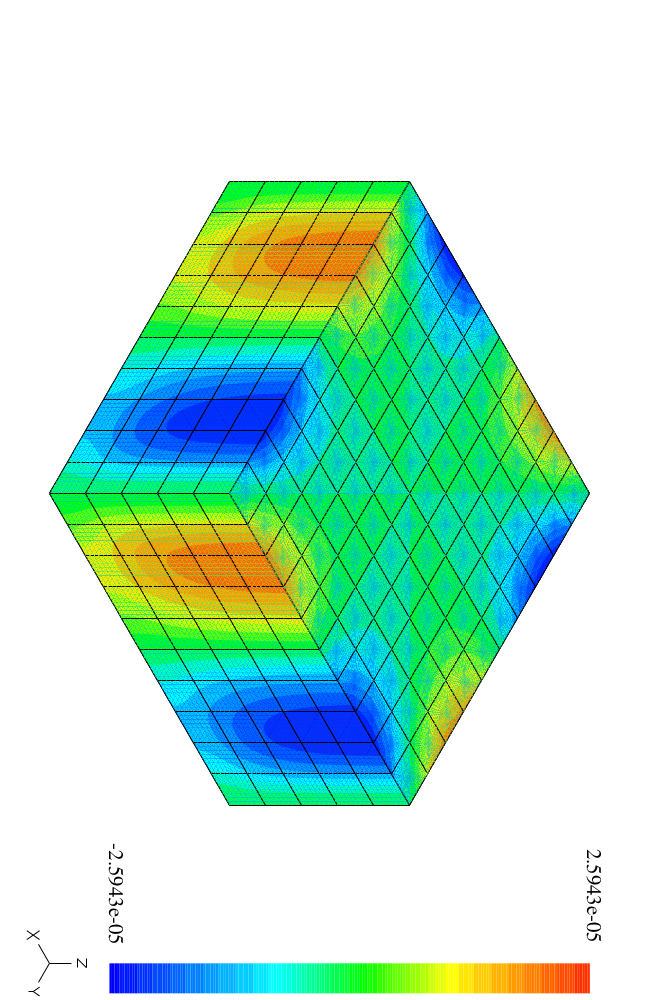} \\
		 (c) & (g) \\[-0.25ex] \hline
		 \includegraphics[height=7cm, angle=90]{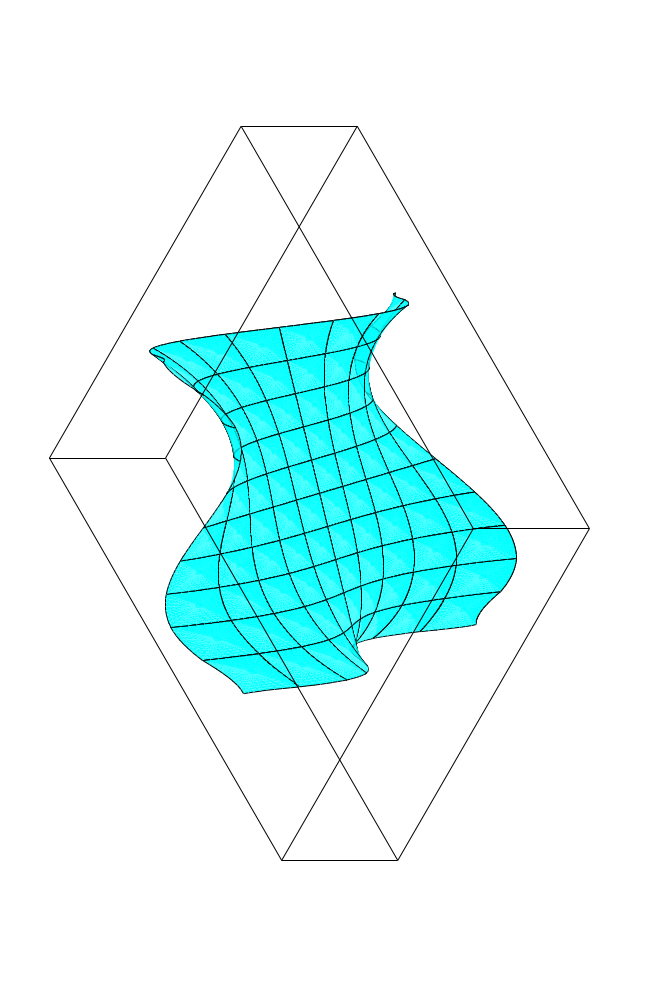}  & 
		 \includegraphics[height=7cm, angle=90]{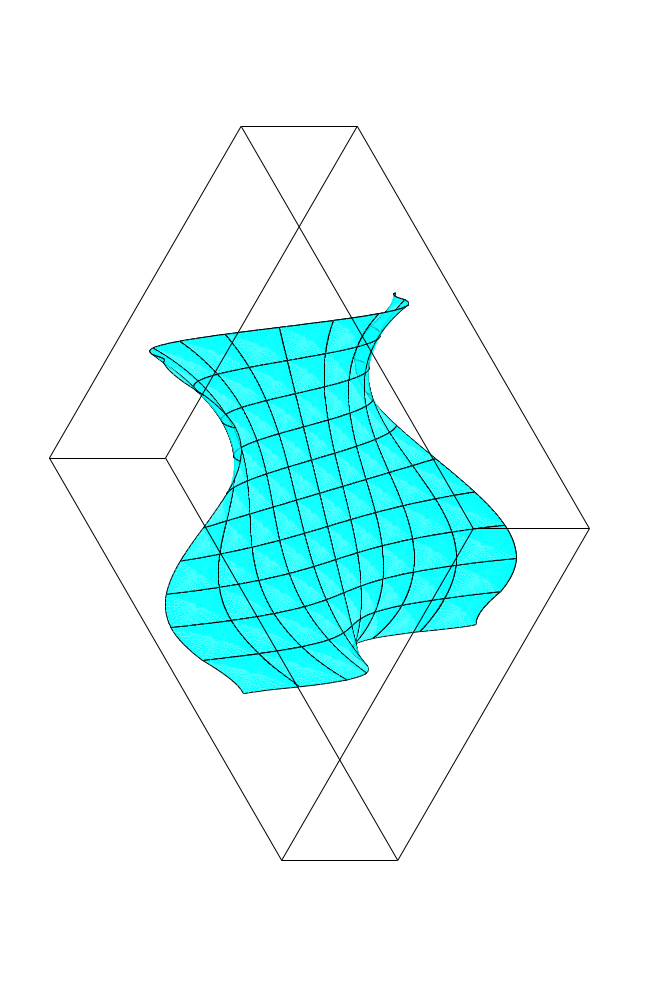}\\
		 (d) & (h) \\[-0.25ex] \hline
		\end{tabular}
		\caption{ \scriptsize Cube $[-1,+1]^{3}$ with \emph{clamped boundary conditions} subject to the to the the numerically determined  \emph{divergence-free body force}, $\curl\vek{A}_{h}$ \eqref{eq:closed-form-divergence-free-body-force}, Step1. The left column shows results for $\mu=1$ while results for $\mu=10^{-4}$  are shown in the right column.
		All results are drawn on the mid plane $z=0$. The solution values in the three-dimensional vector graphs (b),(d) (f) and (h) are scaled with the factor 0.25.\newline 
		(a) iso-level map, displacement component $u_{1\,x}$, $\min{(u_{1\,x})}=-1.0302\cdot10^{-1}$ and $\max{(u_{1\,x})}=1.0302\cdot10^{-1}$. \\
		(b) 3D three-dimensional vector graph, displacement component $u_{1\,x}$,\\
		(c) iso-level map, pressure reaction $p_{1}$,	$\min{(p_{1})}=-2.5943\cdot10^{-1}$ and $\max{(p_{1})}=2.5943\cdot10^{-1}$,\\
		(d) 3D three-dimensional vector graph, pressure reaction $p_{1}$,\\
		(e) iso-level map, displacement component $u_{1\,x}$, $\min{(u_{x})}=-1.030\cdot10^{-1}$ and $\max{(u_{1\,x})}=1.030\cdot10^{-1}$,\\
		(f) 3D three-dimensional vector graph, displacement component $u_{1\,x}$,\\
		(g) iso-level map, pressure reaction $p_{1}$, $\min{(p_{1})}=-2.5943\cdot10^{-5}$ and $\max{(p_{1})}=2.5943\cdot10^{-5}$,\\
		(h) 3D three-dimensional vector graph, pressure reaction $p_{1}$.
			\label{tab:solenoidal-body-force}}
	\end{table}

	\begin{table}[tbh]
		\begin{tabular}{|c|c|}\hline
		\includegraphics[height=7cm, angle=90]{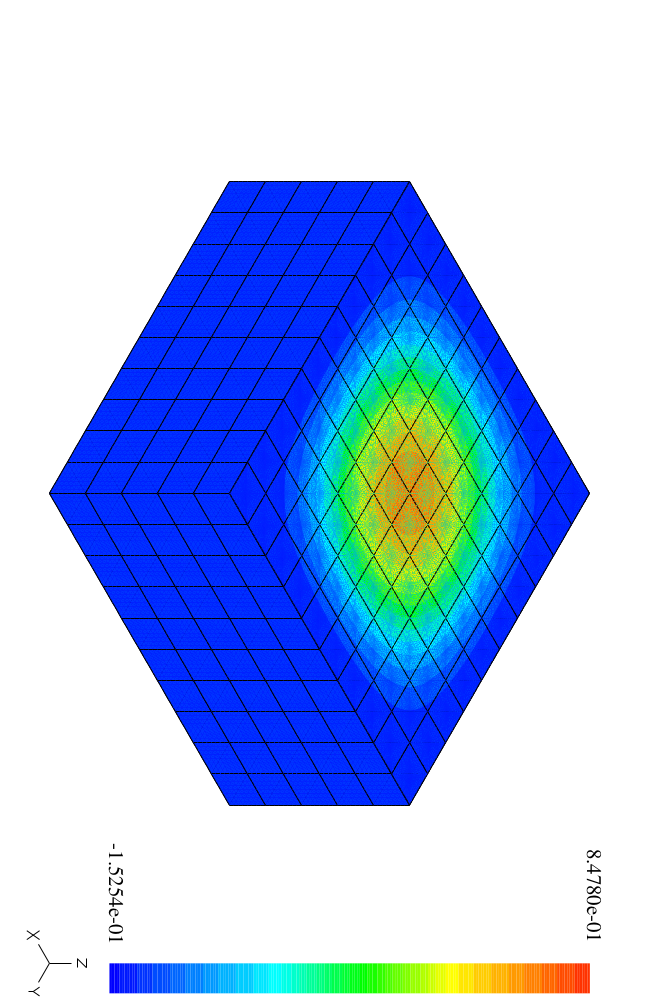} &
	    \includegraphics[height=7cm, angle=90]{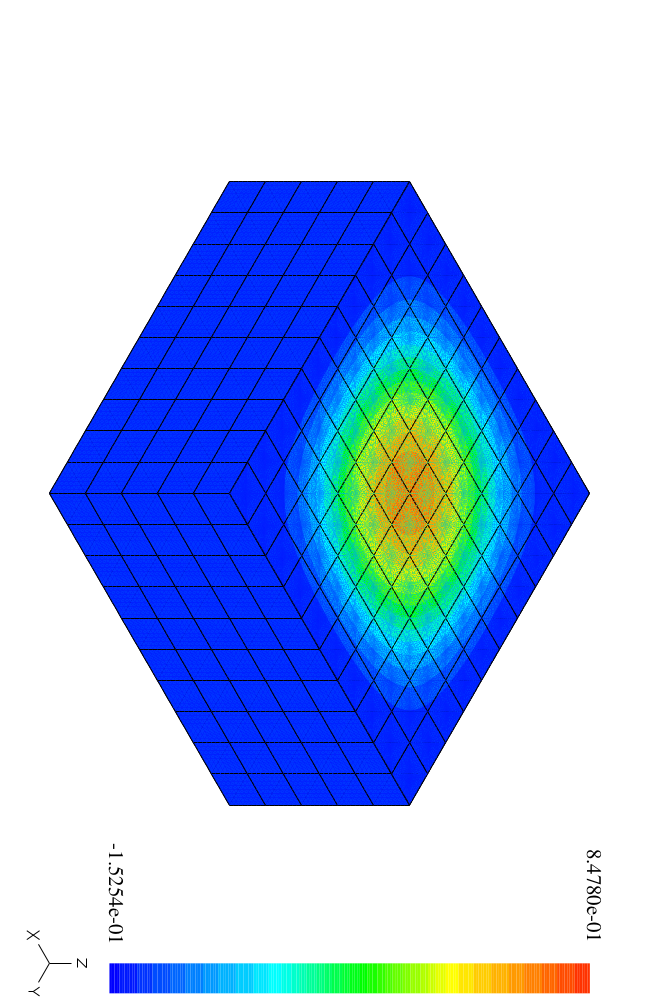}  \\
		(a) & (c) \\[-0.25ex] \hline
		\includegraphics[height=7.4cm, angle=90]{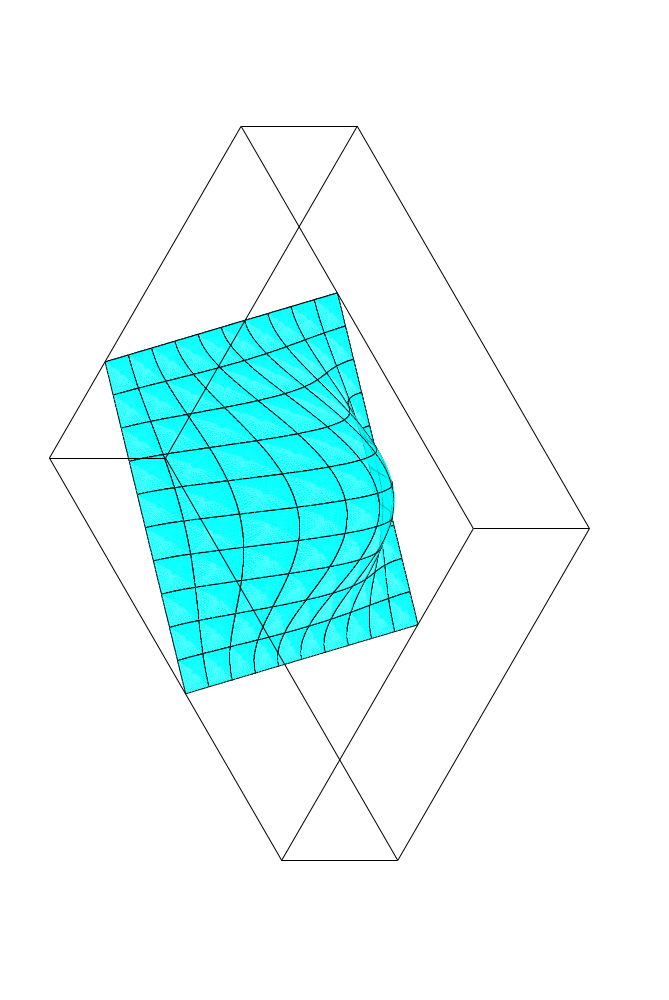}&
		\includegraphics[height=7.4cm, angle=90]{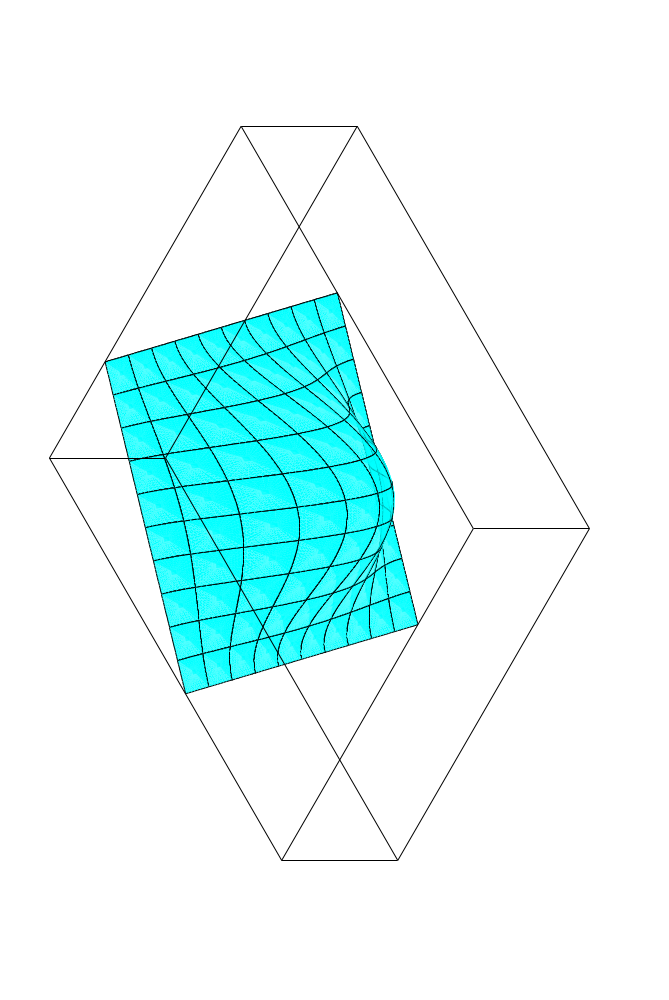} \\
		(b) & (d) \\[-0.25ex] \hline
		\end{tabular}
		\caption{Cube $[-1,+1]^{3}$ with \emph{clamped boundary conditions} subject to the numerically determined \emph{rotation-free body force}, \mbox{$\vek{f}-\curl\vek{A}_{h}$} \eqref{eq:closed-form-rotation-free-body-force}, Step2. 
		The left column shows results for $\mu=1$ while results for $\mu=10^{-4}$ are shown in the right column. All results are drawn on the mid plane $z=0$. The solution values in the three-dimensional vector graphs (b) and(d) are scaled with the factor 0.25.\newline 
		The pressure reaction $p_{2}$, $\min{(p_{2})}=-1.525\cdot10^{-1}$ and $\max{(p_{2})}=8.478\cdot10^{-1}$ is the same for both cases $\mu=1$ and $\mu=10^{-4}$;\newline
		(a) iso-level map, pressure reaction $p_{2}$, \\
		(b) three-dimensional vector graph, pressure reaction $p_{2}$,\\ 
		(c) iso-level map, pressure reaction $p_{2}$,\\
		(d) three-dimensional vector graph, pressure reaction $p_{2}$.\\
		\label{tab:rotation-free-body-force}}
    \end{table}

    	\begin{table}[tbh]
		\begin{tabular}{|c|c|}\hline
		 \includegraphics[height=7cm, angle=90]{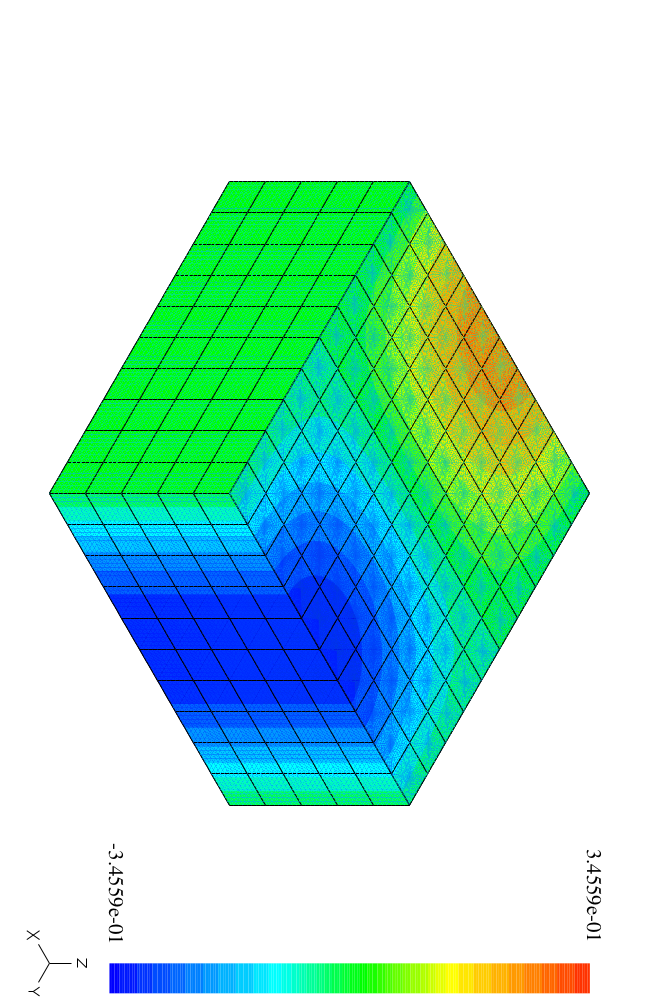}  & 
		 \includegraphics[height=7cm, angle=90]{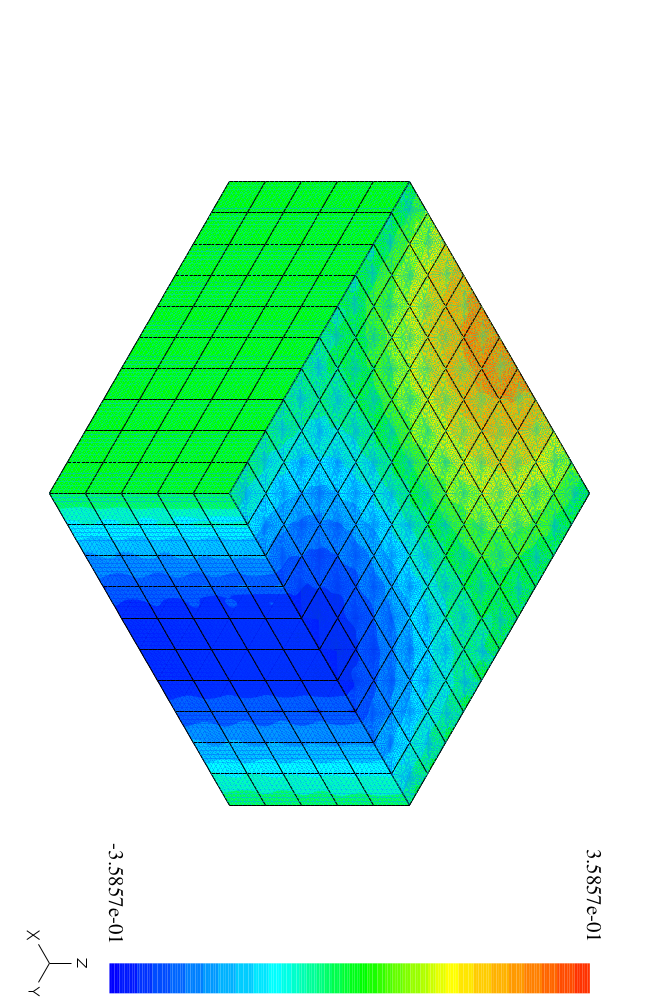}\\
		 (a) & (e) \\[-0.25ex] \hline
		 \includegraphics[height=7cm, angle=90]{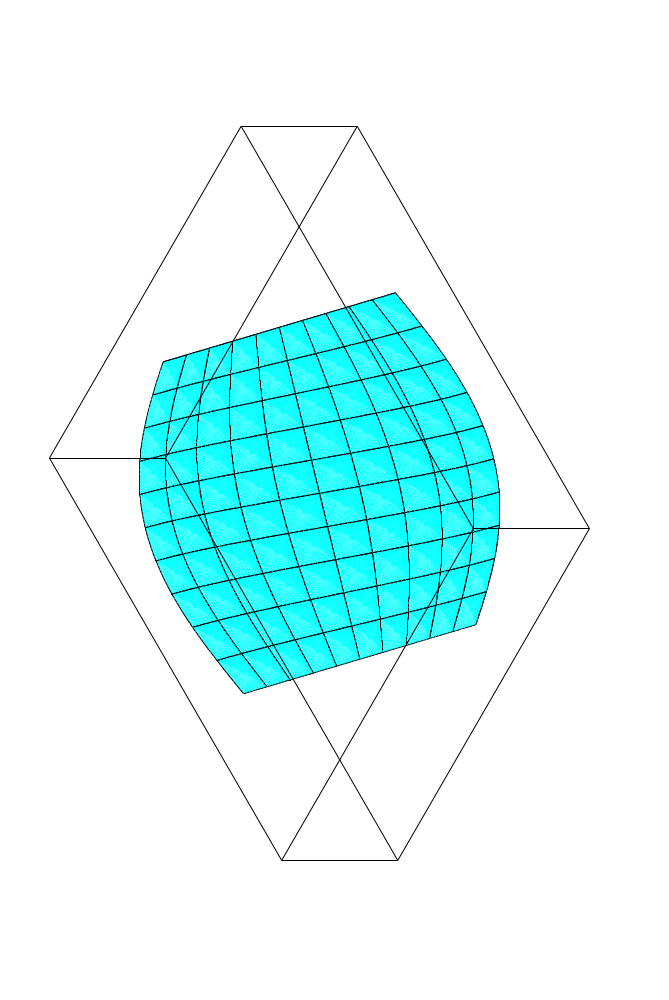}  & 
		 \includegraphics[height=7cm, angle=90]{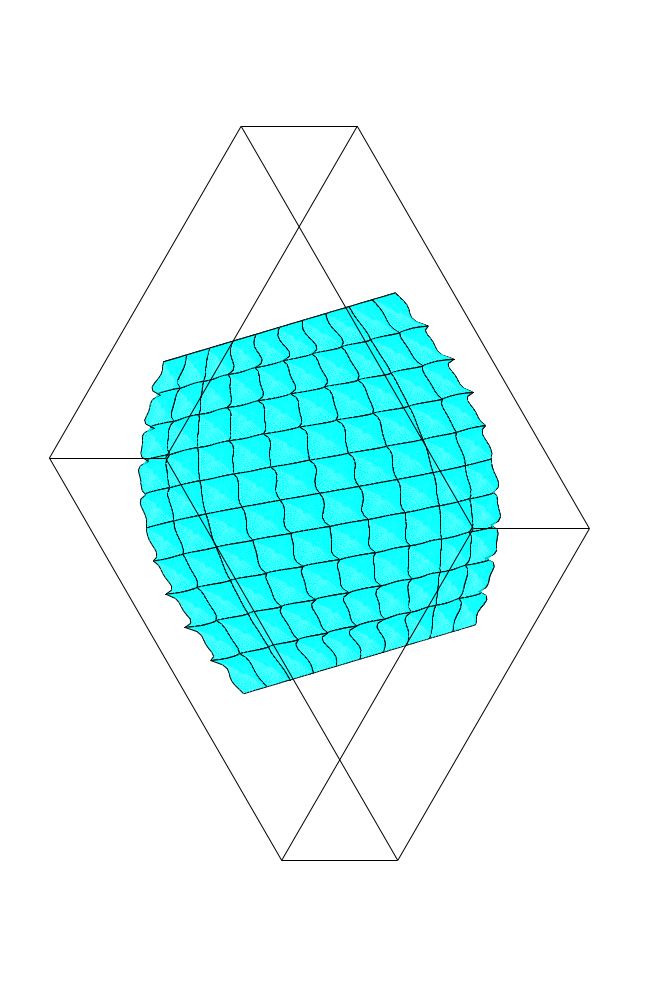}\\
		 (b) & (f) \\[-0.25ex] \hline
		 \includegraphics[height=7cm, angle=90]{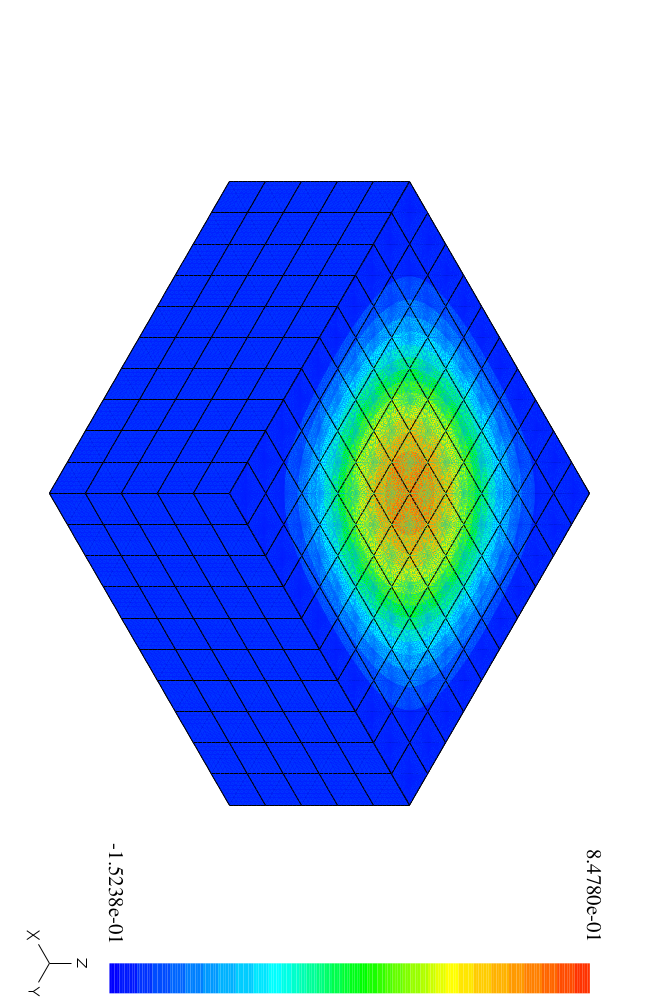}  & 
		 \includegraphics[height=7cm, angle=90]{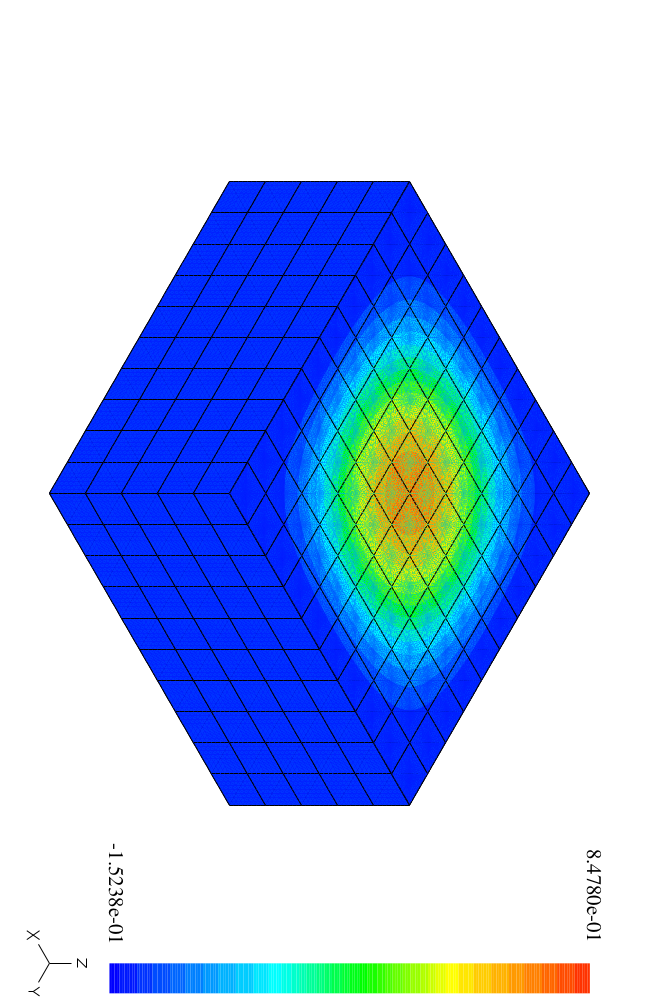} \\
		 (c) & (g) \\[-0.25ex] \hline
		 \includegraphics[height=7cm, angle=90]{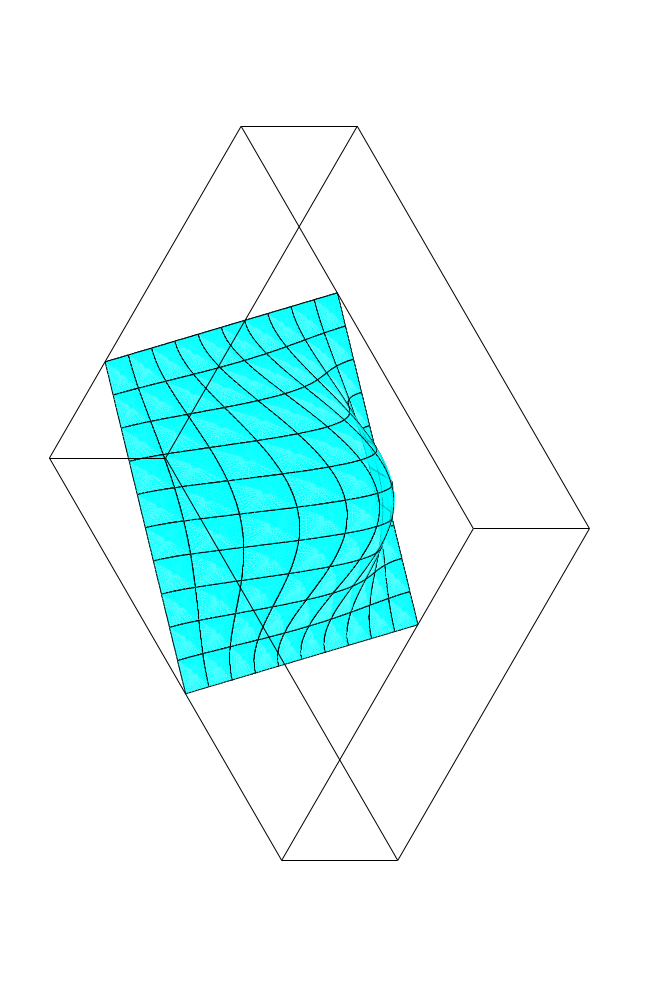}  & 
		 \includegraphics[height=7cm, angle=90]{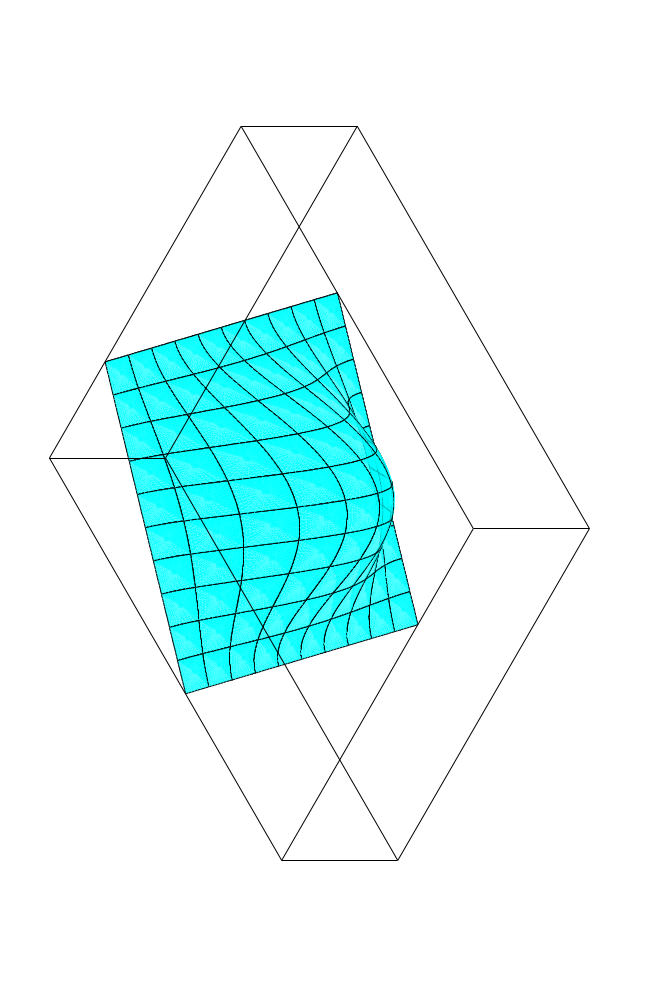}\\
		 (d) & (h) \\[-0.25ex] \hline
		\end{tabular}
		\caption{\scriptsize Cube $[-1,+1]^{3}$ with \emph{no-penetration boundary conditions}, subject to the \emph{total body force}  \eqref{eq:Body force with a solenoidal and a rotation-free part}. The left column shows results for $\mu=1$  while results for $\mu=10^{-4}$ are shown in the right column.
		All results are drawn on the mid plane $z=0$. The solution values in the three-dimensional vector graphs (b),(d) (f) and (h) are scaled with the factor 0.25. The loss of pressure robustness computing the displacement solution is illustrated comparing Figures (a) and (b) with Figures (e) and (f), respectively.
		\newline 
		(a) iso-level map, displacement component $u_{x}$, $\min{(u_{x})}=-3.456\cdot10^{-1}$ and $\max{(u_{x})}=3.456\cdot10^{-1}$. \\
		(b) 3D three-dimensional vector graph, displacement component $u_{x}$,\\
		(c) iso-level map, pressure reaction $p$,$\min{(p)}=-1.524\cdot10^{-1}$ and $\max{(p)}=8.478\cdot10^{-1}$,\\
		(d) 3D three-dimensional vector graph, pressure reaction $p$,\\
		(e) iso-level map, displacement component $u_{x}$, $\min{(u_{x})}=-3.586\cdot10^{-1}$ and $\max{(u_{x})}=3.586\cdot10^{-1}$,\\
		(f) 3D three-dimensional vector graph, displacement component $u_{x}$,\\
		(g) iso-level map, pressure reaction $p$, $\min{(p)}=-1.524\cdot10^{-2}$ and $\max{(p)}=8.478\cdot10^{-1}$.\\
		(h) 3D three-dimensional vector graph, pressure reaction $p$.\\
			\label{tab:total-body-force-no-penetration}}
	\end{table}	
	
	\begin{table}[tbh]
		\begin{tabular}{|c|c|}\hline
		 \includegraphics[height=7cm, angle=90]{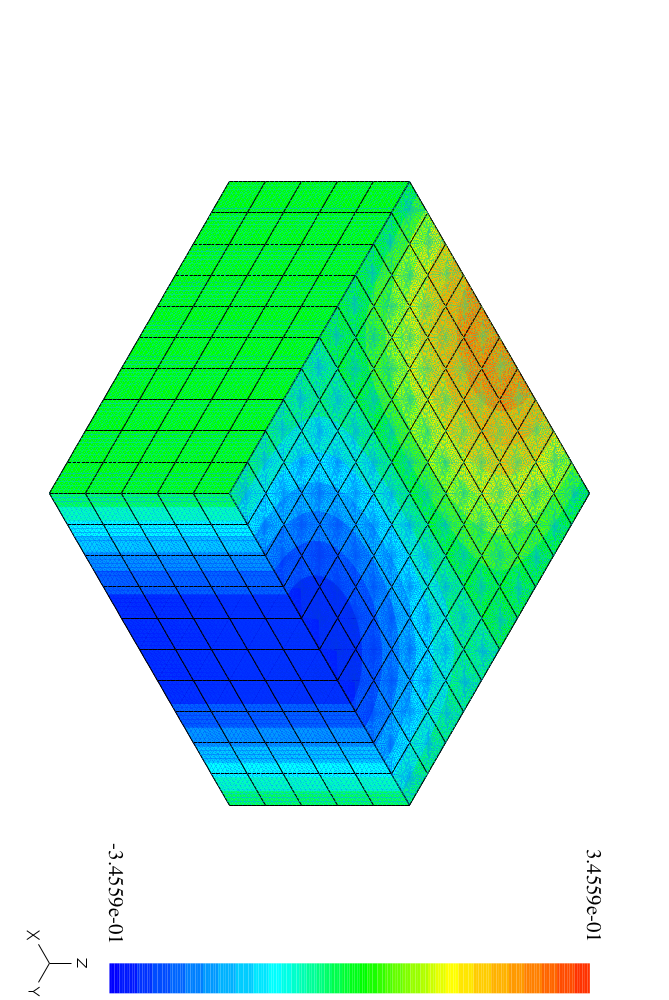}  & 
		 \includegraphics[height=7cm, angle=90]{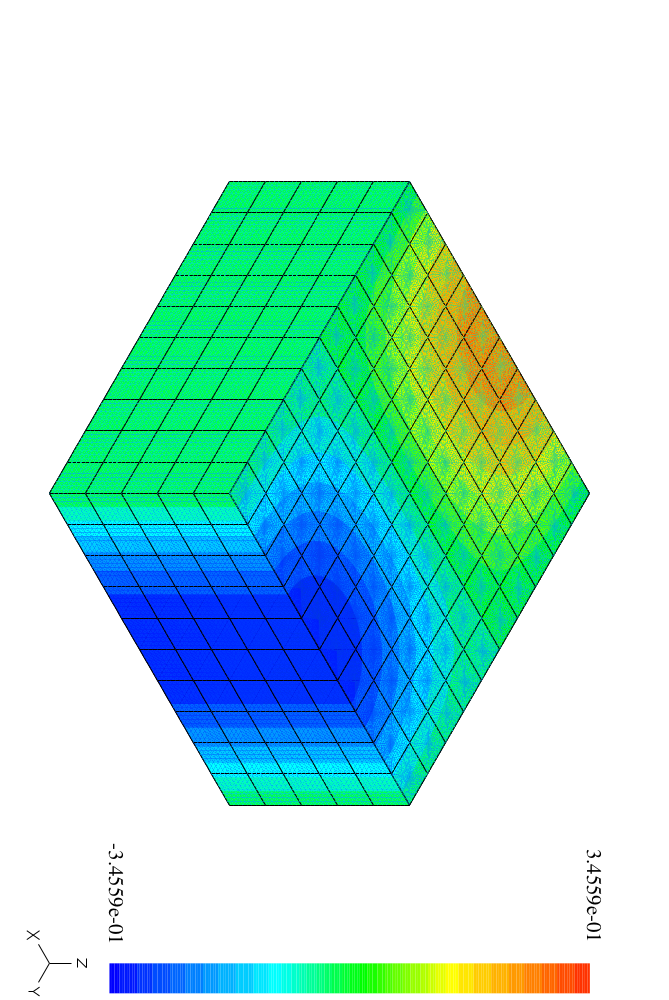}\\
		 (a) & (e) \\[-0.25ex] \hline
		 \includegraphics[height=7cm, angle=90]{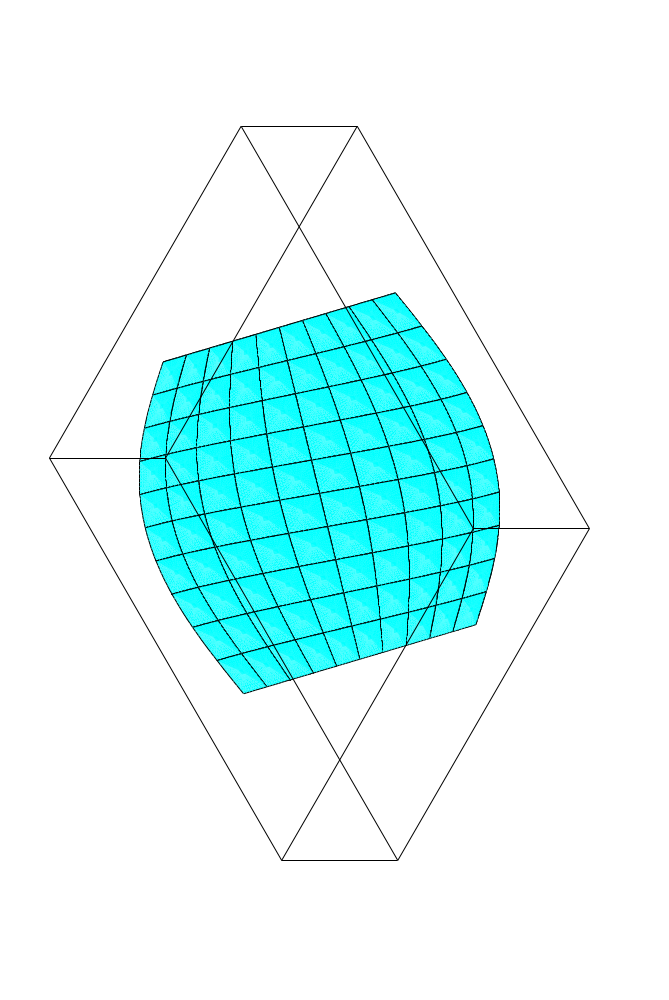}  & 
		 \includegraphics[height=7cm, angle=90]{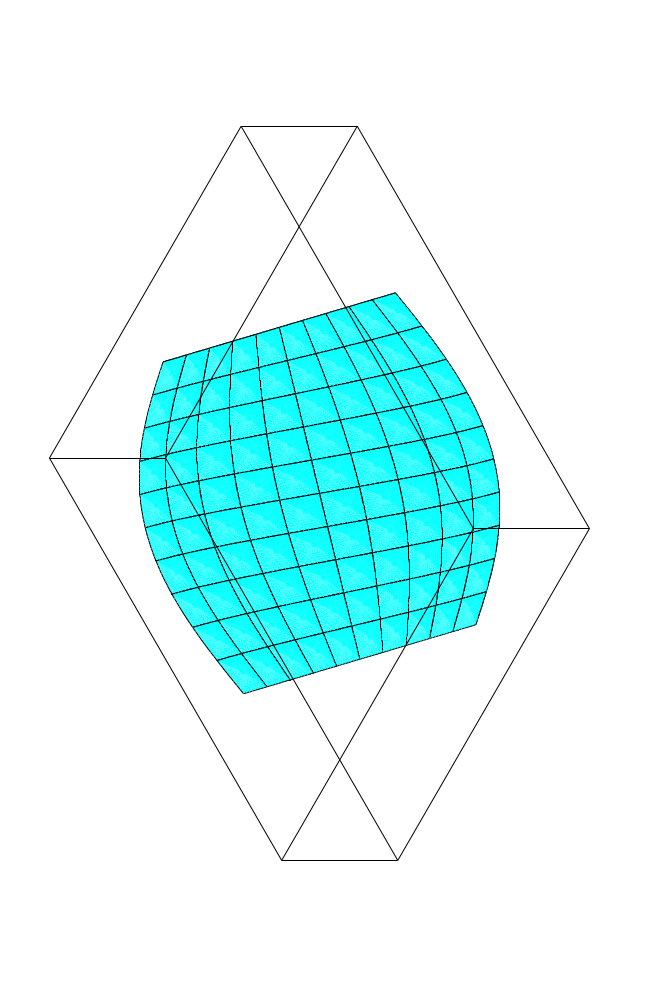}\\
		 (b) & (f) \\[-0.25ex] \hline
		 \includegraphics[height=7cm, angle=90]{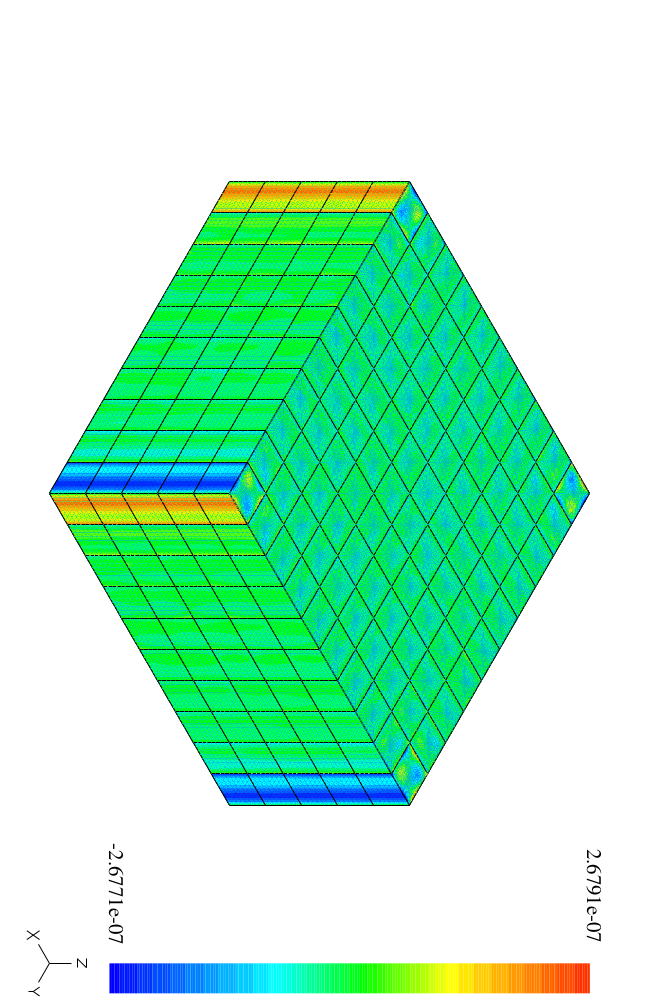}  & 
		 \includegraphics[height=7cm, angle=90]{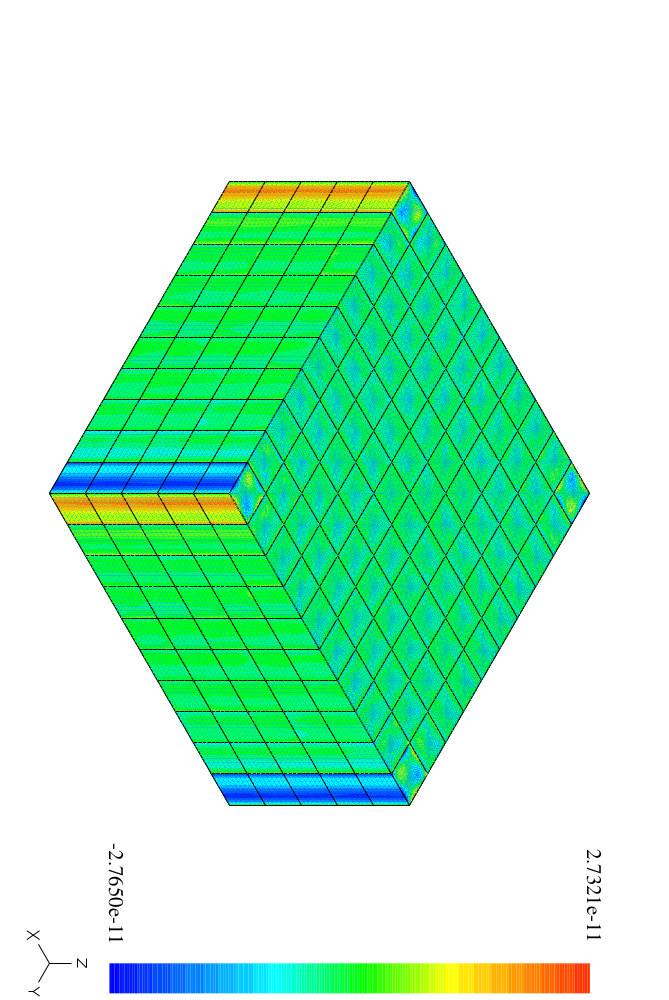} \\
		 (c) & (g) \\[-0.25ex] \hline
		 \includegraphics[height=7cm, angle=90]{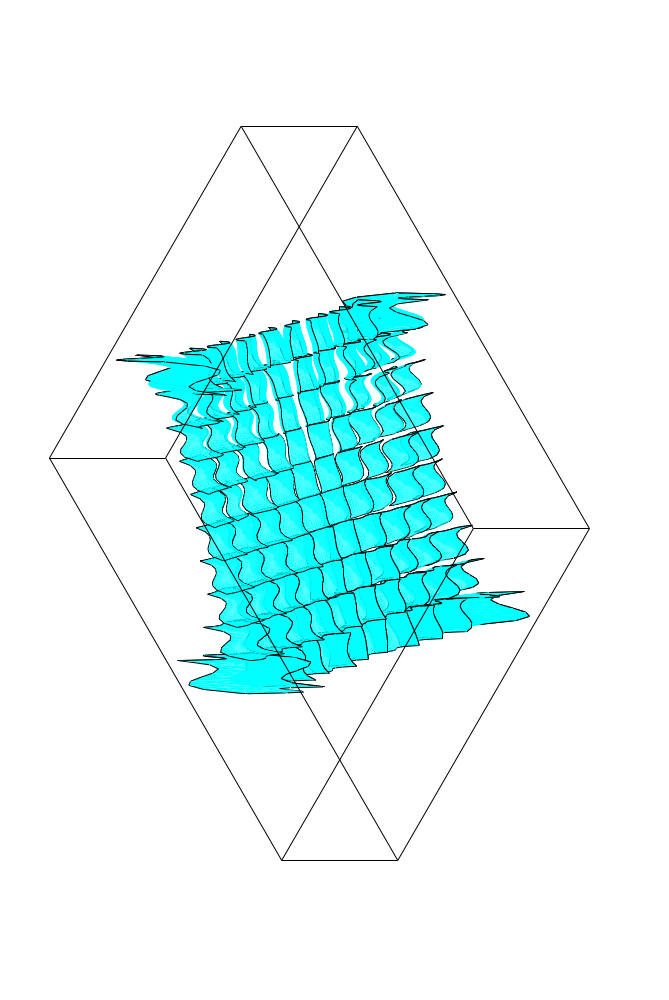}  & 
		 \includegraphics[height=7cm, angle=90]{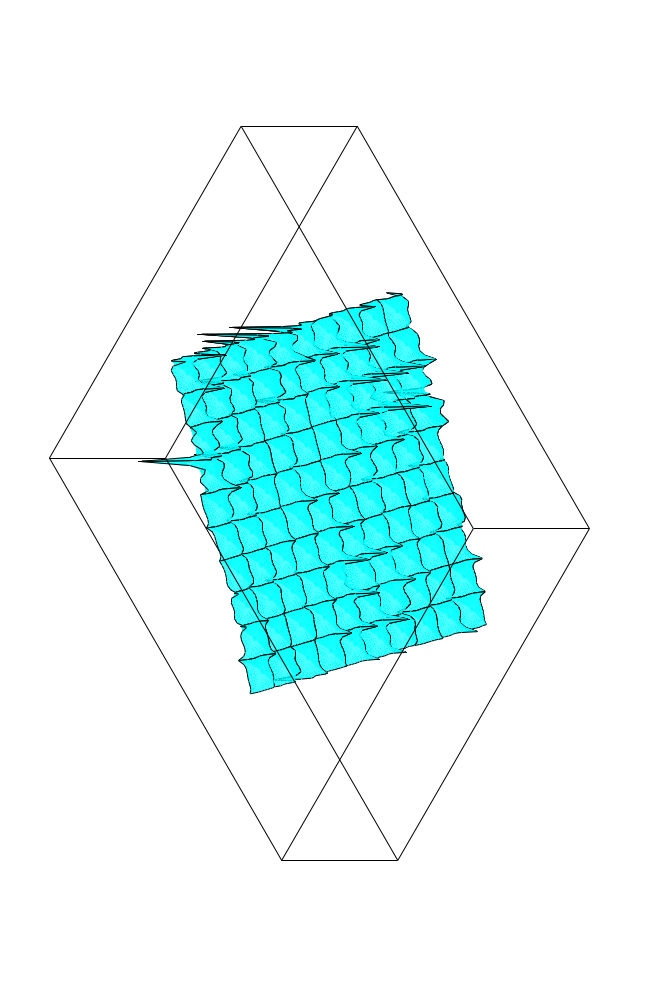}\\
		 (d) & (h) \\[-0.25ex] \hline
		\end{tabular}
		\caption{ \scriptsize Cube $[-1,+1]^{3}$ with \emph{no-penetration boundary conditions} subject to the to the the numerically determined  \emph{divergence-free body force}, $\curl\vek{A}_{h}$ \eqref{eq:closed-form-divergence-free-body-force}, Step1. The left column shows results for $\mu=1$ while results for $\mu=10^{-4}$  are shown in the right column.
		All results are drawn on the mid plane $z=0$. The solution values in the three-dimensional vector graphs (b),(d) (f) and (h) are scaled with the factor 0.25. The pressure reactions shown in Figures (c) and (d) for $\mu=1$ and in  Figures (g) and  (h) for $\mu=10^{-4}$ as numerically zero. They should be compared
		to the corresponding Figures (c) and (d), and (g) and (h) in Table~\ref{tab:solenoidal-body-force}, respectively.
		\newline 
		(a) iso-level map, displacement component $u_{1\,x}$, $\min{(u_{1\,x})}=-3.456\cdot10^{-1}$ and $\max{(u_{1\,x})}=3.456\cdot10^{-1}$. \\
		(b) 3D three-dimensional vector graph, displacement component $u_{1\,x}$,\\
		(c) iso-level map, pressure reaction $p_{1}$,	$\min{(p_{1})}=-2.667\cdot10^{-7}$ and $\max{(p_{1})}=2.679\cdot10^{-7}$,\\
		(d) 3D three-dimensional vector graph, pressure reaction $p_{1}$,\\
		(e) iso-level map, displacement component $u_{1\,x}$, $\min{(u_{x})}=-3.456\cdot10^{-1}$ and $\max{(u_{1\,x})}=3.456\cdot10^{-1}$,\\
		(f) 3D three-dimensional vector graph, displacement component $u_{1\,x}$,\\
		(g) iso-level map, pressure reaction $p_{1}$, $\min{(p_{1})}=-2.765\cdot10^{-11}$ and $\max{(p_{1})}=2.732\cdot10^{-11}$,\\
		(h) 3D three-dimensional vector graph, pressure reaction $p_{1}$.\\
					\label{tab:solenoidal-body-force-no-penetration}}
	\end{table}	
    
 	\begin{table}[tbh]
		\begin{tabular}{|c|c|}\hline
		\includegraphics[height=7cm, angle=90]{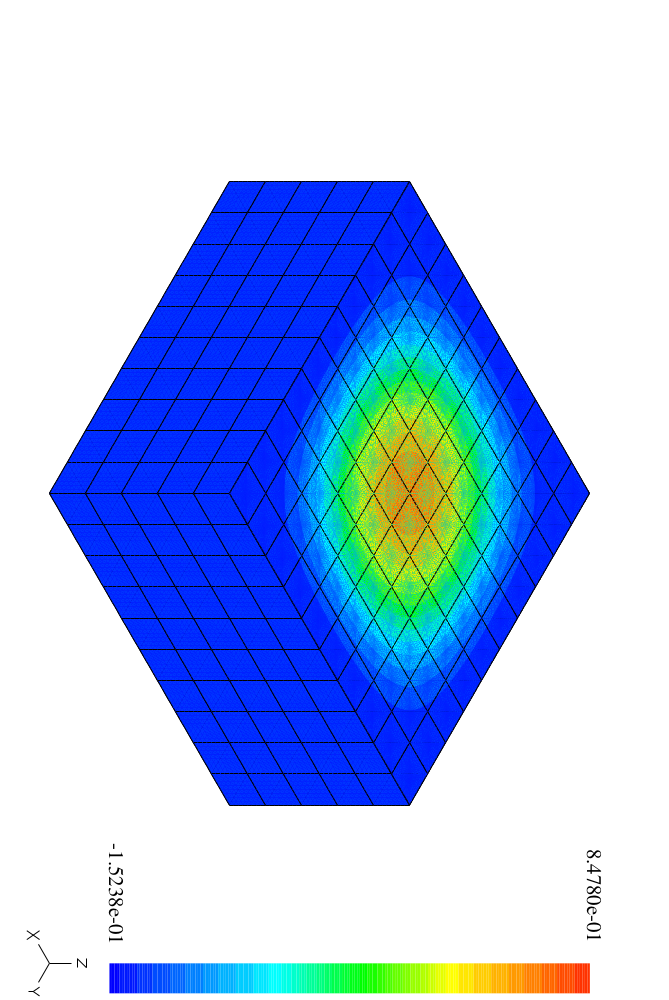} &
	    \includegraphics[height=7cm, angle=90]{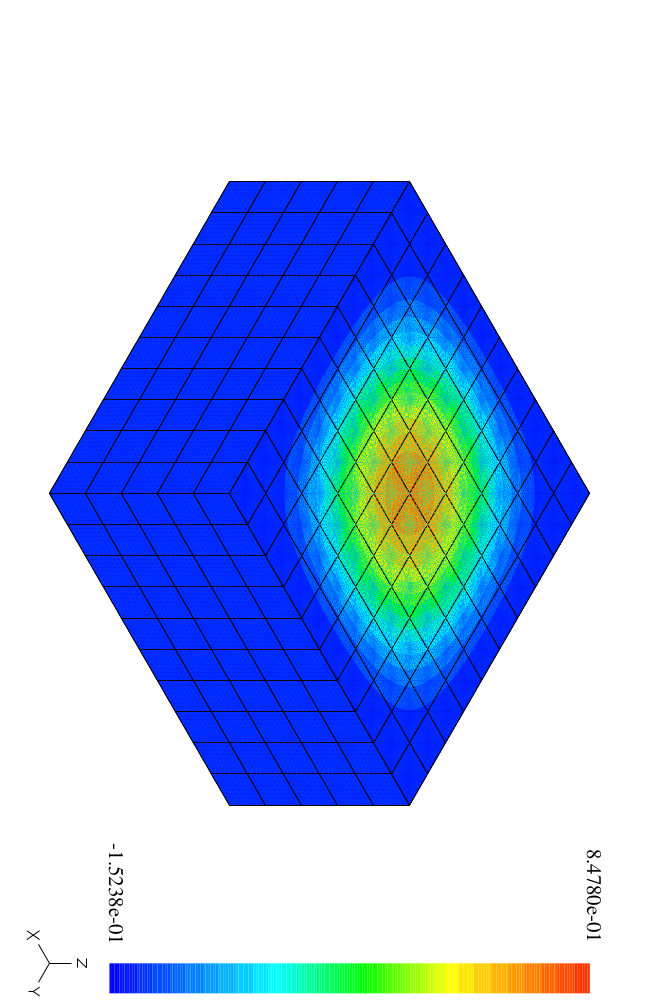}  \\
		(a) & (c) \\[-0.25ex] \hline
		\includegraphics[height=7.4cm, angle=90]{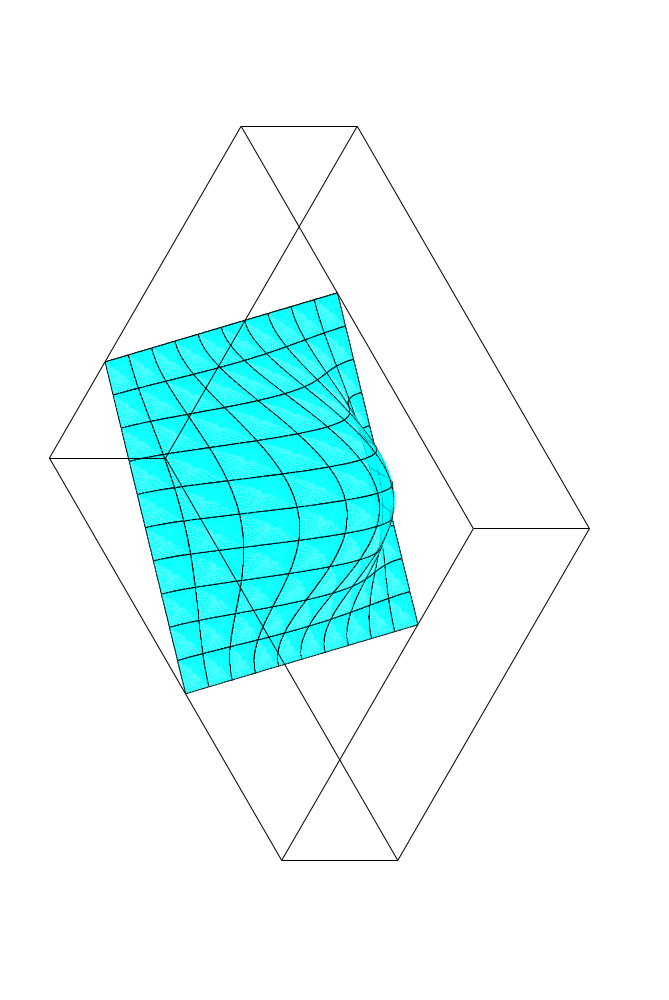}&
		\includegraphics[height=7.4cm, angle=90]{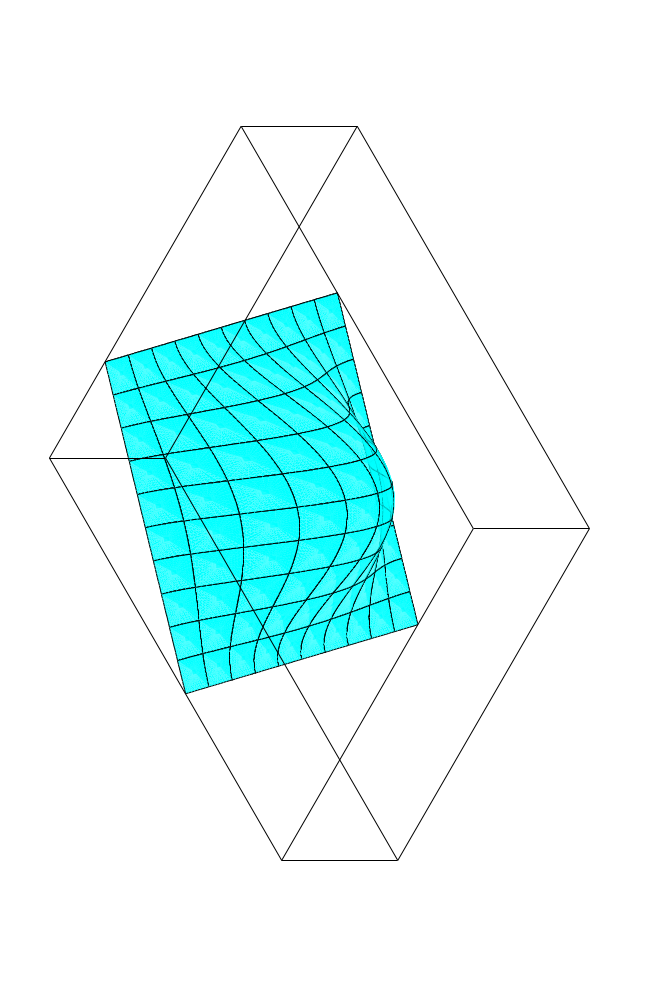} \\
		(b) & (d) \\[-0.25ex] \hline
		\end{tabular}
		\caption{Cube $[-1,+1]^{3}$ with \emph{no-penetration boundary conditions} subject to the numerically determined \emph{rotation-free body force}, \mbox{$\vek{f}-\curl\vek{A}_{h}$} \eqref{eq:closed-form-rotation-free-body-force}, Step2. 
		The left column shows results for $\mu=1$ while results for $\mu=10^{-4}$ are shown in the right column. All results are drawn on the mid plane $z=0$. The solution values in the three-dimensional vector graphs (b) and(d) are scaled with the factor 0.25.\newline 
		The pressure reaction $p_{2}$, $\min{(p_{2})}=-1.524\cdot10^{-1}$ and $\max{(p_{2})}=8.478\cdot10^{-1}$ is the same for both cases $\mu=1$ and $\mu=10^{-4}$;\newline
		(a) iso-level map, pressure reaction $p_{2}$, \\
		(b) three-dimensional vector graph, pressure reaction $p_{2}$,\\ 
		(c) iso-level map, pressure reaction $p_{2}$,\\
		(d) three-dimensional vector graph, pressure reaction $p_{2}$.\\
		\label{tab:rotation-free-body-force-no-penetration}}
    \end{table}

   \begin{table}[tbh]
		\begin{tabular}{|c|c|}\hline
		\includegraphics[height=7cm, angle=90]{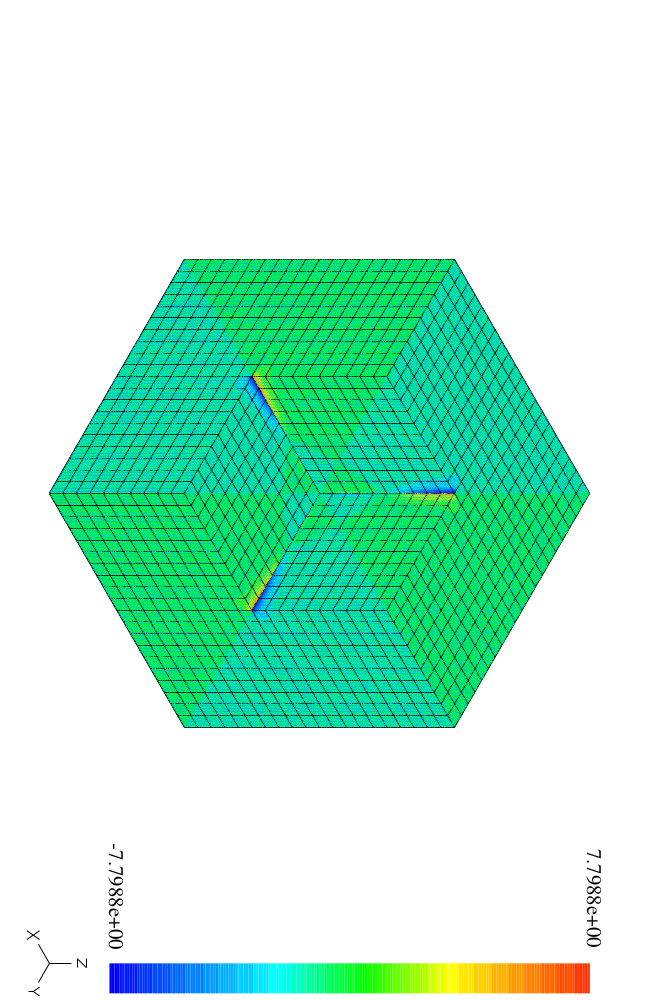} &
	    \includegraphics[height=7cm, angle=90]{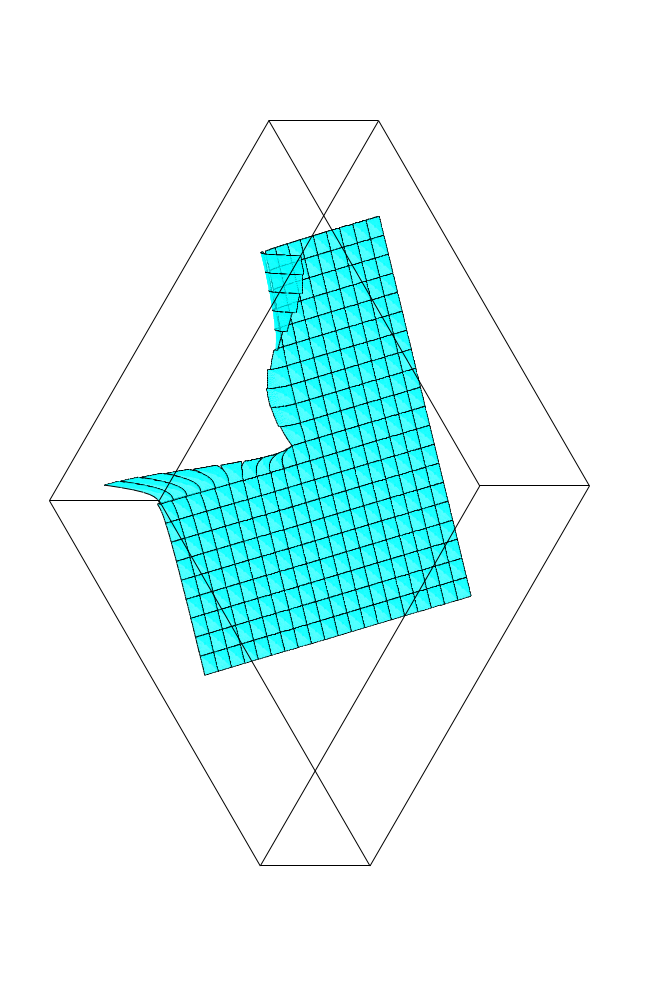}  \\
		(a) & (b) \\[-0.25ex] \hline
		\end{tabular}
		\caption{Three dimensional L-shaped domain $\Omega = [-1,+1]^{3}\setminus[0,+1]^{3}$ with \emph{no-penetration boundary conditions} subject to the
		analytical \emph{divergence-free body force}, \mbox{$\vek{f}_{1}=\curl\vek{A}$} given by \eqref{eq:divergence-free body force L-shaped domain}. 
		The results for $\mu=1$ are shown.\newline 
		(a) iso-level map, pressure reaction $p_{1}$, $\min{(p_{1})}=-7.799$ and $\max{(p_{1})}=7.799$\\
		(b) three-dimensional vector graph, pressure reaction $p_{1}$, drawn on the mid plane $z=0$. The solution values are scaled with the factor 0.25.\\ 
		\label{tab:Lshaped-divergence-free-body-force-no-penetration}}
    \end{table}
    	
    \clearpage
	\appendix
	
	\section{Proof of the elasticity - curl-curl connection}
	\label{sec:proof elasticity curlcurl connection}
	\footnotesize
	
	\begin{lemma}
		\label{lem:prework_curlcurl_elasticity_neumann_bc}
		For $\vek{u},\vek{v}\in \vek{H}^1(\Omega)$ with $\div{\vek{u}}=0$ there holds
		\begin{align}
		\int_{\Gamma}2\ten{e}(\vek{u})\vek{n}\cdot \vek{v}\,d\Gamma = \int_{\Gamma}\curl \vek{u}\times\vek{n}\cdot\vek{v}-\curl \vek{v}\times\vek{n}\cdot\vek{u}+2\ten{e}(\vek{v})\vek{n}\cdot \vek{u}-\frac{4}{3}u_n\div{\vek{v}}\,d\Gamma.
		\end{align}
		\begin{proof}\!\!:
		We will need the following point-wise identities for not necessarily divergence-free functions 
		$\vek{w}\in \vek{C}^2(\Omega)$,
		\begin{align}
		-\bdiv{[2\ten{e}(\vek{w})]}&= \tfrac{2}{3}\,\mathrm{grad\,}\div{\,\vek{w}}-\bdiv{[\grad \vek{w}^{\T}]} -
		\bdiv[\grad\vek{w}]=-\tfrac{1}{3}\,\mathrm{grad\,}\div{\vek{w}}-\bdiv[\grad\vek{w}],\label{eq:diveps_id}\\
		\curl\curl \vek{w} &= \mathrm{grad\,}\div{\,\vek{w}}-\bdiv[\grad\vek{w}] \overset{\eqref{eq:diveps_id}}{=} -\bdiv[2\ten{e}(\vek{w})]+\tfrac{4}{3}\mathrm{grad\,}\div{\,\vek{w}}.\label{eq:curlcurl_id}
		\end{align}
		Since $\vek{C}^2(\Omega)\cap \vek{H}^1(\Omega)$ is dense in $\vek{H}^1(\Omega)$ we will prove the stated identity first for arbitrary\linebreak[2] \mbox{$\vek{w}\in \vek{C}^2(\Omega)\cap \vek{H}^1(\Omega)$}. As $\vek{H}^1(\Omega)\subset \vek{H}(\mathrm{div},\Omega)$ and the existence of commuting interpolation operators \cite{Schoeberl01} we can approximate $\vek{u}$ with a sequence of divergence-free $\vek{w}$. Further, $\vek{v}$ gets approximated by $\vek{z}\in \vek{C}^2(\Omega)\cap \vek{H}^1(\Omega)$. With integration by parts there holds
		\begin{align*}
		\int_{\Gamma}2\ten{e}(\vek{w})\vek{n}\cdot \vek{z}\,ds &= \int_{\Omega}\bdiv[2\ten{e}(\vek{w})]\cdot\vek{z}+2\ten{e}(\vek{w}):\ten{e}(\vek{z})\,d\Omega\\ 
		& \overset{\eqref{eq:curlcurl_id}}{=} \int_{\Omega}2\ten{e}(\vek{w}):\ten{e}(\vek{z})-\curl\curl\vek{w}\cdot\vek{z}+\frac{4}{3}\mathrm{grad\,}\div{\vek{w}}\cdot\vek{z}\,d\Omega\\
		&=\int_{\Omega}2\ten{e}(\vek{w}):\ten{e}(\vek{z})-\curl\vek{w}\cdot\curl\vek{z}\,d\Omega+\int_{\Gamma}\curl\vek{w}\times\vek{n}\cdot\vek{z}\,d\Gamma\\
		&=\int_{\Omega}2\ten{e}(\vek{w}):\ten{e}(\vek{z})-\vek{w}\cdot\curl\curl\vek{z}\,d\Omega+\int_{\Gamma}\curl\vek{w}\times\vek{n}\cdot\vek{z}-\curl\vek{z}\times\vek{n}\cdot\vek{w}\,d\Gamma\\
		&\overset{\eqref{eq:curlcurl_id}}{=}\int_{\Omega}2\ten{e}(\vek{w}):\ten{e}(\vek{z})+\vek{w}\cdot(\bdiv[2\ten{e}(\vek{z})]-\frac{4}{3}\mathrm{grad\,}\div{\vek{z}})\,d\Omega+\int_{\Gamma}\curl\vek{w}\times\vek{n}\cdot\vek{z}-\curl\vek{z}\times\vek{n}\cdot\vek{w}\,d\Gamma\\
		&=\int_{\Omega}\frac{4}{3}\div{\vek{w}}\,\div{\vek{z}}\,d\Omega+\int_{\Gamma}\curl\vek{w}\times\vek{n}\cdot\vek{z}-\curl\vek{z}\times\vek{n}\cdot\vek{w}+2\ten{e}(\vek{z})\vek{n}\cdot\vek{w}-\frac{4}{3}w_n\div{\vek{z}}\,d\Gamma\\
		&=\int_{\Gamma}\curl\vek{w}\times\vek{n}\cdot\vek{z}-\curl\vek{z}\times\vek{n}\cdot\vek{w}+2\ten{e}(\vek{z})\vek{n}\cdot\vek{w}-\frac{4}{3}w_n\div{\vek{z}}\,d\Gamma
		\end{align*}
		concluding the proof by density.
		\end{proof}
	\end{lemma}
	
	Now, we can prove the statement of Lemma~\ref{th:The elasticity - curl-curl connection}:
	\begin{proof}\!\!:
		Since $\vek{C}^2(\Omega)\cap \vek{H}^1(\Omega)$ is dense in $\vek{H}^1(\Omega)$ we will prove the stated identities first for arbitrary\linebreak[2] \mbox{$\vek{w}\in \vek{C}^2(\Omega)\cap \vek{H}^1(\Omega)$}. We start with \eqref{eq:curlcurl_epseps} by integration by parts,
		\begin{equation*}
		\int_{\Omega}\curl \vek{u}\cdot \curl\vek{w}\,d\Omega
		= \int_{\Omega} \vek{u}\cdot \curl\curl\vek{w}\,d\Omega-\int_{\Gamma} \vek{u}\times\vek{n}\cdot\curl\vek{w} \,d\Gamma,
		\end{equation*}
		where the argument of the boundary term  in the local triad $(\vek{t}_1,\vek{t}_2,\vek{n})$ at $\vek{x}\in\Gamma$ reads,
		\begin{equation*}
		\vek{u}\times\vek{n}\cdot\curl\vek{w} =  
		\begin{bmatrix} u_{t_2} & -u_{t_1} & 0 \end{bmatrix} \cdot
		\begin{Bmatrix}
		\dfrac{\partial w_n}{\partial t_2}    -\dfrac{\partial w_{t_2}}{\partial n}\\[1.5ex] 
		\dfrac{\partial w_{t_1}}{\partial n}  -\dfrac{\partial w_{n}}{\partial t_1}\\[1.5ex] 
		\dfrac{\partial w_{t_2}}{\partial t_1}-\dfrac{\partial w_{t_1}}{\partial t_2}
		\end{Bmatrix}.
		\end{equation*}
		Next, we use identity \eqref{eq:curlcurl_id}, integration by parts, and our assumption that $\div{\vek{u}}=0$ to deduce,
		\begin{align*}
		\int_{\Omega}\curl \vek{u}\cdot\curl\vek{w}\,d\Omega &=\int_{\Omega}\vek{u}\cdot\left(-\bdiv[2\ten{e}(\vek{w})]+\tfrac{4}{3}\,\mathrm{grad\,}\div{\,\vek{w}}\right)\,d\Omega-\int_{\Gamma}\vek{u}\times\vek{n}\cdot\curl\vek{w}\,d\Gamma\\
		&=\int_{\Omega}2\ten{e}(\vek{u}):\ten{e}(\vek{w})  - \tfrac{4}{3}\div{\vek{u}}\div{\vek{w}}\,d\Omega
		-\int_{\Gamma}\vek{u}\times\vek{n}\cdot\curl\vek{w}+\vek{u}\cdot2\ten{e}(\vek{w})\vek{n}-\tfrac{4}{3} u_{n}\,\div{\vek{w}}\,d\Gamma\\
		&=\int_{\Omega}2\ten{e}(\vek{u}):\ten{e}(\vek{w})\,d\Omega
		-\int_{\Gamma}\vek{u}\times\vek{n}\cdot\curl\vek{w}+\vek{u}\cdot2\ten{e}(\vek{w})\vek{n}-\tfrac{4}{3} u_{n}\,\div{\vek{w}}\,d\Gamma,
		\end{align*}
		where $u_{n}:=\vek{n}\cdot\vek{u}$. The additional boundary terms, 
		\begin{align*}
		\vek{u}\cdot 2\ten{e}(\vek{w})\vek{n}-\tfrac{4}{3}u_n\,\div{\,\vek{w}} &=  \vek{u}\cdot\left(\frac{\partial\vek{w}}{\partial\vek{n}}+
		\mathrm{grad\,}w_n\right)-2u_n\,\div{\,\vek{w}}\\
		&=
		\begin{bmatrix}
		u_{t_1} & u_{t_2} & u_{n}
		\end{bmatrix}
		\cdot
		\begin{Bmatrix}
		\dfrac{\partial w_{t_1}}{\partial{n}}+\dfrac{\partial w_{n}}{\partial{t_1}}\\[1.5ex]
		\dfrac{\partial w_{t_2}}{\partial{n}}+\dfrac{\partial w_{n}}{\partial{t_2}}\\[1.5ex]
		\dfrac{\partial w_{n}}{\partial{n}}+\dfrac{\partial w_{n}}{\partial{n}}
		\end{Bmatrix} 
		-2u_n\,\left(\dfrac{\partial w_{t_1}}{\partial{t_1}}+\dfrac{\partial w_{t_2}}{\partial{t_2}}+\dfrac{\partial w_{n}}{\partial{n}}\right),
		\end{align*}
		are likewise expanded in the local triad $(\vek{t}_1,\vek{t}_2,\vek{n})$ at $\vek{x}\in\Gamma$ and the sum can be written as,

		\begin{align*}
		\vek{n}\times\vek{u}\cdot\curl\vek{w}+\vek{u}\cdot 2\ten{e}(\vek{w})\vek{n}-\tfrac{4}{3}u_n\,\div{\,\vek{w}} = 2(\vek{u}\cdot \nabla w_n-u_n\div{\vek{w}}).
		\end{align*} 
		If now $\vek{n}\times\vek{u}=\vek{n}\times\vek{w}=\vek{0}$, or $u_n=w_n=0$ the two boundary terms cancel out,
		\begin{align*}
		\vek{n}\times\vek{u}\cdot\curl\vek{w}+\vek{u}\cdot 2\ten{e}(\vek{w})\vek{n}-\tfrac{4}{3}u_n\,\div{\,\vek{w}}=0.
		\end{align*}
		Thus in both cases the boundary terms are zero and we arrive at \eqref{eq:curlcurl_epseps},
		\begin{align*}
		\int_{\Omega}\curl \vek{u}\cdot\curl\vek{w}\,d\Omega &= \int_{\Omega}2\ten{e}(\vek{u}):\ten{e}(\vek{w})\,d\Omega.
		\end{align*}
		In the case of $\curl\vek{u}\times \vek{n}=\vek{0}$ on $\Gamma$ we use Lemma~\ref{lem:prework_curlcurl_elasticity_neumann_bc} for \eqref{eq:curlcurl_epseps_bnd}.
		
		Now, by density, we conclude that \eqref{eq:curlcurl_epseps_bnd} and \eqref{eq:curlcurl_epseps} holds for all $\vek{v}\in \vek{H}^1(\Omega)$.
	\end{proof}
	
	\section{Stability of curl-curl problem}
	\label{sec:stability curlcurl problem}
	\footnotesize
	To prove inf-sup stability of \eqref{eq:auxiliary curl-curl elasticity formulation} we use Brezzi's Theorem for mixed saddle point problems \cite{Bre74} with the bilinear forms
	\begin{align*}
	a:\vek{W}(\Omega)\times \vek{W}(\Omega)\to\mathbb{R},&\quad a(\vek{u},\vek{v})=\int_{\Omega}\mu\,\curl\vek{u}\cdot\curl\vek{v}\,d\Omega,\\
	b:\vek{W}(\Omega)\times \vek{P}(\Omega)\to\mathbb{R},&\quad b(\vek{u},p)=\int_{\Omega}\mathrm{grad\,} p\cdot \vek{u}\,d\Omega,
	\end{align*}
	where we consider the two cases of $\vek{W}(\Omega)=\vek{H}_0(\mathrm{curl},\Omega)$, $P(\Omega)=H^1_0(\Omega)$ and $\vek{W}(\Omega)=\vek{H}(\mathrm{curl},\Omega)$, $P(\Omega)=H^1(\Omega)\backslash\mathbb{R}$, respectively, fulfilling both the compatibility condition $\mathrm{grad\,} P(\Omega) \subset \vek{W}(\Omega)$, and their corresponding Sobolev norms.
	
	Continuity of the bilinear forms and right-hand sides, $\vek{f}\in \vek{L}^2(\Omega)$ and $\vek{s}_{\rm a}\in \vek{L}^2(\Gamma)$, is obvious
	\begin{align*}
	& |a(\vek{u},\vek{v})|\leq \mu\|\vek{u}\|_{\bm{H}(\mathrm{curl},\Omega)}\|\vek{v}\|_{\bm{H}(\mathrm{curl},\Omega)}\quad &&\forall \vek{u},\vek{v}\in \vek{W}(\Omega),\\
	& |b(\vek{u},p)| \leq \|\vek{u}\|_{\bm{H}(\mathrm{curl},\Omega)}\|p\|_{H^1(\Omega)}\quad &&\forall \vek{u}\in \vek{W}(\Omega), p\in P(\Omega),\\
	&\left|\int_{\Omega}\vek{f}\cdot\vek{u}\,d\Omega\right|\leq \|\vek{f}\|_{\bm{L}^2(\Omega)}\|\vek{u}\|_{\bm{H}(\mathrm{curl},\Omega)}\quad &&\forall \vek{u}\in \vek{W}(\Omega),\\
	& \left|\int_{\Gamma}\vek{n}\times \vek{u}\cdot\vek{s}_{\rm a}\,d\Gamma\right|\leq \|\vek{s}_{\rm a}\|_{\bm{L}^2(\Gamma)}\|\vek{u}\|_{\bm{H}(\mathrm{curl},\Omega)}\quad &&\forall \vek{u}\in \vek{W}(\Omega).
	\end{align*}
	The coercivity of the bilinear form $a(\cdot,\cdot)$ on the kernel, with $\alpha_1>0$,
	\begin{align}
	&a(\vek{u},\vek{u})\geq\alpha_1\|\vek{u}\|^2_{\bm{H}(\mathrm{curl},\Omega)}\quad \forall \vek{u}\in \vek{W}_0(\Omega),\\
	&\vek{W}_0(\Omega)=\{\vek{u}\in \vek{W}(\Omega)\,|\, b(\vek{u},p)=0\quad \forall p\in P(\Omega)  \}
	\end{align}
	follows by the Friedrich's-type inequality, $c_F>0$,
	\begin{align}
	\|\vek{u}\|_{\bm{L}^2(\Omega)} \leq c_F\|\curl \vek{u}\|_{\bm{L}^2(\Omega)}\qquad \forall \vek{u}\in \vek{W}_0(\Omega).
	\end{align}
	For the LBB condition to hold, $\beta_1>0$,
	\begin{align}
	\sup_{\bm{u}\in \bm{W}(\Omega)}\frac{b(\vek{u},p)}{\|\vek{u}\|_{\bm{H}(\mathrm{curl},\Omega)}}\geq\beta_1\|p\|_{H^1(\Omega)}\quad\forall p\in P(\Omega)
	\end{align}
	we choose for given $p\in Q(\Omega)$ $\vek{u}=\mathrm{grad\,} p \in \vek{W}(\Omega)$ such that
	\begin{align*}
	\sup_{\bm{u}\in \bm{W}(\Omega)}\frac{b(\vek{u},p)}{\|\vek{u}\|_{\bm{H}(\mathrm{curl},\Omega)}} \geq \frac{\|\mathrm{grad\,} p\|^2_{L^2(\Omega)}}{\|\mathrm{grad\,} p\|_{L^2(\Omega)}} = \|\mathrm{grad\,} p\|_{L^2(\Omega)}.
	\end{align*}
	Thus, by Brezzi, the curl-curl problem is uniquely solvable. The solution fulfills the stability estimate
	\begin{align}
	\|\vek{u}\|_{\bm{H}(\mathrm{curl},\Omega)}+\|p\|_{H^1(\Omega)}\leq c\left(\|\vek{f}\|_{\bm{L}^2(\Omega)}+\|\vek{s}_{\rm a}\|_{\bm{L}^2(\Gamma)}\right).
	\end{align}

	\begin{remark}
		In the discrete case one has to choose the polynomial order for $p$ one higher than for $\vek{u}$.
	\end{remark}
	
	\section{The assumed body force on the cube}
	\label{sec:The assumed body force on the cube}
	\footnotesize
    \subsection{The divergence-free part}
    \footnotesize
	 The divergence-free part of the body force is derived from the divergence-free vector potential \mbox{$\vek{A}=(1-x^{2})(1-y^{2})\vek{e}_{z}$} as $\vek{f}_{1}=\curl\vek{A}$, yielding, 
	\begin{gather}
	\vek{f}_{1}=-2(1-x^{2})y\,\vek{e}_{x}+ 2(1-y^{2})x\,\vek{e}_{y}, \nonumber \\[-0.5ex]
	\label{eq:closed-form-divergence-free-body-force} \\[-0.5ex]
	\div{\,\vek{f}_{1}}=4xy-4yx=0, \qquad 
	\curl\vek{f}_{1}=\left[2(1-y^{2})+2(1-x^{2})\right]\vek{e}_{z}. \nonumber
	\end{gather}
	It is readily checked that $\vek{f}_{1}$ is parallel to the six faces of the cube, i.e. 
	\mbox{$\vek{n}\cdot\vek{f}_{1}=0$} for $\vek{x}\in\Gamma$.

	\subsection{The rotation-free part} 
	\footnotesize
	The rotation-free part of the body force  is derived from the potential $\phi=(1-x^{2})^{2}(1-y^{2})^{2}(1-z^{2})^{2}$.  It is determined as,
	\begin{gather}
	\vek{f}_{2}=\mathrm{grad\,}\phi=\phi_{,x}\,\vek{e}_{x}+\phi_{,y}\,\vek{e}_{y}+\phi_{,z}\,\vek{e}_{z}, \quad
	\begin{cases}
	\phi_{,x} & = -4x(1-x^{2})(1-y^{2})^{2}(1-z^{2})^{2},\\
	\phi_{,y} & = -4y(1-y^{2})(1-x^{2})^{2}(1-z^{2})^{2},\\
	\phi_{,z} & = -4z(1-z^{2})(1-x^{2})^{2}(1-y^{2})^{2},
	\end{cases} \nonumber \\
	\curl\vek{f}_{2}=(\phi_{,zy}-\phi_{,yz})\vek{e}_{x}+(\phi_{,xz}-\phi_{,zx})\vek{e}_{y}+(\phi_{,yx}-\phi_{,xy})\vek{e}_{z} = \vek{0}, \label{eq:closed-form-rotation-free-body-force}\\[1ex]
	\div{\vek{f}_{2}}= \phi_{,xx}+\phi_{,yy}+\phi_{,zz} \neq 0 \quad\text{for}\quad \vek{x}\notin\Gamma\quad 
	\begin{cases}
	\phi_{,xx} &= 4(3x^{2}-1)(1-y^{2})^{2}(1-z^{2})^{2},\\
	\phi_{,yy} &= 4(3y^{2}-1)(1-x^{2})^{2}(1-z^{2})^{2},\\
	\phi_{,zz} &= 4(3z^{2}-1)(1-x^{2})^{2}(1-y^{2})^{2},
	\end{cases}
	\nonumber 
	\end{gather}
	By construction, the potential vanishes at the faces of the cube, i.e., $\phi(\vek{x})=0$ for $\vek{x}\in\Gamma$. 
	This fact implies that the pressure reaction it induces vanishes on the boundary of the cube domain.
	Further, $\vek{f}_{2}$ is not divergence-free in 
	$\Omega$, i.e. it is not a harmonic function.  Finally, it is readily shown that $\vek{f}_{2}$ is
	perpendicular and parallel to the faces of the cube, i.e., $\vek{n}\times\vek{f}_{2}=\vek{0}$ and   
	$\vek{n}\cdot\vek{f}_{2}=0$ for $\vek{x}\in\Gamma$, respectively.

\end{document}